# ROUGH EVOLUTION EQUATIONS

By Massimiliano Gubinelli and Samy Tindel

*Université de Paris-Dauphine and Université de Nancy*

We generalize Lyons' rough paths theory in order to give a pathwise meaning to some nonlinear infinite-dimensional evolution equation associated to an analytic semigroup and driven by an irregular noise. As an illustration, we discuss a class of linear and nonlinear 1d SPDEs driven by a space–time Gaussian noise with singular space covariance and Brownian time dependence.

CONTENTS

**1. Introduction.** This paper can be seen as part of an ongoing project whose aim is to give a pathwise definition to stochastic PDEs. Indeed, the rough path theory [5, 13, 16, 17] and its variants [4, 6] have now reached a certain level of maturity, leading to a proper definition of differential equations driven by irregular signals and in particular by a fractional Brownian motion [2]. Starting from this observation, we have tried in [12] to define and solve the following general problem: let $\mathcal{B}$ be a separable Banach space, and $A: D(A) \to \mathcal{B}$ the infinitesimal generator of an analytical semigroup $\{S_t; t \geq 0\}$ on $\mathcal{B}$, inducing the family $\{\mathcal{B}_\alpha; \alpha \in \mathbb{R}\}$ with $\mathcal{B}_\alpha = D((-A)^\alpha)$. Let also $f$ be a function from $\mathcal{B}$ to $\mathcal{L}(\mathcal{B}_{-\alpha}, \mathcal{B}_{-\alpha})$ for a given $\alpha > 0$ and $x$ a noisy input, considered as a function from $\mathbb{R}_+$ to $\mathcal{B}_{-\alpha}$. Then, for $T > 0$, consider the equation

$$dy_t = Ay_t\, dt + f(y_t)\, dx_t, \qquad t \in [0, T], \tag{1}$$

with an initial condition $y_0 \in \mathcal{B}$. The main example we have in mind is the case of the 1-dimensional heat equation in $[0, 1]$, namely $\mathcal{B} = L^2([0, 1])$, $A = \Delta$ with Dirichlet boundary conditions, the usual Sobolev spaces $\mathcal{B}_\alpha = H_\alpha = W^{2\alpha,2}$, and $x$ a fractional Brownian motion with Hurst parameter $H$ taking values in $\mathcal{B}_{-\alpha}$. Notice, in particular, that we wish to consider a noise $x$ which is irregular in both time and space. Then, in [12], we gave a local existence and uniqueness result for equation (1), by considering it in its mild form

$$y_t = S_t y_0 + \int_0^t S_{ts} f(y_s)\, dx_s, \tag{2}$$

where we let $S_{ts} = S_{t-s}$ and interpreting the integral in this mild formulation as a Young integral. Once the equation is set under the form (2), the main problem one is faced with is to quantify the regularization of the semigroup $S_{ts}$ on the term $f(y_s)\, dx_s$, and then to elaborate the right fixed point argument in order to solve the equation. The general results of [12] could be applied in the case of the stochastic heat equation driven by a fractional Brownian motion with Hurst parameter $H > 1/2$. They should be compared with the reference [18], where a nonlinear fractional SPDE is solved thanks to some fractional calculus methods, but where $x$ is a smooth noise in space.

In the current article, we would like to go one step further with respect to [12], and set the basis of a real rough path expansion in order to define and solve equation (2), which would allow to consider, in the case of the heat equation in $[0,1]$, a fractional Brownian motion with Hurst parameter $H \le 1/2$. This task is quite long and involved, but let us summarize at this point some of the ideas we have followed.

(1) We will recast equation (2) in a suitable way for expansions according to the following simple observation: we have tried to solve our evolution equation by means of its infinite-dimensional setting, since it allows to consider $x$ and $y$ as functions of a unique parameter $t \in [0,T]$, which makes its rough path type analysis easier (see [11] and [22] for a multiparametric setting). However, when we come to the applications to the heat equation, we will consider the evolution equation in $[0,T] \times [0,1]$ under the form

$$y(t,\xi) = \int_0^1 G_t(\xi,\eta) y_0(\eta)\, d\eta + \int_0^t \int_0^1 G_{t-s}(\xi,\eta)\sigma(y_s(\eta)) x(ds,d\eta), \tag{3}$$

where $G$ stands for the fundamental solution to the heat equation, $\sigma : \mathbb{R} \to \mathbb{R}$ is a regular function, and $x(ds,d\eta)$ is understood as the distributional derivative of a real-valued continuous process on $[0,T] \times [0,1]$. This definition of our equation is of course equivalent to (2) when $f$ is considered as the pointwise nonlinear operator $[f(y_t)](\xi) \equiv \sigma(y_t(\xi))$. Now, when written under its multiparametric form (3), the equation is also equivalent to

$$y(t,\xi) = \int_0^1 G_t(\xi,\eta) y_0(\eta)\, d\eta + \int_0^t \int_0^1 G_{t-s}(\xi,\eta) x(ds,d\eta) \sigma(y_s(\eta)),$$

and it happens that this simple reformulation is much more convenient for our future expansions than the original one. When we go back to the original infinite-dimensional setting, we can recast (2) into

$$y_t = S_t y_0 + \int_0^t S_{ts}\, dx_s f(y_s), \tag{4}$$

where $f$ is now a smooth function from $\mathcal{B}$ to $\mathcal{B}$, and $x$ will be understood as a Hölder-continuous process taking values in a space of *deregularizing* operators from $\mathcal{B}$ to a distributional space $\mathcal{B}_{-\zeta}$ for a certain $\zeta > 0$. The product $dx_s f(y_s)$ will then be regularized again by the action of $S_{ts}$, in a way which will be quantified later on. Notice that the form (4) of our evolution equation is a little unusual in the SPDE theory, but makes sense in our context.

(2) Instead of considering Riemann sums like in [12] or like in the original Lyons' theory [16], our analysis will be based on the theory of *generalized differentials*, called $k$-increments, contained in [6]. Roughly speaking, this theory is based on the fact that an elementary operator, called $\delta$, can transform an integral $\int_s^t dg_u[h_u - h_s]$, seen as a function of the variables $s$ and $t$, into a finite difference product $(g_t - g_s)(h_t - h_s)$. Furthermore, under some additional regularity properties on $g$ and $h$, the operator $\delta$ can be inverted, and its inverse $\Lambda$, called *sewing map* (from [4]), will be the building stone of our extension of the notion of integral. Notice that, whenever $g$ and $h$ are Hölder-continuous with Hölder exponent $> 1/2$, this extension coincides with the usual Young integral. When we consider an integral of the form $\int_s^t dg_u\,\phi(g_u)$ for a Hölder-continuous function $g$ with Hölder exponent in $(1/3, 1/2]$ admitting a Lévy area, our definition of integral also coincides with Lyons' one, as shown in [6]. In fact, if the usual rough path theory gives a richer point of view on the algebraic structure of the path $x$, it is worth mentioning that our approach has at least two advantages:

(1) Once our unusual setting is assimilated, it becomes quite easy to figure out how a given expansion in terms of $x$ can be leaded. And indeed, it will become clear throughout the paper, that the $k$-increments theory provides a tool allowing some natural computations for our generalized integrals.
(2) The only step where a discretization procedure is needed is the construction of the $\Lambda$ map alluded to above, and this avoids some of the cumbersome calculations which are one of the main ingredients of the rough path theory.

We hope that this paper will advocate for the use of the $k$-increments theory, which obviously does not exclude the other approaches [4, 16].

(3) The fact that we are dealing with an evolution problem will force us to change some of the algebraic structure we will rely on, especially if one wants to take advantage of the regularizing effect of $S_t$. This will lead us to introduce an operator $a_{ts} = S_{ts} - \mathrm{Id}$ for $t \geq s$, and a modified $\delta$ operator, called $\hat{\delta}$, defined by $\hat{\delta} = \delta - a$. The whole increment theory will have to be built again based on this modified operator, and we will see that it is really suitable for the evolution setting induced by (4). In particular, we will be able to define analogs of the Lévy area and of the higher-order iterated integrals, which are of course harder to express than in the finite-dimensional case, but can be written, in the bilinear case [that is $\sigma(r) = r$ in (3)], as

$$X_{ts}^2 = \int_s^t S_{tu}\,dx_u \int_s^u S_{uv}\,dx_v\,S_{vs}, \qquad X_{ts}^3 = \int_s^t S_{tu}\,dx_u\,X_{us}^2, \qquad \text{etc.} \tag{5}$$

Obviously, a convenient definition of iterated integrals is the key to reach the case of a Hölder-continuous noise of order $\leq 1/2$.

(4) The whole integration theory can be expressed in an abstract way, by just supposing a certain set of assumptions on some incremental operators like $X^2$ and $X^3$. However, we will try to check these assumptions in some interesting cases, like the infinite-dimensional fractional Brownian motion for our Young type integration, or the infinite-dimensional Brownian motion for our step 3 expansion, based on $X^2$ and $X^3$. Notice that the rough expansions for the fractional Brownian motion should be investigated in details too, but one is faced with an additional problem in this situation: on one hand, a Stratonovich type integration requires a lot of regularity in space for the noise, due to the well-known presence of some trace terms. On the other hand, the Skorokhod integral does not fulfill the algebraic requirements we ask for our integral extension. A discussion of these problems and some ideas to solve them will be included at the end of the paper, but for sake of conciseness, we will postpone a complete development of this part to a subsequent paper, and stick here to the Brownian case.

This paper is structured as follows. In Section 2, we recall the basic setup of [6] which allows to embed the theory of rough paths in a theory of integration of *generalized differentials*, called here *k-increments*. We wrote it with the aim of having a self-contained and pedagogical introduction to the topic. However, we give also a new and very elementary proof of the existence of the basic integration map $\Lambda$ of [6]. In Section 3, we introduce and study a modified coboundary induced by the operator $\hat{\delta}$ on the complex of increments, using the additional data provided by an analytic semigroup $S$, in such a way that the new complex can be shown to act simply on convolution integrals of the form appearing in equation (4) and on their iterated versions. This new complex maintain many of the properties of the original complex (e.g., its cohomology is trivial) and it is shown that when equipped with Hölder-like norms which measures "smallness" of the increments, it admits a map, called $\hat{\Lambda}$ here, which is the main tool for building an integration (or better, convolution) theory over those 1-increments which are good enough (again, in a suitable sense, to be specified in due time). A key feature of this perturbed complex is that, due to the convolution with the semigroup $S$, "space" and "time" regularity of increments depends on each other: we can gain space regularity by loosing some time regularity and vice-versa. This property will be essential for the solution of the evolution problem by fixed-point arguments. In Section 4.2, we use the theory outlined in Section 3.2 to define the convolution integral in the Young sense and solve a class of non-linear evolution problems, reobtaining some results of the work [12]. Notice that we will also improve some of our previous results contained in [12], in the sense that we will be able to construct global solutions to our evolution equations in the Young context. In Section 5, we study the bilinear evolution problem

$$y_t = S_t y_0 + \int_0^t S_{ts}\, dx_s\, y_s. \tag{6}$$

We will also introduce a notion of rough-path suitable for noises driving evolution equations. By exploiting this pathwise technique, we are able to obtain automatically the flow semigroup of the equation and we will show how to express this semigroup as a convergent series of iterated-integrals which are the lift of the step-3 rough path used in the construction of the solution. In Section 6, we turn to a nonlinear case of evolution system, namely the case of the quadratic type equation

$$y_t = S_t y_0 + \int_0^t S_{ts}\, dx_s B(y_s \otimes y_s),$$

where $B$ stands for the pointwise multiplication of functions. This requires the additional careful introduction of a collection of a priori increments indexed by planar trees, and an associate notion of controlled path. Finally, all our results will be applied in the concrete case of the stochastic heat equation on the circle, in a setting recalled at Section 3.4. The case of a fractional Brownian case is handled in the special situation of the Young theory, while we stick to the example of an infinite-dimensional Brownian motion in the rougher situation. We build the rough path associated to this latter noise and provide concrete conditions where the theory outlined in the previous sections can be fruitfully applied. A systematic study of the regularity properties of the incremental operators defined as $X^2$ or $X^3$ in (5) will also be provided at Sections 6.5 and 6.6, thanks to some Feynman diagram techniques.

**2. Algebraic integration in one dimension.** The integration theory introduced in [6] is based on an algebraic structure, which turns out to be useful for computational purposes, but has also its own interest. Since this setting is quite nonstandard, compared with the one developed in [16], and since we will elaborate on it throughout the paper, we will recall briefly here its main features. We also provide an elementary proof of the existence of the $\Lambda$ map.

2.1. *Increments.* As mentioned in the Introduction, the extended integral we deal with is based on the notion of increment, together with an elementary operator $\delta$ acting on them. However, this simple structure gives rise to a nice topological structure that we will describe briefly here: first of all, for an arbitrary real number $T > 0$, a vector space $V$, and an integer $k \geq 1$, we denote by $\mathcal{C}_k(V)$ the set of functions $g : [0,T]^k \to V$ such that $g_{t_1 \cdots t_k} = 0$ whenever $t_i = t_{i+1}$ for some $i \leq k-1$. Such a function will be called a $(k-1)$-*increment*, and we will set $\mathcal{C}_*(V) = \bigcup_{k \geq 1} \mathcal{C}_k(V)$. The operator $\delta$ alluded to above can be seen as a coboundary operator acting on $k$-increments, inducing a cochain complex $(\mathcal{C}_*, \delta)$, and is defined as follows on $\mathcal{C}_k(V)$:

$$\delta : \mathcal{C}_k(V) \to \mathcal{C}_{k+1}(V), \qquad (\delta g)_{t_1 \cdots t_{k+1}} = \sum_{i=1}^{k+1} (-1)^i g_{t_1 \cdots \hat{t}_i \cdots t_{k+1}}, \tag{7}$$

where $\hat{t}_i$ means that this particular argument is omitted. Then a fundamental property of $\delta$, which is easily verified, is that $\delta\delta = 0$, where $\delta\delta$ is considered as an operator from $\mathcal{C}_k(V)$ to $\mathcal{C}_{k+2}(V)$. We will denote $\mathcal{ZC}_k(V) = \mathcal{C}_k(V) \cap \operatorname{Ker}\delta|_{\mathcal{C}_k(V)}$ and $\mathcal{BC}_k(V) := \mathcal{C}_k(V) \cap \operatorname{Im}\delta|_{\mathcal{C}_{k-1}(V)}$, respectively the spaces of $k$-*cocycles* and of $k$-*coboundaries*, following standard conventions of homological algebra.

Some simple examples of actions of $\delta$, which will be the ones we will really use throughout the paper, are obtained by letting $g \in \mathcal{C}_1$ and $h \in \mathcal{C}_2$. Then, for any $t, u, s \in [0, T]$, we have

$$(\delta g)_{ts} = g_t - g_s \quad \text{and} \quad (\delta h)_{tus} = h_{ts} - h_{tu} - h_{us}. \tag{8}$$

Furthermore, it is readily checked that the complex $(\mathcal{C}_*, \delta)$ is *acyclic*, that is, $\mathcal{ZC}_{k+1}(V) = \mathcal{BC}_k(V)$ for any $k \geq 1$, or otherwise stated, the sequence

$$0 \to \mathbb{R} \to \mathcal{C}_1(V) \xrightarrow{\delta} \mathcal{C}_2(V) \xrightarrow{\delta} \mathcal{C}_3(V) \xrightarrow{\delta} \mathcal{C}_4(V) \to \cdots \tag{9}$$

is exact. In particular, the following basic property, which we label for further use, holds true.

LEMMA 2.1. *Let $k \geq 1$ and $h \in \mathcal{ZC}_{k+1}(V)$. Then there exists a (nonunique) $f \in \mathcal{C}_k(V)$ such that $h = \delta f$.*

PROOF. This elementary proof is included in [6]; see also Proposition 3.1 below. Let us just mention that $f_{t_1\cdots t_k} = h_{t_1\cdots t_k 0}$ is a possible choice. □

REMARK 2.2. Observe that Lemma 2.1 implies that all the elements $h \in \mathcal{C}_2(V)$, such that $\delta h = 0$, can be written as $h = \delta f$ for some (nonunique) $f \in \mathcal{C}_1(V)$. Thus, we get a heuristic interpretation of $\delta|_{\mathcal{C}_2(V)}$: it measures how much a given 1-increment is far from being an *exact* increment of a function (i.e., a finite difference).

Notice that our future discussions will mainly rely on $k$-increments with $k \leq 2$, for which we will use some analytical assumptions. Namely, we measure the size of these increments by Hölder norms defined in the following way: for $f \in \mathcal{C}_2(V)$ let

$$\|f\|_\mu \equiv \sup_{s,t\in[0,T]} \frac{|f_{ts}|}{|t-s|^\mu} \quad \text{and} \quad \mathcal{C}_2^\mu(V) = \{f \in \mathcal{C}_2(V); \|f\|_\mu < \infty\}.$$

In the same way, for $h \in \mathcal{C}_3(V)$, set

$$\begin{aligned} \|h\|_{\gamma,\rho} &= \sup_{s,u,t\in[0,T]} \frac{|h_{tus}|}{|u-s|^\gamma |t-u|^\rho}, \\ \|h\|_\mu &\equiv \inf\left\{\sum_i \|h_i\|_{\rho_i,\mu-\rho_i}; h = \sum_i h_i, 0 < \rho_i < \mu\right\}, \end{aligned} \tag{10}$$

where the last infimum is taken over all sequences $\{h_i \in \mathcal{C}_3(V)\}$ such that $h = \sum_i h_i$ and for all choices of the numbers $\rho_i \in (0,\mu)$. Then $\|\cdot\|_\mu$ is easily seen to be a norm on $\mathcal{C}_3(V)$, and we set

$$\mathcal{C}_3^\mu(V) := \{h \in \mathcal{C}_3(V); \|h\|_\mu < \infty\}.$$

Eventually, let $\mathcal{C}_3^{1+}(V) = \bigcup_{\mu>1} \mathcal{C}_3^\mu(V)$, and remark that the same kind of norms can be considered on the spaces $\mathcal{ZC}_3(V)$, leading to the definition of some spaces $\mathcal{ZC}_3^\mu(V)$ and $\mathcal{ZC}_3^{1+}(V)$.

With this notation in mind, the following proposition is a basic result which is at the core of our approach to pathwise integration.

PROPOSITION 2.3 (The sewing map $\Lambda$). *There exists a unique linear map $\Lambda : \mathcal{ZC}_3^{1+}(V) \to \mathcal{C}_2^{1+}(V)$ (the sewing map) such that*

$$\delta\Lambda = \mathrm{Id}_{\mathcal{ZC}_3(V)}.$$

*Furthermore, for any $\mu > 1$, this map is continuous from $\mathcal{ZC}_3^\mu(V)$ to $\mathcal{C}_2^\mu(V)$ and we have*

$$\|\Lambda h\|_\mu \le \frac{1}{2^\mu - 2}\|h\|_\mu, \qquad h \in \mathcal{ZC}_3^{1+}(V). \tag{11}$$

PROOF. For sake of completeness, we include a proof of this result here, which is more elementary than the one provided in [6], and which will be generalized at Theorem 3.5. For the sake of notation, we will omit the dependence in $V$ in our functional spaces, and write for instance $\mathcal{C}_3$ instead of $\mathcal{C}_3(V)$. Let then $h$ be an element of $\mathcal{ZC}_3^\mu \subset \mathcal{ZC}_3^{1+}$ for some $\mu > 1$.

*Step 1*: Let us first prove the uniqueness of the 1-increment $M \in \mathcal{C}_2^\mu$ such that $\delta M = h$. Indeed, let $M, \hat{M}$ be two elements of $\mathcal{C}_2^\mu$ satisfying $\delta M = \delta\hat{M} = h$ and set $Q = M - \hat{M}$. Then $\delta Q = 0$ and $Q \in \mathcal{C}_2^\mu$. Invoking Lemma 2.1, there exists an element $q \in \mathcal{C}_1$ such that $Q = \delta q$, but since $\mu > 1$, $q$ is a function on $[0,T]$ with zero derivative, that is a constant and then $Q = 0$.

*Step 2*: Let us construct now a process $M \in \mathcal{C}_2^\mu$, with $\mu > 1$, satisfying $\delta M = h$. Since $\delta h = 0$, invoking again Lemma 2.1, we know that there exists a $B \in \mathcal{C}_2$ such that $\delta B = h$. Pick $s, t \in [0,T]$, such that $s < t$ in order to fix ideas, and for $n \ge 0$, consider the dyadic partition $\{r_i^n; i \le 2^n\}$ of the interval $[s,t]$, where

$$r_i^n = s + \frac{(t-s)i}{2^n} \qquad \text{for } i \le 2^n. \tag{12}$$

Then for $n \ge 0$ set

$$M_{ts}^n = B_{ts} - \sum_{i=0}^{2^n-1} B_{r_{i+1}^n, r_i^n}. \tag{13}$$

Then it is readily checked that $M_{ts}^0 = 0$. Furthermore, we have

$$M_{ts}^{n+1} - M_{ts}^n = \sum_{i=0}^{2^n-1} (B_{r_{2i+2}^{n+1}, r_{2i}^{n+1}} - B_{r_{2i+1}^{n+1}, r_{2i}^{n+1}} - B_{r_{2i+2}^{n+1}, r_{2i+1}^{n+1}})$$
$$= \sum_{i=0}^{2^n-1} (\delta B)_{r_{2i+2}^{n+1}, r_{2i+1}^{n+1}, r_{2i}^{n+1}} = \sum_{i=0}^{2^n-1} h_{r_{2i+2}^{n+1}, r_{2i+1}^{n+1}, r_{2i}^{n+1}},$$

and since $h \in \mathcal{C}_3^\mu$ with $\mu > 1$, we obtain

$$|M_{ts}^n - M_{ts}^{n+1}| \le \frac{\|h\|_\mu (t-s)^\mu}{2^{n(\mu-1)}},$$

which yields that $M_{ts} \equiv \lim_{n\to\infty} M_{ts}^n$ exists, and satisfies inequality (11).

*Step 3*: Let us consider now a general sequence $\{\pi_n; n \ge 1\}$ of partitions $\{r_0^n, r_1^n, \ldots, r_{k_n}^n, r_{k_n+1}^n\}$ of $[s,t]$, with $s = r_0^n < r_1^n < \cdots < r_{k_n}^n < r_{k_n+1}^n = t$. We assume that $\pi_n \subset \pi_{n+1}$, and $\lim_{n\to\infty} k_n = \infty$. Set

$$M_{ts}^{\pi_n} = B_{ts} - \sum_{l=0}^{k_n} B_{r_{l+1}^n, r_l^n}. \tag{14}$$

It is easily seen that there exists $1 \le l \le k_n$ such that

$$|r_{l+1}^n - r_{l-1}^n| \le \frac{2|t-s|}{k_n}. \tag{15}$$

Pick now such an index $l$, and let us transform $\pi^n$ into $\hat{\pi}$, where

$$\hat{\pi} = \{r_0^n, r_1^n, \ldots, r_{l-1}^n, r_{l+1}^n, \ldots, r_{k_n}^n, r_{k_n+1}^n\}.$$

Then, as in the previous step,

$$M_{ts}^{\hat{\pi}} = M_{ts}^{\pi_n} - (\delta B)_{r_{l+1}^n, r_l^n, r_{l-1}^n} = M_{st}^{\pi_n} - h_{r_{l+1}^n, r_l^n, r_{l-1}^n},$$

using the definition of the space $\mathcal{C}_3^\mu$ and the bound (15) we have

$$|M_{ts}^{\hat{\pi}} - M_{ts}^{\pi_n}| \le 2^\mu \|h\|_\mu \left(\frac{t-s}{k_n}\right)^\mu.$$

Repeating now this operation until we end up with the trivial partition $\hat{\pi}_0 \equiv \{s,t\}$, for which $M_{st}^{\hat{\pi}_0} = 0$, we obtain

$$|M_{ts}^{\pi_n}| \le 2^\mu \|h\|_\mu |t-s|^\mu \sum_{j=1}^{k_n} j^{-\mu} \le 2^\mu \|h\|_\mu |t-s|^\mu \sum_{j=1}^{\infty} j^{-\mu} \equiv c_{\mu,h} |t-s|^\mu.$$

Hence, there exists a subsequence $\{\pi_m; m \ge 1\}$ of $\{\pi_n; n \ge 1\}$ such that $M_{ts}^{\pi_m}$ converges to an element $M_{ts}$, satisfying $M_{ts} \le c_{\mu,h} |t-s|^\mu$. With the same considerations as in [13], it can also be checked that the limit $M$ does not depend on the particular sequence of partitions we have chosen, and thus coincides with the one constructed at Step 2.

*Step 4*: It remains to show that $\delta M = h$. Consider then $0 \le s < u < t \le T$, and two sequences of partitions $\pi^n_{us}$ and $\pi^n_{tu}$ of $[s,u]$ and $[u,t]$, respectively, whose meshes tend to 0 as $n \to \infty$. Set also $\pi^n_{ts} = \pi^n_{tu} \cup \pi^n_{us}$. From the previous step, one can construct easily some subsequences $\pi^m_{tu}, \pi^m_{us}, \pi^m_{ts}$, with $\pi^m_{ts} = \pi^m_{tu} \cup \pi^m_{us}$, such that

$$\lim_{m\to\infty} M^{\pi^m_{tu}}_{tu} = M_{tu}, \qquad \lim_{m\to\infty} M^{\pi^m_{us}}_{us} = M_{us}, \qquad \lim_{m\to\infty} M^{\pi^m_{ts}}_{ts} = M_{ts}.$$

Call now $k^m_{ts}$ (resp. $k^m_{tu}, k^m_{us}$) the number of points of the partition $\pi^m_{ts}$ (resp. $\pi^m_{tu}, \pi^m_{us}$). Then a direct computation, using definition (14), shows that for any $0 \le i \le 2^n$ we have

$$\begin{aligned} & M^{\pi^m_{ts}}_{ts} - M^{\pi^m_{su}}_{tu} - M^{\pi^m_{ut}}_{us} \\ &\quad = (\delta B)_{tus} - \left( \sum_{l=0}^{k^m_{ts}+k^m_{us}+1} B_{r^m_{l+1} r^m_l} - \sum_{l=0}^{k^m_{tu}} B_{r^m_{l+1} r^m_l} - \sum_{l=k^m_{tu}+1}^{k^m_{tu}+k^m_{us}+1} B_{r^m_l r^m_{l+1}} \right) \\ &\quad = (\delta B)_{tus} = h_{tus}. \end{aligned}$$

Taking the limit $m \to \infty$ in the latter relation, we get $(\delta M)_{tus} = h_{tus}$, which ends the proof. □

We can now give an algorithm for a canonical decomposition of the preimage of the space $\mathcal{Z}\mathcal{C}^{1+}_3(V)$, or in other words, of a function $g \in \mathcal{C}_2(V)$ whose increment $\delta g$ is smooth enough.

Corollary 2.4. *Take an element $g \in \mathcal{C}_2(V)$, such that $\delta g \in \mathcal{C}^\mu_3(V)$ for $\mu > 1$. Then $g$ can be decomposed in a unique way as*

$$g = \delta f + \Lambda \delta g,$$

*where $f \in \mathcal{C}_1(V)$.*

Proof. Elementary; see [6]. □

At this point, the connection of the structure we introduced with the problem of integration of irregular functions can be still quite obscure to the noninitiated reader. However, something interesting is already going on and the previous corollary has a very nice consequence which is the subject of the following property.

COROLLARY 2.5 (Integration of small increments). *For any 1-increment $g \in \mathcal{C}_2(V)$, such that $\delta g \in \mathcal{C}_3^{1+}$, set $\delta f = (\mathrm{Id} - \Lambda\delta)g$. Then*

$$(\delta f)_{ts} = \lim_{|\Pi_{ts}| \to 0} \sum_{i=0}^{n} g_{t_{i+1} t_i},$$

*where the limit is over any partition $\Pi_{ts} = \{t_0 = t, \ldots, t_n = s\}$ of $[t, s]$ whose mesh tends to zero. The 1-increment $\delta f$ is the indefinite integral of the 1-increment $g$.*

PROOF. Just consider the equation $g = \delta f + \Lambda\delta g$ and write

$$S_\Pi = \sum_{i=0}^{n} g_{t_{i+1} t_i} = \sum_{i=0}^{n} (\delta f)_{t_{i+1} t_i} + \sum_{i=0}^{n} (\Lambda\delta g)_{t_{i+1} t_i}$$
$$= (\delta f)_{ts} + \sum_{i=0}^{n} (\Lambda\delta g)_{t_{i+1} t_i}.$$

Then observe that, due to the fact that $\Lambda\delta g \in \mathcal{C}_2^{1+}(V)$, the last sum converges to zero. □

2.2. *Computations in $\mathcal{C}_*$.* For sake of simplicity, let us assume, until Section 3, that $V = \mathbb{R}$, and set $\mathcal{C}_k(\mathbb{R}) = \mathcal{C}_k$. Then the complex $(\mathcal{C}_*, \delta)$ is an (associative, noncommutative) graded algebra once endowed with the following product: for $g \in \mathcal{C}_n$ and $h \in \mathcal{C}_m$ let $gh \in \mathcal{C}_{n+m-1}$ the element defined by

$$(gh)_{t_1, \ldots, t_{m+n-1}} = g_{t_1, \ldots, t_n} h_{t_n, \ldots, t_{m+n-1}}, \qquad t_1, \ldots, t_{m+n+1} \in [0, T]. \tag{16}$$

In this context, the coboundary $\delta$ act as a graded derivation with respect to the algebra structure. In particular, we have the following useful properties.

PROPOSITION 2.6. *The following differentiation rules hold true:*

(1) *Let $g, h$ be two elements of $\mathcal{C}_1$. Then*

$$\delta(gh) = \delta g\, h + g\, \delta h. \tag{17}$$

(2) *Let $g \in \mathcal{C}_1$ and $h \in \mathcal{C}_2$. Then*

$$\delta(gh) = \delta g\, h + g\, \delta h, \qquad \delta(hg) = \delta h\, g - h\, \delta g.$$

PROOF. We will just prove (17), the other relations being equally trivial: if $g, h \in \mathcal{C}_1$, then

$$[\delta(gh)]_{ts} = g_t h_t - g_s h_s = g_t(h_t - h_s) + (g_t - g_s) h_s = g_t (\delta h)_{ts} + (\delta g)_{ts} h_s,$$

which proves our claim. □

The iterated integrals of smooth functions on $[0,T]$ are obviously particular cases of elements of $\mathcal{C}$ which will be of interest for us, and let us recall some basic rules for these objects: consider $f,g \in \mathcal{C}_1^\infty$, where $\mathcal{C}_1^\infty$ is the set of smooth functions from $[0,T]$ to $\mathbb{R}$. Then the integral $\int dg\, f$, which will be denoted by $\mathcal{J}(dg\, f)$, can be considered as an element of $\mathcal{C}_2^\infty$. That is, for $s,t \in [0,T]$, we set

$$\mathcal{J}_{ts}(dg\, f) = \left(\int dg\, f\right)_{ts} = \int_s^t dg_u\, f_u.$$

The multiple integrals can also be defined in the following way: given a smooth element $h \in \mathcal{C}_2^\infty$ and $s,t \in [0,T]$, we set

$$\mathcal{J}_{ts}(dg\, h) \equiv \left(\int dg\, h\right)_{ts} = \int_s^t dg_u\, h_{us}.$$

In particular, the double integral $\mathcal{J}_{ts}(df^3\, df^2\, f^1)$ is defined, for $f^1, f^2, f^3 \in \mathcal{C}_1^\infty$, as

$$\mathcal{J}_{ts}(df^3\, df^2\, f^1) = \left(\int df^3\, df^2\, f^1\right)_{ts} = \int_s^t df_u^3\, \mathcal{J}_{us}(df^2\, f^1)$$

and if $f^1, \dots, f^{n+1} \in \mathcal{C}_1^\infty$, we set

$$\mathcal{J}_{ts}(df^{n+1}\, df^n \cdots df^2\, f^1) = \int_s^t df_u^{n+1}\, \mathcal{J}_{us}(df^n \cdots df^2\, f^1), \tag{18}$$

which defines the iterated integrals of smooth functions recursively.

The following relations between multiple integrals and the operator $\delta$ will also be useful in the remainder of the paper.

Proposition 2.7. *Let $f,g$ be two elements of $\mathcal{C}_1^\infty$. Then, recalling the convention (16), it holds that*

$$\delta f = \mathcal{J}(df), \qquad \delta(\mathcal{J}(dg\, f)) = 0,$$
$$\delta(\mathcal{J}(dg\, df)) = (\delta g)(\delta f) = \mathcal{J}(dg)\mathcal{J}(df)$$

*and, in general,*

$$\delta(\mathcal{J}(df^n \cdots df^1)) = \sum_{i=1}^{n-1} \mathcal{J}(df^n \cdots df^{i+1})\mathcal{J}(df^i \cdots df^1).$$

Proof. Here again, the proof is elementary, and we will just show the third of these relations: we have, for $s,t \in [0,T]$,

$$\mathcal{J}_{ts}(dg\, df) = \int_s^t dg_u\, (f_u - f_s) = \int_s^t dg_u\, f_u - K_{ts},$$

with $K_{ts} = (g_t - g_s) f_s$. The first term of the right-hand side is easily seen to be in $\mathcal{Z}\mathcal{C}_2$. Thus,

$$\delta(\mathcal{J}(dg\, df))_{tus} = -(\delta K)_{tus} = [g_t - g_u][f_u - f_s],$$

which gives the announced result. □

2.3. *Dissection of an integral.* The purpose of this section is not to provide an account on all the computations contained in [6]. However, we will go into some semi-heuristic considerations that, hopefully, will shed some light on the way we will solve rough PDEs later on: with the notation of Section 2.2 in mind, we will try to give, intuitively speaking, a meaning to the integral $\int \varphi(x)\, dx = \mathcal{J}(dx\, \varphi(x))$ for a nonsmooth function $x \in \mathcal{C}_1$. Notice that, in the sequel, $x$ should be considered as a vector valued function, since the whole theory can be handled via the Doss–Soussman methodology in the real case. However, we will present the main ideas of the algorithm below as if $x$ were real valued, the generalization from $\mathbb{R}$ to $\mathbb{R}^n$ being just a matter of (cumbersome) notation.

2.3.1. *The Young case.* The first idea one can have in mind in order to define $\mathcal{J}(dx\, \varphi(x))$ is to perform an expansion around the increment $dx$: indeed, in the smooth case, we have

$$\mathcal{J}(dx\, \varphi(x)) = \delta x\, \varphi(x) + \mathcal{J}(dx\, d\varphi(x)). \tag{19}$$

If we wish to extend the right-hand side of (19) to a nonsmooth case, we see that the first term is harmless, since it is defined independently of the regularity of $x$, by

$$[\delta x\, \varphi(x)]_{ts} = [x_t - x_s]\varphi(x_s) \qquad \text{for } s, t \in [0, T].$$

The last term of (19) is more problematic and we proceed to its *dissection* by the application of $\delta$: invoking Proposition 2.7, we get, in the smooth case, that

$$\begin{gathered} \delta(\mathcal{J}(dx\, d\varphi(x))) = \delta x\, \delta(\varphi(x)), \text{ that is,} \\ [\delta(\mathcal{J}(dx\, d\varphi(x)))]_{tus} = [\delta x]_{tu}[\delta(\varphi(x))]_{us}. \end{gathered} \tag{20}$$

Now the r.h.s. of (20) is well defined independently of the regularity of $x$. Thus, if $\delta x \delta(\varphi(x)) \in \mathcal{C}_3^{1+}$, which happens when $x \in \mathcal{C}_1^{\alpha}$ with $\alpha > \frac{1}{2}$ and $\varphi \in C^1(\mathbb{R})$, then Proposition 2.3 can be applied, and $\Lambda[\delta x\, \delta(\varphi(x))]$ is defined unambiguously. Hence, owing to (20), we set

$$\mathcal{J}(dx\, d\varphi(x)) = \Lambda(\delta x\, \delta(\varphi(x)))$$

and

$$\mathcal{J}(dx\, \varphi(x)) = \delta x\, \varphi(x) + \Lambda(\delta x\, \delta(\varphi(x))) = (\mathrm{Id} - \Lambda\delta)[\delta x\, \varphi(x)], \tag{21}$$

where the last equality is due to Proposition 2.6 and to the fact that $\delta\delta x = 0$. Notice once again that this construction is valid whenever $x \in \mathcal{C}_1^\alpha$ with $\alpha > \frac{1}{2}$ and $\varphi \in C^1(\mathbb{R})$, and it is easily shown, along the same lines as in the proof of Proposition 2.3 that the integral $\mathcal{J}(dx\,\varphi(x))$ defined by (21) corresponds to the usual Young integral.

2.3.2. *Case of a $\alpha$-Hölder path with $\frac{1}{3} < \alpha < \frac{1}{2}$.* The construction (21) does not work if $x \notin \mathcal{C}_1^{1/2+}$. However, if $x \in \mathcal{C}_1^\alpha$ with $\alpha > \frac{1}{3}$, we can proceed further in the expansion of equation (19) by observing that, still in the smooth case, we have, for $s,t \in [0,T]$,

$$\int_s^t [d\varphi(x)]_u = \int_s^t dx_u\,\varphi'(x_u) = [x_t - x_s]\varphi'(x_s) + \int_s^t dx_u \int_s^u dx_v\,\varphi''(x_v),$$

or according to the notation of Section 2.2,

$$\delta\varphi(x) = \mathcal{J}(d\varphi(x)) = \mathcal{J}(dx\,\varphi'(x)) = \delta x\,\varphi'(x) + \mathcal{J}(dx\,d\varphi'(x)). \tag{22}$$

Injecting this equality in equation (19), thanks to (18), we obtain

$$\mathcal{J}(dx\,\varphi(x)) = \delta x\,\varphi(x) + \mathcal{J}(dx\,dx)\varphi'(x) + \mathcal{J}(dx\,dx\,d\varphi'(x)). \tag{23}$$

Let us assume now that we are given a process $\mathcal{J}(dx\,dx) \in \mathcal{C}_2$, usually (and somewhat improperly) called the Lévy area of $x$, such that

$$\delta(\mathcal{J}(dx\,dx)) = \delta x\,\delta x \quad \text{and} \quad \mathcal{J}(dx\,dx) \in \mathcal{C}_2^{2\alpha}. \tag{24}$$

This assumption is of course not automatically satisfied, but it can be checked for instance in the Brownian and fractional Brownian cases. Then the right-hand side of (23) is again well defined independently of the regularity of $x$, except for the last term. However, recast equation (23) as

$$-\mathcal{J}(dx\,dx\,d\varphi'(x)) = -\mathcal{J}(dx\,\varphi(x)) + \delta x\,\varphi(x) + \mathcal{J}(dx\,dx)\varphi'(x),$$

and apply again $\delta$ to both sides of this last expression. Invoking Proposition 2.7 and recalling that $\delta(\mathcal{J}(dx\,dx)) = \delta x\,\delta x$, we obtain

$$\begin{aligned} -\delta\mathcal{J}(dx\,dx\,d\varphi'(x)) &= -\delta x\,\delta\varphi(x) + \delta x\,\delta x\,\varphi'(x) - \mathcal{J}(dx\,dx)\delta\varphi'(x) \\ &= -\delta x[\delta\varphi(x) - \delta x\,\varphi'(x)] - \mathcal{J}(dx\,dx)\delta\varphi'(x). \end{aligned} \tag{25}$$

Everything in the r.h.s. of equation (25) is well defined at this stage, and if we assume that all the terms belong to $\mathcal{C}_3^\mu$ with $\mu > 1$ [which can be justified via Taylor's expansions whenever $x \in \mathcal{C}_1^\alpha$ with $\alpha > \frac{1}{3}$ and $\varphi \in C^2(\mathbb{R})$], we can conclude that

$$\begin{aligned} \mathcal{J}(dx\,\varphi(x)) &= \delta x\,\varphi(x) + \mathcal{J}(dx\,dx)\varphi'(x) \\ &\quad + \Lambda[\mathcal{J}(dx\,dx)\delta\varphi'(x) + \delta x(\delta\varphi(x) - \delta x\,\varphi'(x))], \end{aligned}$$

or stated otherwise

$$\mathcal{J}(dx\,\varphi(x)) = (\mathrm{Id} - \Lambda\delta)[\delta x\,\varphi(x) + \mathcal{J}(dx\,dx)\varphi'(x)],$$

where we used the fact that $\delta\mathcal{J}(dx\,dx) = \delta x\,\delta x$ to put in evidence the fact that we are actually integrating (in the sense of Corollary 2.5) the 1-increment $\delta x\,\varphi(x) + \mathcal{J}(dx\,dx)\varphi'(x)$ which can be thought of as a *corrected* version of the more natural integrand $\delta x\,\varphi(x)$. It is worth noticing at that point that this integral has now to be understood as an integral over the (step-2) rough path $(x, \mathcal{J}(dx\,dx))$ introduced in [6] and it coincides with the notion of integral over a rough path given by Lyons in [17].

REMARK 2.8. This algorithm has an obvious extension to higher orders if we assume that a reasonable definition of the iterated integrals $\mathcal{J}(dx\,dx\cdots dx)$ can be given. To proceed further, however, we need the notion of *geometric* rough path (for more details on this notion see [17]) which must be exploited crucially to show that some terms are small enough and belong to the domain of $\Lambda$. For a more general approach, which does not rely on geometric rough-path, see [10].

**3. Algebraic integration associated to a semigroup.** The aim of this section is to set the basis for our future computations: after recalling some basic facts about analytic semigroups, we will define a set of increments $\hat{\mathcal{C}}_*$ and a modified operator $\hat{\delta}$ adapted to our evolution setting. Then we will give some basic calculus rules for $(\hat{\mathcal{C}}_*, \hat{\delta})$ and eventually, we will fix the notation for the main application we have chosen, that is the stochastic heat equation.

3.1. *Analytical semigroups.* As in [12], we will be able to develop our integration theory in the abstract setting of analytical semigroups on Banach spaces, whose basic features can be summarized as follows: let $(\mathcal{B}, |\cdot|)$ be a separable Banach space, and $(A, \mathrm{Dom}(A))$ be a nonbounded linear operator on $\mathcal{B}$. We will assume in the sequel that (see [20, Sections 2.5 and 2.6]) $A$ is the generator of an analytical semigroup $\{S_t; t \geq 0\}$, satisfying

$$|S_t| \leq Me^{-\lambda t} \qquad \text{for some constants } M, \lambda > 0 \text{ and for all } t \geq 0,$$

where $|\cdot|$ also stands for the operator norm on $\mathcal{B}$. Set now $A_o = -A$. This allows us, in particular, to define the fractional powers $(A_o^\alpha, \mathrm{Dom}(A_o^\alpha))$ for any $\alpha \in \mathbb{R}$.

For $\alpha \geq 0$, let $\mathcal{B}_\alpha$ be the space $\mathrm{Dom}(A_o^\alpha)$ with the norm $|x|_{\mathcal{B}_\alpha} = |A_o^\alpha x|$. Since $A_o^{-\alpha}$ is continuous, it follows that the norm $|\cdot|_{\mathcal{B}_\alpha}$ is equivalent to the graph norm of $A_o^\alpha$. If $\alpha = 0$, then $\mathcal{B}_\alpha = \mathcal{B}$ and $A_o^0 = \mathrm{Id}$. If $\alpha < 0$, let $\mathcal{B}_\alpha$ be the completion of $\mathcal{B}$ with respect to $|x|_{\mathcal{B}_\alpha} = |A_o^\alpha x|$, which means in particular that $\mathcal{B}_\alpha$ is a larger space than $\mathcal{B}$. We will also set $\mathcal{B}_{-\infty} = \bigcup_{\alpha\in\mathbb{R}} \mathcal{B}_\alpha$.

Among the important facts about these spaces, note the following ones: For any $\alpha \in \mathbb{R}$ and any $\rho \geq 0$,

$$A_o^{-\rho} \text{ maps } \mathcal{B}_\alpha \text{ onto } \mathcal{B}_{\alpha+\rho} \qquad \text{for all } \alpha \in \mathbb{R}, \rho \geq 0, \tag{26}$$

$$|x|_{\mathcal{B}_\alpha} \leq C_{\alpha,\rho} |x|_{\mathcal{B}_\rho} \qquad \text{for all } x \in \mathcal{B}_\alpha \text{ and all } \alpha \leq \rho. \tag{27}$$

Moreover, for all $\alpha, \beta \in \mathbb{R}$,

$$A_o^\alpha A_o^\beta = A_o^{\alpha+\beta} \qquad \text{on } \mathcal{B}_\gamma \tag{28}$$

with $\gamma = \max\{\alpha, \beta, \alpha+\beta\}$. The semigroup $(S_t)_{t \geq 0}$ also satisfies

$$S_t \text{ may be extended to } \mathcal{B}_\alpha \qquad \text{for all } \alpha < 0 \text{ and all } t > 0, \tag{29}$$

$$S_t \text{ maps } \mathcal{B}_\alpha \text{ to } \mathcal{B}_\rho \qquad \text{for all } \alpha \in \mathbb{R}, \rho \geq 0, t > 0, \tag{30}$$

$$\text{for all } t > 0,\ \alpha \geq 0 \qquad |A_o^\alpha S_t| \leq M_\alpha t^{-\alpha} e^{-\lambda t}, \tag{31}$$

$$\text{for } 0 < \alpha \leq 1,\ x \in \mathcal{B}_\alpha \qquad |S_t x - x| \leq C_\alpha t^\alpha |A_o^\alpha x|. \tag{32}$$

We will denote with $\mathcal{L}(\mathcal{B}, \mathcal{B}')$ the space of continuous linear operators from the Banach space $\mathcal{B}$ to the Banach space $\mathcal{B}'$. We let $\mathcal{L}(\mathcal{B}) = \mathcal{L}(\mathcal{B}, \mathcal{B})$. In order to be coherent with our previous notation, we also set $S_{t-s} = S_{ts}$ for a generic semigroup $S$, and $0 \leq s < t \leq T$.

3.2. *Convolutional increments.* Let us turn now to the main concern of this section, that is the definition of a complex $(\hat{\mathcal{C}}_*, \hat{\delta})$ which behaves nicely for the definition of our evolution problem.

Notice that, due to the fact that the operator $S_{t_1 t_2}$ is well defined only for $t_1 > t_2$, our integration domains will be of the form $\mathcal{S}_n$, where $\mathcal{S}_n$ stands for the $n$-simplex

$$\mathcal{S}_n = \{(t_1, \ldots, t_n) : T \geq t_1 \geq t_2 \geq \cdots \geq t_n \geq 0\}.$$

Let then $V$ be a separable Banach space. The basic family of increments we will work with is $\{\hat{\mathcal{C}}_n(V); n \geq 0\}$, where $\hat{\mathcal{C}}_n(V)$ denotes the space of continuous functions from $\mathcal{S}_n$ to $V$. Observe that an operator $\delta : \hat{\mathcal{C}}_n(V) \to \hat{\mathcal{C}}_{n+1}(V)$ can be defined just like in (7). In particular, if $A \in \hat{\mathcal{C}}_1(V)$ and $B \in \hat{\mathcal{C}}_2(V)$, the relation (8) is still valid. However, let us see now why $\delta$ is not adapted to the resolution of equation (4).

What made $\delta$ an interesting operator in Section 2 was the simple fact that, if $F \in \hat{\mathcal{C}}_1^\infty(\mathbb{R})$, then, for $t, s \in [0, T]^2$, we have

$$[\delta F]_{ts} = \int_s^t f_u \, du \qquad \text{with } f = F'. \tag{33}$$

However, if $S_t$ is the semigroup defined at Section 3.1, and if we set

$$\hat{F}_t = \int_0^t S_{tu} f_u \, du \qquad \text{for } t \geq 0, f \in \hat{\mathcal{C}}_1^\infty(\mathcal{B}),$$

then the same kind of relation does not hold true for $\hat{F}$. Indeed, for $s \le t$, if we define the operator $a_{ts}: \mathcal{B} \to \mathcal{B}$ as

$$a_{ts} = S_{ts} - \mathrm{Id}, \tag{34}$$

where $\mathrm{Id}: \mathcal{B} \to \mathcal{B}$ is the identity operator, then it is easily seen that

$$[\delta \hat{F}]_{ts} = \hat{F}_t - \hat{F}_s = a_{ts}\hat{F}_s + \int_s^t S_{tu} f_u \, du,$$

and hence, in order to get a similar relation to (33) in this new context, one should consider an operator $\hat{\delta}: \hat{\mathcal{C}}_n(\mathcal{B}) \to \hat{\mathcal{C}}_{n+1}(\mathcal{B})$, defined by

$$\begin{aligned} [\hat{\delta} A]_{t_1 \cdots t_{n+1}} &= [\delta A]_{t_1 \cdots t_{n+1}} - a_{t_1 t_2} A_{t_2 \cdots t_{n+1}} \\ &\qquad \text{for } A \in \hat{\mathcal{C}}_n(\mathcal{B}), (t_1 \cdots t_{n+1}) \in \mathcal{S}_{n+1}. \end{aligned} \tag{35}$$

In the remainder of the paper, we will write $\hat{\delta} A = \delta A - aA$, where we made use of the convention (16). As in Section 2.1, one can define, for $n \ge 1$,

$$\mathcal{Z}\hat{\mathcal{C}}_n(\mathcal{B}) = \hat{\mathcal{C}}_n(\mathcal{B}) \cap \ker(\hat{\delta}) \quad \text{and} \quad \mathcal{B}\hat{\mathcal{C}}_n(\mathcal{B}) = \hat{\mathcal{C}}_n(\mathcal{B}) \cap \mathrm{Im}(\hat{\delta}).$$

Then the perturbed operator $\hat{\delta}$ preserves some important properties of the original coboundary $\delta$.

PROPOSITION 3.1. *The couple* $(\hat{\mathcal{C}}_*, \hat{\delta})$ *is an acyclic cochain complex:* $\mathcal{Z}\hat{\mathcal{C}}_{n+1} = \mathcal{B}\hat{\mathcal{C}}_n$ *for any* $n \ge 0$.

PROOF. Let us prove first that $\hat{\delta}$ is a coboundary, that is, $\hat{\delta}\hat{\delta} = 0$. Indeed, if $F \in \hat{\mathcal{C}}_n$ according to the fact that $\delta\delta = 0$ and thanks to the forthcoming Lemma 3.2, we have

$$\begin{aligned} \hat{\delta}\hat{\delta} F &= (\delta - a)[(\delta - a)F] = \delta\delta F - \delta(aF) - a\delta F + aaF \\ &= -\delta a F + a\delta F - a\delta F + aaF = aaF - \delta a F. \end{aligned}$$

Furthermore, it is readily checked that

$$(\delta a)_{tus} = a_{tu} a_{us}, \qquad (t, u, s) \in \mathcal{S}_3,$$

which gives $\hat{\delta}\hat{\delta} F = 0$.

The fact that $\mathrm{Im}\, \hat{\delta}_{|\hat{\mathcal{C}}_n} = \ker \hat{\delta}_{|\hat{\mathcal{C}}_{n+1}}$ can be proved along the same lines as for the $(\mathcal{C}_*, \delta)$ complex [6]: pick $A \in \hat{\mathcal{C}}_{n+1}$ such that $\hat{\delta} A = 0$, and set $B_{t_1 \ldots t_n} = A_{t_1 \ldots t_n s}$, with $s = 0$. Then

$$\begin{aligned} [\hat{\delta} B]_{t_1 \cdots t_{n+1}} &= [\delta A]_{t_1 \cdots t_{n+1} s} + (-1)^{n+1} A_{t_1 \cdots t_{n+1}} - a_{t_1 t_2} A_{t_2 \cdots t_n s} \\ &= [\hat{\delta} A]_{t_1 \cdots t_{n+1} s} + (-1)^{n+1} A_{t_1 \cdots t_{n+1}} = (-1)^{n+1} A_{t_1 \cdots t_{n+1}}. \end{aligned}$$

Thus, setting $C=(-1)^{n+1}B$, we get $\hat\delta C=A$. $\square$

The cochain complex $(\hat{\mathcal{C}}_*,\hat\delta)$ will be the structure at the base of all the constructions in this paper. Let us also mention at this point that, when the meaning is obvious, we will transpose the notation of Section 2 to our infinite-dimensional setting. Furthermore, whenever this does not lead to an ambiguous situation, we will write $\hat{\mathcal{C}}_n$ instead of $\hat{\mathcal{C}}_n(\mathcal{B})$.

Let us give now a simple and useful extension of Proposition 2.6, which has already been used in the last proposition.

LEMMA 3.2. *Let $L\in\hat{\mathcal{C}}_{n-1}(\mathcal{B})$ and $M\in\hat{\mathcal{C}}_2(\mathcal{L}(\mathcal{B}))$. Then*

$$\delta(ML)=\delta ML-M\delta L.$$

PROOF. Let $G_{t_1\cdots t_n}=M_{t_1t_2}L_{t_2\cdots t_n}$. Then

$$\begin{aligned}
[\delta G]_{t_1\cdots t_{n+1}}&=\sum_{i=1}^{n+1}(-1)^iG_{t_1\cdots\hat t_i\cdots t_{n+1}}\\
&=-M_{t_2t_3}L_{t_3\cdots t_{n+1}}+M_{t_1t_3}L_{t_3\cdots t_{n+1}}+\sum_{i=3}^{n+1}(-1)^iM_{t_1t_2}L_{t_2\cdots\hat t_i\cdots t_{n+1}}\\
&=[\delta M]_{t_1t_2t_3}L_{t_3\cdots t_{n+1}}+M_{t_1t_2}\sum_{i=2}^{n+1}(-1)^iL_{t_2\cdots\hat t_i\cdots t_{n+1}},
\end{aligned}$$

which yields our claim. $\square$

3.3. *Computations in $\hat{\mathcal{C}}_*$.* Here again, like in Section 2, we will try to move from a smooth setting to an irregular one. And we will start by giving the equivalent, in our new setting, of Proposition 2.3, which will require first the introduction of some analytical structures on the spaces $\hat{\mathcal{C}}_n$.

3.3.1. *Hölder type spaces.* First of all, we have to define some Hölder type subspaces of $\hat{\mathcal{C}}_k$, $k\le 3$, related to the spaces $\mathcal{B}_\alpha$, $\alpha\in\mathbb{R}$: for $\mu\ge 0$ and $g\in\hat{\mathcal{C}}_2(\mathcal{B}_\alpha)$, we set

$$\|g\|_{\mu,\alpha}\equiv\sup_{t,s\in\mathcal{S}_2}\frac{|g_{ts}|_{\mathcal{B}_\alpha}}{|t-s|^\mu}\quad\text{and}\quad\hat{\mathcal{C}}_2^{\mu,\alpha}=\{g\in\hat{\mathcal{C}}_2(\mathcal{B}_\alpha);\|g\|_{\mu,\alpha}<\infty\},\tag{36}$$

and the definition above also induces some seminorms on $\mathcal{C}_1$: for $\gamma>0,\alpha\in\mathbb{R}$, we say that $f\in\hat{\mathcal{C}}_1^{\gamma,\alpha}$ if

$$\|f\|_{\gamma,\alpha}\equiv\|\hat\delta f\|_{\gamma,\alpha}<\infty.$$

Another useful subspace of $\hat{\mathcal{C}}_1$ will be $\hat{\mathcal{C}}_1^{0,\alpha}$, the space of bounded paths in $\mathcal{B}_\alpha$ with the supremum norm $\|f\|_{0,\alpha}=\sup_{t\in[0,T]}|f_t|_{\mathcal{B}_\alpha}$.

As far as $\hat{\mathcal{C}}_3$ is concerned, $\hat{\mathcal{C}}_3^{\mu,\alpha}$ can be defined in the following way: set

$$\begin{aligned}
\|h\|_{\gamma,\rho,\alpha} &= \sup_{t,u,s\in\mathcal{S}_3} \frac{|h_{tus}|_{\mathcal{B}_\alpha}}{|t-u|^\gamma |u-s|^\rho}, \\
\|h\|_{\mu,\alpha} &\equiv \inf\left\{\sum_i \|h_i\|_{\rho_i,\mu-\rho_i,\alpha}; h=\sum_i h_i, 0<\rho_i<\mu\right\},
\end{aligned} \tag{37}$$

where the last infimum is taken over all sequences $\{h_i\}_i$ such that $h=\sum_i h_i$ and for all choices of the numbers $\rho_i\in(0,\mu)$. Then $\|\cdot\|_{\mu,\alpha}$ is again easily seen to be a norm, and we set

$$\hat{\mathcal{C}}_3^{\mu,\alpha} = \{h\in\hat{\mathcal{C}}_3(\mathcal{B}_\alpha); \|h\|_{\mu,\alpha}<\infty\}.$$

In order to avoid ambiguities, we shall also denote in the sequel by $\mathcal{N}[f;\hat{\mathcal{C}}_j^\kappa]$ the $\kappa$-Hölder norm on the space $\hat{\mathcal{C}}_j$, for $j=1,2,3$. For $\zeta\in\hat{\mathcal{C}}_1(V)$, we also set $\mathcal{N}[\zeta;\hat{\mathcal{C}}_1^0(V)]=\sup_{s\in[0;T]}\|\zeta_s\|_V$.

Eventually, we will need to introduce a slight extension of the spaces we have just defined above: for $j=1,2$, let $\mathcal{E}_j^{\mu,\alpha}$ be defined by

$$\mathcal{E}_j^{\mu,\alpha} = \bigcap_{\varepsilon\le\mu\wedge 1^-} \hat{\mathcal{C}}_j^{\mu-\varepsilon,\alpha+\varepsilon}, \tag{38}$$

where $\varepsilon\le\mu\wedge 1^-$ stands for the condition $\varepsilon\in[0,\mu]\cap[0,1)$, and where the intersection is considered along any arbitrary family $\{0\le\varepsilon_1<\cdots<\varepsilon_n\le\mu\wedge 1^-\}$ for $n\ge 1$. Obviously, some families of operators will play an important role in the sequel, and this will lead us to the following specific definitions for operator-valued increments.

Definition 3.3. For $\mu\ge 0$ and $\alpha,\beta\in\mathbb{R}$, we will call $\hat{\mathcal{C}}_2^\mu\mathcal{L}^{\beta,\alpha}$ the space $\hat{\mathcal{C}}_2^\mu(\mathcal{L}(\mathcal{B}_\beta;\mathcal{B}_\alpha))$, and will denote by $\mathcal{E}_2^\mu\mathcal{L}^{\beta,\alpha}$ the space

$$\mathcal{E}_2^\mu\mathcal{L}^{\beta,\alpha} = \bigcap_{\varepsilon\le\mu\wedge 1^-} \hat{\mathcal{C}}_2^{\mu-\varepsilon}\mathcal{L}^{\beta,\alpha+\varepsilon},$$

where the intersection is still considered along any arbitrary finite family $\{0\le\varepsilon_1<\cdots<\varepsilon_n\le\mu\wedge 1^-\}$ for $n\ge 1$. The natural norm on $\hat{\mathcal{C}}_2^\mu\mathcal{L}^{\beta,\alpha}$ will be defined by

$$\|A\|_{\mu,\beta,\alpha} = \sup_{t,s\in\mathcal{S}_2} \frac{\|A_{ts}\|_{\mathrm{op}}}{|t-s|^\mu}, \tag{39}$$

and when we consider some Hilbert–Schmidt operators, the corresponding spaces will be denoted by $\hat{\mathcal{C}}_2^\mu\mathcal{L}_{\mathrm{HS}}^{\beta,\alpha}$ and $\mathcal{E}_2^\mu\mathcal{L}_{\mathrm{HS}}^{\beta,\alpha}$.

3.3.2. *The convolution sewing map and related properties.* Here is a first proposition showing how the analytical structures introduced above interact with our previous algebraic notation.

PROPOSITION 3.4. *If $\mu > 1$, then for any $\alpha \in \mathbb{R}$, $\mathcal{Z}\hat{\mathcal{C}}_2^{\mu,\alpha} = \{0\}$.*

PROOF. Take $h \in \mathcal{Z}\hat{\mathcal{C}}_2^{\mu,\alpha}$. Then, according to Proposition 3.1, there exists $f \in \hat{\mathcal{C}}_1$ such that $h = \hat{\delta} f$. Consider the telescopic sum

$$h_{ts} = (\hat{\delta} f)_{ts} = \sum_{i=0}^{n} S_{tt_{i+1}} (\hat{\delta} f)_{t_{i+1} t_i},$$

with respect to the partition $\Pi_{ts}^n = \{t_{0 \le i \le n+1} : t_0 = s, t_{n+1} = t\}$ of the interval $[s,t]$. Since $\hat{\delta} f \in \mathcal{Z}\hat{\mathcal{C}}_2^{\mu,\alpha}$ with $\mu > 1$, we have

$$|(\hat{\delta} f)_{ts}|_{\mathcal{B}_\alpha} \le \sum_{i=0}^{n} |(\hat{\delta} f)_{t_{i+1} t_i}|_{\mathcal{B}_\alpha} \le \|\hat{\delta} f\|_{\mu,\alpha} \sum_{i=0}^{n} |t_{i+1} - t_i|^{\mu}$$

which converges to zero as the size of the partition goes to zero. Since $t, s$ are arbitrary, we have $\hat{\delta} f = h = 0$ in $\hat{\mathcal{C}}_2^{\mu,\alpha}$. $\square$

We can now state and prove the equivalent of Proposition 2.3 in our evolution setting, which is the main aim of this section.

THEOREM 3.5. *Let $\mu > 1$, $\alpha \in \mathbb{R}$. There exists a unique sewing map $\hat{\Lambda} : \mathcal{Z}\hat{\mathcal{C}}_3^{\mu,\alpha} \to \mathcal{E}_2^{\mu,\alpha}$ such that $\hat{\delta}\hat{\Lambda} = \mathrm{Id}_{\mathcal{Z}\hat{\mathcal{C}}_3}$. Furthermore, for any $0 \le \varepsilon \le \mu \wedge 1^-$, there exists a strictly positive constant $c_{\mu,\varepsilon}$ such that*

$$\|\hat{\Lambda} h\|_{\mu-\varepsilon,\alpha+\varepsilon} \le c_{\mu,\varepsilon} \|h\|_{\mu,\alpha}, \tag{40}$$

*for any $h \in \mathcal{Z}\hat{\mathcal{C}}_3^{\mu,\alpha}$.*

PROOF. Like in the proof of Proposition 2.3, we will divide our computations in two steps below.

*Step 1*: The uniqueness part of our theorem simply stems from the fact that if we have $\hat{\delta} a = h$ and $\hat{\delta} a' = h$ with $a, a' \in \hat{\mathcal{C}}_2^{\mu,\alpha}$, then $b = a - a' \in \mathcal{Z}\hat{\mathcal{C}}_2^{\mu,\alpha}$ and since $\mu > 1$, by Proposition 3.4, we must have $b = 0$.

*Step 2*: The existence part can be adapted from Proposition 2.3, and we will construct a process $M \in \mathcal{E}_2^{\mu,\alpha}$ such that $\hat{\delta} M = h$ starting from any $B \in \hat{\mathcal{C}}_2(\mathcal{B}_\alpha)$ satisfying $\hat{\delta} B = h$ (this increment $B$ exists thanks to Lemma 3.1). Now, similar to (13), we will set, for a given $n \ge 1$, and $(t,s) \in \mathcal{S}_2$,

$$M_{ts}^n = B_{ts} - \sum_{i=0}^{2^n - 1} S_{tr_{i+1}^n} B_{r_{i+1}^n, r_i^n},$$

where $s,t$ and $r_i^n$ have been defined at (12). Then $M_{ts}^0 = 0$ and

$$\begin{aligned}
&M_{ts}^{n+1} - M_{ts}^n \\
&\quad = \sum_{i=0}^{2^n-1} (S_{tr_{2i+2}^{n+1}} B_{r_{2i+2}^{n+1}, r_{2i}^{n+1}} - S_{tr_{2i+2}^{n+1}} B_{r_{2i+2}^{n+1}, r_{2i+1}^{n+1}} - S_{tr_{2i+1}^{n+1}} B_{r_{2i+1}^{n+1}, r_{2i}^{n+1}}) \\
&\quad = \sum_{i=0}^{2^n-1} S_{tr_{2i+2}^{n+1}} (B_{r_{2i+2}^{n+1}, r_{2i}^{n+1}} - B_{r_{2i+2}^{n+1}, r_{2i+1}^{n+1}} - B_{r_{2i+1}^{n+1}, r_{2i}^{n+1}}) \\
&\qquad\quad - S_{tr_{2i+2}^{n+1}} [S_{r_{2i+2}^{n+1} r_{2i+1}^{n+1}} - \mathrm{Id}] B_{r_{2i+1}^{n+1}, r_{2i}^{n+1}}.
\end{aligned}$$

Thus, according to the definition (35) of $\hat{\delta}$, we get

$$\begin{aligned}
M_{ts}^{n+1} - M_{ts}^n &= \sum_{i=0}^{2^n-1} S_{tr_{2i+2}^{n+1}} [(\delta B)_{r_{2i+2}^{n+1}, r_{2i+1}^{n+1}, r_{2i}^{n+1}} - a_{r_{2i+2}^{n+1}, r_{2i+1}^{n+1}} B_{r_{2i+1}^{n+1}, r_{2i}^{n+1}}] \\
&= \sum_{i=0}^{2^n-1} S_{tr_{2i+2}^{n+1}} (\hat{\delta} B)_{r_{2i+2}^{n+1}, r_{2i+1}^{n+1}, r_{2i}^{n+1}} = \sum_{i=0}^{2^n-1} S_{tr_{2i+2}^{n+1}} h_{r_{2i+2}^{n+1}, r_{2i+1}^{n+1}, r_{2i}^{n+1}}.
\end{aligned}$$

Hence, for any $\varepsilon < \mu$, we get, invoking (31),

$$\begin{aligned}
|A^{\alpha+\varepsilon}(M_{ts}^{n+1} - M_{ts}^n)| &\le c_\varepsilon \sum_{i=0}^{2^n-1} |t - r_{i+1}^n|^{-\varepsilon} |h|_{\mu,\alpha} |t-s|^\mu \\
&\le \frac{c_\varepsilon |t-s|^{\mu-\varepsilon} |h|_{\mu,\alpha}}{2^{n(\mu-1)}} \int_0^1 u^{-\varepsilon}\, du,
\end{aligned}$$

which gives, like in Proposition 2.3, that $M_{ts} \equiv \lim_{n\to\infty} M_{ts}^n$ exists, and is an element of $\mathcal{E}_2^{\mu,\alpha}$. Now, the fact that $\hat{\delta} M = h$ can be shown analogously to the case of Proposition 2.3, and the proof of (40) is straightforward. □

A direct consequence of the existence of the $\hat{\Lambda}$-map is a result of convergence of finite sums.

Corollary 3.6. *Let $g \in \hat{\mathcal{C}}_2$ such that $\hat{\delta} g \in \hat{\mathcal{C}}_3^{\mu,\alpha}$ for some $\mu > 1$. Then the 1-increment $\hat{\delta} f = (\mathrm{Id} - \hat{\Lambda}\hat{\delta}) g \in \hat{\mathcal{C}}_2^\alpha$ satisfies*

$$(\hat{\delta} f)_{ts} = \lim_{|\Pi_{ts}| \to 0} \sum_{i=0}^n S_{tt_{i+1}} g_{t_{i+1} t_i},$$

*for all $(t,s) \in \mathcal{S}_2$.*

Proof. It follows the lines of the proof of Corollary 2.5. □

We will now define an equivalent of the iterated integrals of Section 2.2 in our convolution context: consider some smooth functions $g \in \hat{\mathcal{C}}_1^\infty(\mathcal{L}(\mathcal{B}_\alpha))$ and $f \in \hat{\mathcal{C}}_1^\infty(\mathcal{B}_\alpha)$, for some $\alpha \in \mathbb{R}$. Then $\mathcal{J}(dg\, f)$ will be defined as an element of $\hat{\mathcal{C}}_2^\infty(\mathcal{B}_\alpha)$ by $\mathcal{J}_{ts}(dg\, f) = \int_s^t dg_v\, f_v$, for $(t,s) \in \mathcal{S}_2$. We will also need some integrals of processes weighted by the semigroup $S$, defined as follows, for $0 \le s < t \le T$:

$$\mathcal{J}_{ts}(\hat{d}g\, f) = \int_s^t S_{tv}\, dg_v\, f_v.$$

Once these elementary blocks have been defined, the iterated integrals

$$\mathcal{J}(d^{*n} f_n \cdots d^{*1} f_1) \quad \text{for } f_n, \ldots, f_2 \in \hat{\mathcal{C}}_1^\infty(\mathcal{L}(\mathcal{B}_\alpha)), f_1 \in \hat{\mathcal{C}}_1^\infty(\mathcal{B}_\alpha), \tag{41}$$

where $d^{*j} f_j$ stands for any of the increments of the form $df_j$ or $\hat{d} f_j$, can be defined recursively along the same lines as in Section 2.2. In particular, the operator-valued increment $\mathcal{J}(\hat{d}g\, S)$ is defined by

$$\mathcal{J}_{ts}(\hat{d}g\, S) = \int_s^t S_{tu}\, dg_u\, S_{us}.$$

The relations between $\hat{\delta}$ and these integrals, which will be useful for our purposes, can be summarized in the following:

PROPOSITION 3.7. *Let $\alpha \in \mathbb{R}$, and $g \in \hat{\mathcal{C}}_1^\infty(\mathcal{L}(\mathcal{B}_\alpha))$, $f \in \hat{\mathcal{C}}_1^\infty(\mathcal{B}_\alpha)$. Then*

$$\hat{\delta}(\mathcal{J}(\hat{d}f)) = 0, \qquad \hat{\delta}(\mathcal{J}(\hat{d}g\, f)) = 0, \qquad \hat{\delta}(\mathcal{J}(\hat{d}g\, \hat{\delta} f)) = \mathcal{J}(\hat{d}g\, S)\hat{\delta}(f)$$

*and*

$$\hat{\delta}(\mathcal{J}(\hat{d}g\, \hat{d}f)) = \mathcal{J}(\hat{d}g\, S)\mathcal{J}(\hat{d}f), \qquad \hat{\delta}(\mathcal{J}(\hat{d}g\, df)) = \mathcal{J}(\hat{d}g)\mathcal{J}(df).$$

PROOF. The proof of these results is elementary. We will give some details about the last relation for sake of completeness. For any $(t,u,s) \in \mathcal{S}_3$, invoking the definition of $\hat{\delta}$, we have

$$\begin{aligned}
[\hat{\delta}(\mathcal{J}(\hat{d}g\, df))]_{tus} &= [\delta(\mathcal{J}(\hat{d}g\, df))]_{tus} - a_{tu}\mathcal{J}_{us}(\hat{d}g\, df) \\
&= \int_s^t S_{tv}\, dg_v \int_s^v df_w - \int_u^t S_{tv}\, dg_v \int_u^v df_w \\
&\quad - \int_s^u S_{uv}\, dg_v \int_s^v df_w - [S_{tu} - \mathrm{Id}] \int_s^u S_{uv}\, dg_v \int_s^v df_w \\
&= \mathcal{J}_{tu}(\hat{d}g)\mathcal{J}_{us}(df),
\end{aligned}$$

which proves the claim. □

3.4. *Fractional heat equation setting.* In this section, we will give the general setting under which we will try to define and solve the stochastic heat equation driven by an infinite-dimensional fractional Brownian motion: as mentioned in the Introduction, the main application we have in mind is the situation where $A = \Delta - \mathrm{Id}$, and $\Delta$ is the Laplace operator on the circle $S$, assimilated to $[0,1]$. This operator can be diagonalized in the trigonometric basis of $L^2([0,1];\mathbb{C})$, namely $\{e_n; n \in \mathbb{Z}\}$, where

$$e_n(x) = e^{2\imath\pi n x}, \qquad x \in [0,1].$$

Associated to these eigenfunctions are the eigenvalues $\lambda_n = -1 - (2\pi n)^2$. We have chosen to deal with $A = \Delta - \mathrm{Id}$ instead of $\Delta$ itself for computational convenience, since this choice avoids the problem of a null eigenvalue for constant functions. Notice that in this case, $A$ is the generator of an analytical semigroup, and all the constructions of Section 3.1 goes through. Then $\mathcal{B}_\alpha$ can be identified with $H_\alpha$, the usual Sobolev space based on $L^2([0,1])$, for the definition of which we refer to Adams [1], and $\{S_t; t \geq 0\}$ stands for the heat semigroup, which admits a kernel $G_t(\xi,\eta)$ for $t > 0$ and $\xi,\eta \in [0,1]$. In this context, set $G_t^\alpha(\xi,\eta)$ for the kernel of the operator $A_o^\alpha S_t$, and $G^\beta(\xi,\eta)$ for the kernel of the operator $A_o^{-\beta}$. Then, for $\alpha \in \mathbb{R}$ and $\beta > 0$, $G_t^\alpha$ and $G^\beta$ admit the following spectral decomposition:

$$\begin{aligned} G_t^\alpha(\xi,\eta) &= \sum_{n\in\mathbb{Z}} \lambda_n^\alpha e^{-t\lambda_n} e_n(\xi)\bar{e}_n(\eta) \quad \text{and} \\ G^\beta(\xi,\eta) &= \sum_{n\in\mathbb{Z}} \lambda_n^{-\beta} e_n(\xi)\bar{e}_n(\eta). \end{aligned} \tag{42}$$

Let us specify now the noise $X$ we will consider: we will try to stick to the existing literature on the topic, and choose a fractional Brownian noise in time, defined on a certain complete probability space $(\Omega,\mathcal{F},P)$, which will be homogeneous in space, with a spatial covariance function $Q$. Namely, $X$ will be a centered Gaussian field indexed by functions on $[0,T]\times[0,1]$, such that if $\phi$ and $\psi$ are smooth enough, then

$$\begin{aligned} &E[X(\phi)X(\psi)] \\ &\quad = c_H \int_{[0,T]^2} \left( \int_{[0,1]^2} Q(\xi-\eta)\phi(u,\xi)\psi(v,\eta)\, d\xi\, d\eta \right) |u-v|^{2H-2}\, du\, dv, \end{aligned} \tag{43}$$

with $c_H = H(2H-1)$, for $H > \frac{1}{2}$. Notice that, in order to simplify our statements, we will generally assume that $Q$ can be decomposed itself on the basis $\{e_n; n \in \mathbb{Z}\}$ in the following way:

$$Q(\xi) = \sum_{n\in\mathbb{Z}} q_n e_n(\xi) \qquad \text{with } q_n = \lambda_n^{-\nu}, \text{ for } \nu \in [0,1), \tag{44}$$

and notice that the case $\nu = 0$ corresponds to the white noise in space, while the case $\nu > 1/2$ corresponds to a noise admitting a density in space. Some explicit construction of such kind of noise, as well as an account on the related stochastic calculus, can be found in [24]. The methodology we will develop in the rough case will also enable us to handle the Brownian motion case, which means a covariance structure given by

$$E[X(\phi)X(\psi)] = \int_{[0,T]} \left( \int_{[0,1]^2} \phi(u,\xi) Q(\xi - \eta) \psi(u,\eta) \, d\xi \, d\eta \right) du. \tag{45}$$

We give here a slight extension of a result result of [6], which will be used below to prove existence of regular versions of some stochastic processes, following the well-known approach of Garsia–Rodemich–Rumsey. The proof is conceptually similar to that appearing in [6] but there is a small technical difficulty due to the fact that convolutional increments are one-sided and which forces us to follow the scheme of the proof of the GRR inequality in Stroock's book rather that which can be found in [6].

In order to state this extension, we shall introduce for the first time a variant of the operator $\hat\delta$, called $\tilde\delta$, acting on operator-valued increments which turns out to be useful in the sequel, and which is defined by

$$\tilde\delta Q = \hat\delta Q - Qa = \delta Q - aQ - Qa, \qquad Q \in \mathcal{C}^*(\mathcal{L}(\mathcal{B})). \tag{46}$$

With this additional notation in hand, our regularity lemma is the following below.

Lemma 3.8. *For any $\gamma > 0$, $\alpha, \beta \in \mathbb{R}$ and $p \ge 1$, there exists a constant $C$ such that for any $R \in \hat{\mathcal{C}}_2(\mathcal{L}^{\beta,\alpha})$, we have*

$$\|R\|_{\gamma,\beta,\alpha} \le C(U_{\gamma+2/p,p,\beta,\alpha}(R) + \|\tilde\delta R\|_{\gamma,\beta,\alpha}), \tag{47}$$

*where*

$$U_{\gamma,p,\beta,\alpha}(R) = \left[ \int_{\mathcal{S}_2} \left( \frac{|R_{ts}|_{\beta,\alpha}}{|t-s|^\gamma} \right)^p dt \, ds \right]^{1/p}.$$

Proof. As in [6], this result is a direct consequence of a more general Lemma 3.9 below by choosing $\Phi(x) = x^p$ and $p(t) = t^{\gamma+2/p}$. □

Lemma 3.9. *Let $p$ and $\Psi$ be strictly increasing, continuous functions on $\mathbb{R}_+$ satisfying $p(0) = \Psi(0) = 0$ and $\Psi(x) \to \infty$ as $x \to \infty$. Then there exist a constant $K$ such that for any $\alpha, \beta \in \mathbb{R}$ and any $R \in \hat{\mathcal{C}}_2(\mathcal{L}^{\beta,\alpha})$ for which*

$$U = \iint_{0<s<t<T} \Phi\left( \frac{|R_{ts}|_{\beta,\alpha}}{p(t-s)} \right) dt \, ds < \infty$$

*and*

$$\sup_{s\le u\le t} |\tilde\delta R_{tus}|_{\beta,\alpha} \le \Phi^{-1}\biggl(\frac{4C}{(t-s)^2}\biggr)p(t-s), \qquad 0\le s\le t\le T$$

*for some constant $C<\infty$, then*

$$|R_{ts}|_{\beta,\alpha} \le 8K\int_0^{t-s}\Phi^{-1}\biggl(\frac{4B}{u^2}\biggr)p(du) + 9K\int_0^{t-s}\Phi^{-1}\biggl(\frac{4C}{u^2}\biggr)p(du)$$

*for all $0\le s\le t\le T$.*

PROOF. The proof follows closely Stroock's proof of the Garsia–Rodemich–Rumsey inequality. First, show the estimate for $T=t=1$ and $s=0$. For a sequence of times $t_n$,$s_n$ such that $t>t_{k+1}>t_k$, $s<s_{k+1}<s_k$ and $t_0=s_0$, we have

$$R_{ts} = R_{tt_0}S_{t_0s} + S_{tt_0}R_{t_0s} + \tilde\delta R_{tt_0s},$$

$$R_{tt_0} = R_{tt_{n+1}}S_{t_{n+1}t_0} + \sum_{k=0}^n S_{tt_{k+1}}R_{t_{k+1}t_k}S_{t_kt_0} + \sum_{k=0}^n \tilde\delta R_{tt_{k+1}t_k}S_{t_kt_0}$$

and

$$R_{t_0s} = S_{t_0s_{n+1}}R_{s_{n+1}s} + \sum_{k=0}^n S_{t_0s_k}R_{s_ks_{k+1}}S_{s_{k+1}s} + \sum_{k=0}^n S_{t_0s_k}\tilde\delta R_{s_ks_{k+1}s}.$$

Next, we choose these times as follows. Let $I(v)=\int_0^v \psi(\frac{|R_{vu}|_{\beta,\alpha}}{p(v-u)})\,du$. For any $s_n$, define $d_n$ by the equation $2p(d_n)=p(s_n)$. Remark that since $\int_0^1 I(t)\,dt = U$ there exists $t_0$ such that $I(t_0)\le U$. We claim that there exists $s_{n+1}\in(0,d_n)$ such that both inequalities

$$I(s_{n+1}) \le \frac{2U}{d_n} \quad\text{and}\quad \Phi\biggl(\frac{|R_{s_ns_{n+1}}|_{\beta,\alpha}}{p(s_n-s_{n+1})}\biggr) \le \frac{2I(s_n)}{d_n}$$

hold. This is always possible since, if we call $A_n\subset(0,d_n)$ (resp. $B_n$) the set of $s_{n+1}$ where the first (resp. the second) fail, we have

$$U \ge \int_{A_n} ds_{n+1}\, I(s_{n+1}) > \frac{2U}{d_n}|A_n| \quad\text{and}$$

$$I(s_n) \ge \int_{B_n} ds_{n+1}\,\Phi\biggl(\frac{|R_{s_ns_{n+1}}|_{\beta,\alpha}}{p(s_n-s_{n+1})}\biggr) > \frac{2I(s_n)}{d_n}|B_n|$$

so we must have $|A_n|<d_n/2$ and $|B_n|<d_n/2$ which means that $(0,d_n)\backslash(A_n\cup B_n)$ has positive measure. Then since

$$p(s_n-s_{n+1}) \le p(s_n) = 2p(d_n) = 4(p(d_n)-p(d_n)/2) \le 4(p(d_n)-p(d_{n+1})),$$

we have

$$
\begin{aligned}
|R_{s_n s_{n+1}}|_{\beta,\alpha} &\leq \Phi^{-1}\left(\frac{2I(s_n)}{d_n}\right)p(s_n - s_{n+1}) \\
&\leq 4\Phi^{-1}\left(\frac{4U}{d_n d_{n-1}}\right)(p(d_n) - p(d_{n+1})) \\
&\leq 4\int_{d_{n+1}}^{d_n} \Phi^{-1}\left(\frac{4U}{u^2}\right)p(du)
\end{aligned}
$$

and

$$
\begin{aligned}
|\tilde{\delta} R_{s_n s_{n+1} s}|_{\beta,\alpha} &\leq \Phi^{-1}\left(\frac{4C}{d_n^2}\right)p(s_n) \leq 4\Phi^{-1}\left(\frac{4C}{d_n^2}\right)(p(d_n) - p(d_{n+1})) \\
&\leq 4\int_{d_{n+1}}^{d_n} \Phi^{-1}\left(\frac{4C}{u^2}\right)p(du)
\end{aligned}
$$

then we have

$$
\begin{aligned}
|R_{t_0 s}|_{\beta,\alpha} &\leq M^2 \sum_{n=0}^{\infty} |R_{s_k s_{k+1}}| + M \sum_{n=0}^{\infty} |\hat{\delta} R_{s_k s_{k+1} s}| \\
&\leq \sum_{n=0}^{\infty} 4M^2 \int_{d_{n+1}}^{d_n} \Phi^{-1}\left(\frac{4U}{u^2}\right)p(du) + \sum_{n=0}^{\infty} 4M \int_{d_{n+1}}^{d_n} \Phi^{-1}\left(\frac{4C}{u^2}\right)p(du) \\
&\leq 4M^2 \int_0^{t-s} \Phi^{-1}\left(\frac{4U}{u^2}\right)p(du) + 4M \int_0^{t-s} \Phi^{-1}\left(\frac{4C}{u^2}\right)p(du),
\end{aligned}
$$

where we used the fact that there exists $M > 1$ such that $|S_\tau|_{\alpha,\alpha} \leq M$ for any $\alpha$ and any $\tau \geq 0$.

Similarly, we find

$$
|R_{tt_0}|_{\beta,\alpha} \leq 4M^2 \int_0^{t-s} \Phi^{-1}\left(\frac{4U}{u^2}\right)p(du) + 4M \int_0^{t-s} \Phi^{-1}\left(\frac{4C}{u^2}\right)p(du).
$$

So using that

$$
|\tilde{\delta} R_{tt_0 s}|_{\beta,\alpha} \leq \Phi^{-1}\left(\frac{4C}{(t-s)^2}\right)p(t-s) \leq \int_0^{t-s} \Phi^{-1}\left(\frac{4C}{u^2}\right)p(du),
$$

we obtain

$$
|R_{ts}|_{\beta,\alpha} \leq 8M^2 \int_0^{t-s} \Phi^{-1}\left(\frac{4U}{u^2}\right)p(du) + 9M \int_0^{t-s} \Phi^{-1}\left(\frac{4C}{u^2}\right)p(du).
$$

Now it is not difficult to extend this to generic $0 < s < t < T$. $\square$

**4. Young theory.** We are now ready to analyze the Young integration in the evolution setting along the same lines as in Section 2.3: we will first define the integral $\mathcal{J}(\hat{d}x\, z)$ for two Young paths $x, z$ in an abstract setting. Then we will solve Young SPDEs, and eventually, check our main assumptions in the fractional heat equation setting of Section 3.4.

4.1. *Young integration.* The extension of the notion of integral weighted by an analytical semigroup will be performed through the following algorithm, which will be used in fact throughout the remainder of the paper:

(1) Assume first that $x$ is a regular operator-valued increment, and $z$ a regular $\mathcal{B}$-valued function and let $\mathcal{J}_{ts}(\hat{d}x\, z) \equiv \int_s^t S_{tu}\, dx_u\, z_u$, for $(t,s) \in \mathcal{S}_2$, as an element of $\hat{\mathcal{C}}_2$.
(2) Through the application of $\hat{\delta}$ and $\hat{\Lambda}$, try to get an expression for $\mathcal{J}(\hat{d}x\, z)$ which depends only on minimal regularity requirements for $x$ and $z$.
(3) Extend the notion of integral using the previous step, and see that it induces the convergence of some well-chosen Riemann sums.

Here is how this general strategy can be implemented here: suppose for the moment that $x$ is a smooth operator-valued function and $z$ a smooth function. Then it is easily checked that

$$\mathcal{J}(\hat{d}x\, z) = \mathcal{J}(\hat{d}x\, S)z + \mathcal{J}(\hat{d}x\, \hat{\delta}z). \tag{48}$$

Note that in this last equality appears for the first time an incremental operator which will play a fundamental role in the sequel, namely the operator $X^1 \in \hat{\mathcal{C}}_2(\mathcal{L}(\mathcal{B}))$ defined by

$$X^1_{ts} = \mathcal{J}_{ts}(\hat{d}x\, S) = \int_s^t S_{tu}\, dx_u\, S_{us}. \tag{49}$$

And here is an important point of our strategy: the noise $x$ does not appear by itself but always inside a convolution of the form (49), so its action is *milded* by the regularizing properties of the semigroup.

Applying $\hat{\delta}$ to the last term of equation (48) and invoking Proposition 3.7, we get

$$\hat{\delta}[\mathcal{J}(\hat{d}x\, \hat{\delta}z)] = \mathcal{J}(\hat{d}x\, S)\, \hat{\delta}z = X^1\, \hat{\delta}z.$$

If the 2-increment $X^1\, \hat{\delta}z$ is small enough, namely if $X^1\, \hat{\delta}z \in \hat{\mathcal{C}}_2^{\mu,\theta}$ for some $\theta$ and some $\mu > 1$, then we can express $\mathcal{J}(\hat{d}x\, z)$ as

$$\mathcal{J}(\hat{d}x\, z) = X^1 z + \hat{\Lambda}[X^1\, \hat{\delta}z] = (\mathrm{Id} - \hat{\Lambda}\hat{\delta})[X^1 z]. \tag{50}$$

The last equality is justified by noting that when $x$ is a smooth incremental operator, we have $\hat{\delta}X^1 = X^1 a$ (i.e. $\tilde{\delta}X^1 = 0$), and thus by Lemma 3.2 one obtains that $\hat{\delta}(X^1 z) = -X^1\hat{\delta}z$.

Let us turn now to the second point of our general strategy, which consists in inverting the process which leads to (50): indeed, if we can define properly the right-hand side of (50), then we will be able to extend the notion of integral by a procedure which is coherent with the basic properties required to any integral $\mathcal{J}$. Notice that this step only relies on the definition of an operator $X^1$ associated to $x$, satisfying $\hat{\delta} X^1 = 0$, and such that $X^1$ is regular enough. This will be formalized in the following theorem (recall that the space $\mathcal{B}_{-\infty}$ has been defined at Section 3.1).

THEOREM 4.1. *Let then $x$ be a path from $[0,T]$ to $\mathcal{B}_{-\infty}$ such that the operator $X^1$ associated to $x$ is well defined as an element of $\hat{\mathcal{C}}_2^\kappa \mathcal{L}^{\beta,\alpha}$, where $\beta, \kappa$ are positive constants, and $\alpha \in \mathbb{R}$. We also assume that $X^1$ satisfies the algebraic relation $\hat{\delta} X^1 = X^1 a$. Let $z \in \hat{\mathcal{C}}_1^{\eta,\beta}$, with $\kappa + \eta > 1$, and set*

$$\mathcal{J}(\hat{d}x\, z) = X^1 z + \hat{\Lambda}[X^1 \hat{\delta} z] = (\mathrm{Id} - \hat{\Lambda}\hat{\delta})[X^1 z]. \tag{51}$$

*Then*

(1) *$\mathcal{J}(\hat{d}x\, z)$ is well defined as an element of $\mathcal{E}_2^{\kappa,\alpha}$.*
(2) *For a constant $c > 0$, we have*

$$\|\mathcal{J}(\hat{d}x\, z)\|_{\kappa,\alpha} \le \|X^1\|_{\kappa,\beta,\alpha}(\|z\|_{0,\beta} + c_\mu \|z\|_{\eta,\beta}),$$

*where the norm $\|\cdot\|_{\kappa,\beta,\alpha}$ has been defined at relation (39).*
(3) *It holds that, for any $0 \le s < t \le T$*

$$\mathcal{J}_{ts}(\hat{d}x\, z) = \lim_{|\Pi_{ts}| \to 0} \sum_{i=0}^{n} S_{t-t_{i+1}} X^1_{t_{i+1},t_i} z_{t_i},$$

*where the limit is over all partitions $\Pi_{ts} = \{t_0 = t, \dots, t_n = s\}$ of $[s,t]$ as the mesh of the partition goes to zero.*

PROOF. Since $X^1 z$ is a well defined element of $\hat{\mathcal{C}}_2^{\kappa,\alpha}$, in order to show that the r.h.s. of equation (51) is well defined, it only remains to check that $X^1 \hat{\delta} z$ is in the domain of $\hat{\Lambda}$. However, since we have assumed that $X^1 \in \hat{\mathcal{C}}_2^\kappa \mathcal{L}^{\beta,\alpha}$ and $\hat{\delta} z \in \hat{\mathcal{C}}_2^{\eta,\beta}$, we obviously get that $X^1 \hat{\delta} z \in \hat{\mathcal{C}}_3^{\kappa+\eta,\alpha}$. Thus, according to Theorem 3.5, $X^1 \hat{\delta} z \in \mathrm{Dom}(\hat{\Lambda})$, yielding the first assertion of our theorem.

Moreover, thanks to the second part of Theorem 3.5, we have

$$\|\hat{\Lambda}[X^1 \hat{\delta} z]\|_{\mu,\alpha} \le c_\mu \|X^1\|_{\kappa,\beta,\alpha} \|z\|_{\eta,\beta}, \tag{52}$$

and it is also readily checked that

$$\|X^1 z\|_{\kappa,\alpha} \le \|X^1\|_{\kappa,\beta,\alpha} \|z\|_{0,\beta}, \tag{53}$$

which shows our second claim, by using equation (52) and equation (53) to estimate the r.h.s. of (51). Eventually, the third part of the theorem is a direct consequence of Corollary 3.6. □

REMARK 4.2. It is worth stressing at this point some elementary properties enjoyed by the extension of the notion of integral given by Theorem 4.1:

- The third part of the theorem states that $\mathcal{J}_{ts}(\hat{d}x\, z)$ is associated to some natural Riemann sums involving the processes $x$ (through $X^1$) and $z$.
- The arguments leading to relation (50) also show that, in case of some smooth processes $x$ and $z$, our integral $\mathcal{J}_{ts}(\hat{d}x\, z)$ coincides with the usual one.

These first properties seem to imply that our integral extension is a reasonable one.

4.2. *Young SPDEs.* Recall that we wish to solve an equation of the form

$$dy_t = Ay_t\, dt + dx_t\, f(y_t), \qquad t \in [0,T], \tag{54}$$

with an initial condition $y_0 = \psi \in \mathcal{B}_\kappa$, where $x$ is an operator-valued process which represents our noise and $f:\mathcal{B} \to \mathcal{B}$ is a (possibly) nonlinear regular map. As mentioned in the Introduction, we will consider equation (54) in the mild sense, that is, we will say that $y$ is a solution to (54) if, for a given $\kappa > 0$ (specified below) we have $y \in \mathcal{C}_1^{\kappa,\kappa}$ and if, for any $t \in [0,T]$, $y_t$ satisfies

$$y_t = S_t\psi + \int_0^t S_{tu}\, dx_u f(y_u) = S_t y_0 + \mathcal{J}_{t0}(\hat{d}x\, f(y)), \tag{55}$$

where the integral $\mathcal{J}_{t0}(\hat{d}x\, f(y))$ is understood in the sense of Theorem 4.1. In fact, we will focus here on a slight extension of the problem given by (55): we will search for a (unique) process $y \in \mathcal{C}_1^{\kappa,\kappa}$ satisfying, for any $(t,s) \in \mathcal{S}_2$,

$$y_t = S_{ts} y_s + \mathcal{J}_{ts}(\hat{d}x\, f(y)), \qquad y_0 = \psi, \tag{56}$$

from which one recovers obviously (55) by taking $s = 0$. Now, (56) can be expressed in terms of convolution increments, since it is equivalent to the following one:

$$\begin{aligned} &[\hat{\delta} y]_{ts} = \mathcal{J}_{ts}(\hat{d}x\, f(y)) = [(\mathrm{Id} - \hat{\Lambda}\hat{\delta})[X^1 f(y)]]_{ts} \qquad \text{for } (t,s) \in \mathcal{S}_2 \quad \text{and} \\ &\quad y_0 = \psi, \end{aligned} \tag{57}$$

which sticks better to the algebraic formalism introduced in the previous sections.

Let us specify also some of the assumptions under which our computations will be performed: first of all, the incremental operator $X^1$ defined by (49) will be assumed to be in the following class.

Hypothesis 1. *Assume that* $X^1 \in \hat{\mathcal{C}}_2^{\tilde{\gamma}} \mathcal{L}^{0,-\kappa} \cap \hat{\mathcal{C}}_2^{\kappa_0} \mathcal{L}^{\kappa,\kappa}$ *for some* $\tilde{\gamma} > \kappa_0 > \kappa > 1/4$ *such that*

$$\tilde{\gamma} + \kappa > 1, \qquad \tilde{\gamma} - \kappa \geq \kappa_0, \qquad \kappa < 1/2.$$

Notice that in the hypothesis above, the condition $\kappa < 1/2$ is somehow redundant. Indeed, if $\tilde{\gamma} \geq \kappa + \kappa_0 \geq 2\kappa$, this forces the relation $\kappa < 1/2$.

As far as the function $f$ is concerned, we will also assume that the following holds true.

Hypothesis 2. *Let* $\kappa$ *be the strictly positive constant defined at Hypothesis 1. We assume that the function* $f : \mathcal{B}_\kappa \to \mathcal{B}_\kappa$ *is locally Lipschitz, and satisfies* $|f(x)|_{\mathcal{B}_\kappa} \leq c_f(1 + |x|_{\mathcal{B}_\kappa})$. *Furthermore, we suppose that* $f$ *can also be seen as a map from* $\mathcal{B}$ *to* $\mathcal{B}$, *and when considered as such, it holds that* $f$ *is globally Lipschitz.*

With these assumptions and notation in mind, we are now able to solve our evolution equation in the Young sense.

Theorem 4.3. *Assume Hypotheses 1 and 2 hold true, and that* $\psi \in \mathcal{B}_\kappa$. *Let* $\hat{\mathcal{C}}_1^{*,\kappa}$ *be the subspace of* $\hat{\mathcal{C}}_1$ *defined by the norm*

$$\|z\|_{*,\kappa} = \|z\|_{0,\kappa} + \|\hat{\delta} z\|_{\kappa,\kappa}. \tag{58}$$

*Then there exists a unique global solution to (57) in* $\hat{\mathcal{C}}_1^{*,\kappa}$. *Furthermore, this solution enjoys the following properties:*

(a) *For any* $t \in [0,T]$, $y_t$ *can be written as* $y_t = S_t \psi + (\hat{\delta} y)_{t0}$.
(b) *Let us call* $\Phi$ *the map* $(\psi, X^1) \mapsto \Phi(\psi, X^1) = y$, *where* $y$ *is the solution to (57). Then* $\Phi$ *is Lipschitz continuous from* $\mathcal{B}_{\kappa_0} \times (\hat{\mathcal{C}}_2^{\tilde{\gamma}} \mathcal{L}^{0,-\kappa} \cap \hat{\mathcal{C}}_2^{\kappa_0} \mathcal{L}^{\kappa,\kappa})$ *to* $\hat{\mathcal{C}}_1^{*,\kappa}$.

Proof. A classical fixed point argument will be sufficient to obtain the global solution. Let us introduce the map $\Gamma : \hat{\mathcal{C}}_1^{*,\kappa} \to \hat{\mathcal{C}}_1^{*,\kappa}$ defined in the following way: if $y \in \hat{\mathcal{C}}_1^{*,\kappa}$, we set $\Gamma(y) = z$, where $z$ satisfies

$$\begin{aligned} [\hat{\delta} z]_{ts} &= \mathcal{J}_{ts}(\hat{d}x\, f(y)) \\ &= X^1 f(y) + [\hat{\Lambda}[X^1 \,\delta f(y)]]_{ts} \qquad \text{for } (t,s) \in \mathcal{S}_2 \quad \text{and} \\ z_0 &= \psi. \end{aligned} \tag{59}$$

Let also $B$ be the ball defined by

$$B = \{y; y_0 = \psi, \|y\|_{*,\kappa} \leq 2(1 + |\psi|_{\mathcal{B}_\kappa})\}. \tag{60}$$

Then the fixed point argument can be decomposed into two usual steps:

(1) Show that, on a small enough interval $[0,T]$, the ball $B$ is left invariant by $\Gamma$.
(2) Prove that $\Gamma$, restricted to the ball $B$, is a contraction.

We will mainly focus, in this proof, on the first of these steps, since it contains most of the technical difficulties associated to our method.

Take $y \in B$, and let us show that $z = \Gamma(y) \in B$ whenever $T$ is small enough (recall that $\mathcal{S}_2$ depends on the parameter $T$). To this purpose, we will first bound the term $\hat{\Lambda}[X^1\,\hat{\delta} f(y)]$ in (59). Recall that

$$[\hat{\delta} f(y)]_{ts} = [\delta f(y)]_{ts} - a_{ts} f(y_s), \tag{61}$$

and let us estimate the terms in the right-hand side of (61) separately: on one hand, recalling the notation of Section 3.3.1, and thanks to the fact that $f$ is Lipschitz on $\mathcal{B}$, we have

$$\begin{aligned} |[\delta f(y)]_{ts}|_{\mathcal{B}} &\le c_f |[\delta y]_{ts}|_{\mathcal{B}} \le c_f(|[\hat{\delta} y]_{ts}|_{\mathcal{B}} + |a_{ts} y_s|_{\mathcal{B}}) \\ &\le c_f[\|y\|_{\kappa,0} + \|y\|_{0,\kappa}]|t-s|^{\kappa} \\ &\le c_f \|y\|_{*,\kappa}|t-s|^{\kappa}, \end{aligned} \tag{62}$$

where $c_f$ is a positive constant which may change from line to line, but which depends only on $f$. On the other hand, according to Hypothesis 2, it is readily checked that $|f(y_s)|_{\mathcal{B}_\kappa} \le c_f(1 + \|y\|_{0,\kappa})$. Thus, invoking (32), we obtain that

$$|a_{ts} f(y_s)|_{\mathcal{B}} \le c_f(1 + \|y\|_{*,\kappa})|t-s|^{\kappa}, \tag{63}$$

and plugging (63) and (62) into (61), we get

$$|[\hat{\delta} f(y)]_{ts}|_{\mathcal{B}} \le c_f(1 + \|y\|_{*,\kappa})|t-s|^{\kappa}. \tag{64}$$

However, we know that $X^1 \in \hat{\mathcal{C}}_2^{\tilde{\gamma}} \mathcal{L}^{0,-\kappa}$, and this fact, together with the last estimate, yields

$$|[X^1\,\hat{\delta} f(y)]_{tus}|_{\mathcal{B}_{-\kappa}} \le c_f \|X^1\|_{\tilde{\gamma},0,-\kappa}(1 + \|y\|_{*,\kappa})|t-u|^{\tilde{\gamma}}|u-s|^{\kappa}.$$

Furthermore, by Hypothesis 1, we have $\tilde{\gamma} + \kappa > 1$. This means that Theorem 3.5 can be applied here to obtain that $\hat{\Lambda}(X^1\,\hat{\delta} f(y)) \in \mathcal{E}_2^{\tilde{\gamma}+\kappa,-\kappa}$. In particular, invoking the definition (38) of the space $\mathcal{E}_2^{\tilde{\gamma}+\kappa,-\kappa}$, and since $2\kappa < 1$ and $\kappa_0 < \tilde{\gamma} - \kappa$, we get $\hat{\Lambda}(X^1\,\hat{\delta} f) \in \hat{\mathcal{C}}_2^{\tilde{\gamma}-\kappa,\kappa} \subseteq \hat{\mathcal{C}}_2^{\kappa_0,\kappa}$. Moreover,

$$\|\hat{\Lambda}(X^1\,\hat{\delta} f)\|_{\kappa_0,\kappa} \le c_{f,\tilde{\gamma},\kappa}\|X^1\|_{\tilde{\gamma},0,-\kappa}(1 + \|y\|_{*,\kappa}). \tag{65}$$

A bound similar to equation (65) can be found for the term $X^1 f(y)$ appearing in the definition of $\hat{\delta} z$ in equation (59). Indeed, owing to the fact that $X^1 \in \hat{\mathcal{C}}_2^{\kappa_0}\mathcal{L}^{\kappa,\kappa}$ and that $f$ has linear growth in $\mathcal{B}_\kappa$, we get

$$\begin{aligned} |X^1_{ts} f(y_s)|_{\mathcal{B}_\kappa} &\le c_f \|X^1\|_{\kappa_0,\kappa,\kappa}(1 + |y_s|_{\mathcal{B}_\kappa})|t-s|^{\kappa_0} \\ &\le c_f \|X^1\|_{\kappa_0,\kappa,\kappa}(1 + \|y\|_{*,\kappa})|t-s|^{\kappa_0}. \end{aligned} \tag{66}$$

Hence, plugging (66) and (65) into (59), one obtains that $\|\hat{\delta}z\|_{\kappa_0,\kappa} \le c_{f,X^1,\tilde{\gamma},\kappa}(1+\|y\|_{*,\kappa})$. Note here a crucial point: starting from $y \in \hat{\mathcal{C}}_1^{\kappa,\kappa}$, we have constructed $z \in \hat{\mathcal{C}}_1^{\kappa_0,\kappa}$ with $\varepsilon = \kappa_0 - \kappa > 0$. This little regularity gain can be used in order to write

$$\|\hat{\delta}z\|_{\kappa,\kappa} \le c_{f,X^1,\tilde{\gamma},\kappa}(1+\|y\|_{*,\kappa})T^{\varepsilon}. \tag{67}$$

Now, the quantity $T^{\varepsilon}$ can be made arbitrarily small as $T \to 0$. Moreover, recall that we still need a bound on $\|z\|_{*,\kappa}$ defined by (58), and thus an estimate on $\|z\|_{0,\kappa}$ is needed at this point. However, it is easily checked that

$$|z_t|_{\mathcal{B}_\kappa} \le |S_t\psi|_{\mathcal{B}_\kappa} + |(\hat{\delta}z)_{t0}|_{\mathcal{B}_\kappa} \le |\psi|_{\mathcal{B}_\kappa} + T^{\kappa_0}\|\hat{\delta}z\|_{\kappa_0,\kappa}. \tag{68}$$

Putting together (67) and (68), we finally get, on $[0,T]$, that

$$\|z\|_{*,\kappa} \le |\psi|_{\mathcal{B}_\kappa} + c(1+\|y\|_{*,\kappa})T^{\varepsilon} \qquad \text{with } c = c_{f,X^1,\tilde{\gamma},\kappa},$$

which yields that, whenever $cT^{\varepsilon} \le 1/2$, the ball $B$ defined by (60) is left invariant by the map $\Gamma$.

Now that the invariance of $B$ has been shown, the contraction property for $\Gamma$ in a small interval $[0,T]$ is a matter of standard arguments, and is left to the reader for sake of conciseness. Let us just mention that $f$ is only supposed to be locally Lipschitz when considered as a function from $\mathcal{B}_\kappa$ to $\mathcal{B}_\kappa$. However, we are able to establish the contraction property here, due to the fact that we are confined to the ball $B$. This gives the existence and uniqueness result for equation (57) in the small interval $[0,T]$ whose size does not depend on the initial condition $\psi$. The construction of a global unique solution from the solution in $[0,T]$ is also quite standard, and its proof will be omitted here. □

4.3. *Application: the fractional heat equation.* Let us see now how the abstract results of Section 4.2 can be applied in the case of the heat equation driven by a fractional Brownian motion defined at Section 3.4. Recall that this means that we wish to solve equation (55) in case $A = \Delta - \mathrm{Id}$, where $\Delta$ is the Laplace operator on the circle, $x$ is a fractional Brownian motion defined by the covariance function (43), $\mathcal{B}_\kappa$ stands for the usual Sobolev space on $[0,1]$, and $f : \mathcal{B}_\kappa \to \mathcal{B}_\kappa$ is defined by $[f(y)](\xi) = \sigma(y(\xi))$ for $\xi \in [0,1]$ and a smooth function $\sigma : \mathbb{R} \to \mathbb{R}$. In other words, we will try to solve the equation

$$y(t,\xi) = \int_0^1 G_t(\xi,\eta)\psi(\eta)\,d\eta + \int_0^t \int_0^1 G_{t-s}(\xi,\eta)X(ds,d\eta)\sigma(y_s(\eta)), \tag{69}$$

where the last integral has to be understood in the sense of Theorem 4.1. Notice that we have chosen here a multiparametric formulation for our equation, for computational purposes. However, as mentioned in the Introduction, this setting can be translated easily into the infinite-dimensional one.

Now, the application of Theorem 4.3 in this context amounts to define an incremental operator $X^1$ related to our problem, and then to show that Hypotheses 1 and 2 are fulfilled.

Let us give then a natural definition of the operator $X^1$ associated to our equation: we will set, for $\psi \in \mathcal{B}$ and $(t,s) \in \mathcal{S}_2$,

$$
\begin{aligned}
[X^1_{ts}\psi](\xi) &= [\mathcal{J}_{ts}(\hat{d}XS)]\psi(\xi) \\
&= \int_s^t \int_0^1 G_{t-u}(\xi,\eta_1) X(du, d\eta_1) \left( \int_0^1 G_{v-s}(\eta_1,\eta_2)\psi(\eta_2)\, d\eta_2 \right),
\end{aligned} \tag{70}
$$

which has to be understood now in the Wiener sense, as a centered Gaussian random variable whose variance is given by (43). In this context, the regularity result we obtain on $X^1$ is the following.

Proposition 4.4. *Let $X$ be an infinite-dimensional fractional Brownian motion defined by the covariance function (43) for a given $H > 1/2$, with $Q$ given by (44) for $\nu \in [0,1)$. Suppose that $H + \bar{\nu}/2 > 3/4$, with the convention $\bar{\nu} = \nu \wedge (1/2)$. Let $X^1$ be the incremental operator defined by (70). Then for any $\tilde{\gamma} < H - 1/4 + \bar{\nu}/2$, $\kappa \in (1/4, 1/2)$, $\kappa_0 = \tilde{\gamma} - \kappa$ and $\gamma < H$ we have*

$$X^1 \in \hat{\mathcal{C}}_2^{\tilde{\gamma}} \mathcal{L}_{\mathrm{HS}}^{0,-\kappa} \cap \hat{\mathcal{C}}_2^{\kappa_0} \mathcal{L}_{\mathrm{HS}}^{\kappa,\kappa} \cap \hat{\mathcal{C}}_2^{\gamma} \mathcal{L}_{\mathrm{HS}}^{\kappa,-\kappa},$$

*almost surely.*

Remark 4.5. The reader will probably notice that the assumption $\kappa > 1/4$ is not necessary in order to prove the proposition above. However, we include it already at this stage, since this restriction is crucial for Proposition 4.10 to hold true.

The proof of Proposition 4.4 relies on the following elementary lemmas, that we label for further use.

Lemma 4.6. *For any $\alpha < \beta$, such that $\alpha + \beta > 1/2$, there exists a constant $C$ such that*

$$\sum_{i,j: i+j=k} \lambda_i^{-\alpha} \lambda_j^{-\beta} \le C \lambda_k^{-\alpha - \bar{\beta} + 1/2},$$

*where $\bar{\beta} = \min(\beta; 1/2)$.*

Lemma 4.7. *Let $a$ and $b$ be two positive constants, and $H > 1/2$. Then the integral*

$$\int_0^1 \int_0^1 |u-v|^{2H-2} |2-u-v|^{-a} |u+v|^{-b}\, du\, dv$$

*is finite whenever $2H - a > 0$ and $2H - b > 0$.*

We leave the easy proof of these results to the reader.

Proof of Proposition 4.4. We need to prove that the r.v. $X^1$ has a version with the claimed regularity. For random operators, up to our knowledge, no standard method is available to prove regularity properties. So we have chosen the following simple (though arguably nonoptimal) strategy in order to obtain a regular version: first, we determine the kernel associated to the operator $X^1$, then using the kernel we estimate its Hilbert–Schmidt norm in some $L^2$ space. This will be enough to apply the modified Garsia–Rodemich–Rumsey Lemma 3.8 and conclude the proof. We will develop now this strategy into several steps, discussing in detail the proof of $X^1 \in \hat{\mathcal{C}}_2^{\tilde{\gamma}} \mathcal{L}_{\mathrm{HS}}^{0,-\kappa}$. The other pathwise statements can be proven similarly.

*Step 1: Definition of a random kernel.* For $(t,s) \in \mathcal{S}_2$, $X^1_{ts}$ is considered as an operator from $\mathcal{B} = L^2([0,1])$ to $\mathcal{B}_{-\kappa}$, and thus $\|X^1_{ts}\|_{\mathrm{HS},\mathcal{B}\to\mathcal{B}_{-\kappa}} = \|A_o^{-\kappa} X^1_{ts}\|_{\mathrm{HS},\mathcal{B}\to\mathcal{B}}$, which is the expression we are going to evaluate. Pick $\psi \in \mathcal{B}$ smooth enough. Applying Fubini's theorem for the fractional Brownian motion, we get

$$\begin{aligned}
[A_o^{-\kappa} X^1_{ts}\psi](\xi) &= A_o^{-\kappa} \int_s^t \int_0^1 G_{t-u}(\xi,\eta_1) X(du, d\eta_1) \\
&\qquad \times \left( \int_0^1 G_{u-s}(\eta_1,\eta_2)\psi(\eta_2)\, d\eta_2 \right) \\
&= \int_s^t \int_0^1 G^{-\kappa}_{t-u}(\xi,\eta_1) X(du, d\eta_1) \left( \int_0^1 G_{u-s}(\eta_1,\eta_2)\psi(\eta_2)\, d\eta_2 \right) \\
&= \int_0^1 K_{ts}(\xi,\eta_2)\psi(\eta_2)\, d\eta_2,
\end{aligned}$$

where the kernel $G^{-\kappa}_{t-u}$ has been defined at Section 3.4, and where $K_{ts}(\xi,\eta)$ is the random kernel on $[0,1]^2$ defined by the Wiener integral

$$K_{ts}(\xi,\eta) = \int_s^t \int_0^1 G^{-\kappa}_{t-u}(\xi,\eta_1) G_{u-s}(\eta_1,\eta) X(du, d\eta_1).$$

Hence, the Hilbert–Schmidt norm of $X^1_{ts}$, seen as an operator from $\mathcal{B}$ to $\mathcal{B}_{-\kappa}$, will be given by

$$\|X^1_{ts}\|^2_{\mathrm{HS}} = \int_0^1 \int_0^1 [K_{ts}(\xi,\eta)]^2\, d\xi\, d\eta. \tag{71}$$

Our next aim will then be to evaluate this last quantity.

*Step 2: $L^2$ computations.* A direct application of (43) gives

$$E[K^2_{ts}(\xi,\eta)]$$

$$= c_H \int_s^t \int_s^t \Big( \int_{[0,1]^2} G_{t-u}^{-\kappa}(\xi, z) G_{u-s}(\eta, z) Q(z - \hat{z}) \times G_{t-v}^{-\kappa}(\xi, \hat{z}) G_{v-s}(\eta, \hat{z})\, dz\, d\hat{z} \Big) |u - v|^{2H-2}\, du\, dv.$$

Furthermore, for $z, \hat{z} \in [0,1]$, it holds that

$$\int_0^1 G_{t-u}^{-\kappa}(\xi, z) G_{t-v}^{-\kappa}(\xi, \hat{z})\, d\xi = G_{2t-u-v}^{-2\kappa}(z, \hat{z}),$$

$$\int_0^1 G_{u-s}(\eta, z) G_{v-s}(\eta, \hat{z})\, d\eta = G_{u+v-2s}(z, \hat{z}).$$

Thus, going back to relation (71), we obtain

$$\begin{aligned} A_{ts} &\equiv E[\|X_{ts}^1\|_{\mathrm{HS}}^2] \\ &= c_H \int_s^t \int_s^t \Big( \int_{[0,1]^2} Q(z - \hat{z}) G_{u+v-2s}(z, \hat{z}) G_{2t-u-v}^{-2\kappa}(z, \hat{z})\, dz\, d\hat{z} \Big) \\ &\qquad \times |u - v|^{2H-2}\, du\, dv \\ &= c_H \int_0^\varepsilon \int_0^\varepsilon F(u, v) |u - v|^{2H-2}\, du\, dv, \end{aligned} \tag{72}$$

where we have set $\varepsilon = t - s$, and with $F : [0, \varepsilon]^2 \to \mathbb{R}_+$ defined by

$$F(u, v) = \int_{[0,1]^2} Q(z - \hat{z}) G_{u+v}(z, \hat{z}) G_{2\varepsilon-u-v}^{-2\kappa}(z, \hat{z})\, dz\, d\hat{z}. \tag{73}$$

Furthermore, plugging the definitions (42) and (44) into (73), and invoking the fact that $\{e_n; n \in \mathbb{Z}\}$ is an orthonormal basis of $L^2([0,1])$, we get

$$F(u, v) = \sum_{m,n,l \in D} \lambda_n^{-\nu} \lambda_l^{-2\kappa} e^{-\lambda_m(u+v)} e^{-\lambda_l(2\varepsilon-u-v)},$$

where $D = \{m, n, l \in \mathbb{Z}^3 : m + n + l = 0\}$. Then

$$A_{ts} = c_H \sum_{m,n,l \in D} \lambda_n^{-\nu} \lambda_l^{-2\kappa} \int_0^\varepsilon \int_0^\varepsilon \frac{e^{-\lambda_m(u+v)} e^{-\lambda_l(2\varepsilon-u-v)}}{|u - v|^{2-2H}}\, du\, dv.$$

Owing now to the fact that $x \mapsto x^a e^{-x}$ is a bounded function on $\mathbb{R}_+$ for any $a > 0$, we obtain, for a constant $c$ which may change from line to line,

$$\begin{aligned} A_{ts} &\le c \sum_{m,n,l \in D} \lambda_n^{-\nu} \lambda_l^{-2\kappa} \lambda_m^{-a} \int_0^\varepsilon \int_0^\varepsilon \frac{du\, dv}{|u - v|^{2-2H} (u + v)^a} \\ &\le c \varepsilon^{2H-a} \sum_{m+n+l=0} \lambda_n^{-\nu} \lambda_l^{-2\kappa} \lambda_m^{-a}, \end{aligned} \tag{74}$$

where we have used Lemma 4.7 under the condition $a<2H$. Let us now analyze the sum. Of course, we can write

$$\sum_{m+n+l=0} \lambda_n^{-\nu}\lambda_l^{-2\kappa}\lambda_m^{-a} \le \sum_{l,k:l+k=0} \lambda_l^{-2\kappa} \sum_{m,n:m+n=k} \lambda_m^{-a}\lambda_n^{-\bar\nu}.$$

Moreover, taking $a=1/2-\bar\nu+\eta$ for some small $\eta>0$ and using Lemma 4.6, we have

$$\sum_{m+n+l=0} \lambda_n^{-\nu}\lambda_l^{-2\kappa}\lambda_m^{-a} \le c \sum_{l,k:l+k=0} \lambda_l^{-2\kappa}\lambda_k^{-a-\bar\nu+1/2} = c \sum_{l,k:l+k=0} \lambda_l^{-2\kappa}\lambda_k^{-\eta},$$

and this sum is always finite under the condition $\kappa>1/4$. Then, going back to (74), we have found that $A_{ts}\le c\varepsilon^{2\tilde\gamma'}$, for any $\tilde\gamma'=H-a/2<H-1/4+\bar\nu/2$, where we recall that $\bar\nu=\inf(\nu;1/2)$.

*Step 3*: *$L^p$ estimates.* We will prove now that, for any $p\ge 1$, we have

$$E[\|X_{ts}^1\|_{\mathrm{HS}}^{2p}] \le c_p(t-s)^{2\tilde\gamma' p} \qquad \text{for } 0\le s<t\le T. \tag{75}$$

Indeed, a simple application of Hölder's inequality yields

$$\begin{aligned} E[\|X_{ts}^1\|_{\mathrm{HS}}^{2p}] &= \int_{[0,1]^{2p}} E\left[\prod_{i=1}^p K_{ts}^2(\xi_i,\eta_i)\right] d\xi_1\,d\eta_1\,\cdots\,d\xi_p\,d\eta_p \\ &\le \int_{[0,1]^{2p}} \prod_{i=1}^p E^{1/p}[K_{ts}^{2p}(\xi_i,\eta_i)]\,d\xi_1\,d\eta_1\,\cdots\,d\xi_p\,d\eta_p, \end{aligned}$$

and since $K_{ts}(\eta_i,\eta_i)$ is a Gaussian variable, we get

$$\begin{aligned} E[\|X_{ts}^1\|_{\mathrm{HS}}^{2p}] &\le c_p \int_{[0,1]^{2p}} \prod_{i=1}^p E[K_{ts}^2(\xi_i,\eta_i)]\,d\xi_1\,d\eta_1\,\cdots\,d\xi_p\,d\eta_p \\ &= c_p\left(\int_{[0,1]} E[K_{ts}^2(\xi,\eta)]\,d\xi\,d\eta\right)^p \\ &= c_p E^p[\|X_{ts}^1\|_{\mathrm{HS}}^2], \end{aligned}$$

which easily yields (75).

*Step 4*: *Conclusion.* Recall that $X^1$ is considered as an element of $\hat{\mathcal{C}}_2(\mathcal{L}_{\mathrm{HS}}^{0,-\kappa})$. We can use now inequality (47), which can be read here as

$$\|X^1\|_{\tilde\gamma,0,-\kappa} \le C[U_{\tilde\gamma+2/p,p,0,-\kappa}(X^1) + \|\tilde\delta X^1\|_{\tilde\gamma,0,-\kappa}], \tag{76}$$

in order to bound $\|X^1\|_{\tilde\gamma,0,-\kappa}$ for any $\tilde\gamma<\tilde\gamma'<H-1/4+\bar\nu/2$. Indeed, if $p$ is large enough, we have that $\tilde\gamma+2/p<\tilde\gamma'$, and the term $U_{\tilde\gamma+2/p,p,0,-\kappa}(X^1)$ is easily handled thanks to (75). This yields

$$E[U_{\tilde\gamma',p,0,-\kappa}(X^1)]<\infty. \tag{77}$$

We are now left with the estimation of $\|\tilde{\delta}X^1\|_{\tilde{\gamma}}$. However, remember that $\tilde{\delta}X^1 = 0$ in case of a regular signal $x$, and it is readily checked that this relation is still valid in the current fractional Brownian setting, so this term is identically zero. Thus, we have obtained that

$$E[\|X^1\|_{\tilde{\gamma},0,-\kappa}] \le cE[U_{\tilde{\gamma}',p,0,-\kappa}(X^1)] < \infty,$$

which implies that $\|X^1\|_{\tilde{\gamma},0,-\kappa} < \infty$ almost surely, concluding the proof.

Along the same lines as in the preceding steps, some $L^2$ bounds state that

$$E[\|X^1_{ts}\|^2_{\mathrm{HS},\mathcal{L}(\mathcal{B}_\kappa,\mathcal{B}_\kappa)}] \le c(t-s)^{2\kappa_0'} \qquad \text{for } 0 \le s < t \le T \tag{78}$$

and

$$E[\|X^1_{ts}\|^2_{\mathrm{HS},\mathcal{L}(\mathcal{B}_{-\kappa},\mathcal{B}_\kappa)}] \le c(t-s)^{2\gamma'} \qquad \text{for } 0 \le s < t \le T, \tag{79}$$

for any $\kappa_0' < H - 1/4 - \kappa + \bar{\nu}/2$ and $\gamma' < H$, respectively. Following the same strategy as before, these bounds are enough to prove the remaining assertions of the proposition. $\square$

Let us see now how this results can be related to our Hypothesis 1. Recall that the restriction $\kappa > 1/4$ is dictated by the fact that we need to work in a space $\mathcal{B}_\kappa$ embedded in the space $C([0,1])$ of continuous functions on $[0,1]$ in order to prove Proposition 4.10 below.

COROLLARY 4.8. *Suppose $X$ is an infinite-dimensional fractional Brownian motion defined by the covariance function (43) for a given $H > 1/2$, with $Q$ given by (44) for $\nu \ge 0$. Assume moreover that $H > 7/8 - \bar{\nu}/2$. Then the incremental operator $X^1$ satisfies Hypothesis 1 for some*

$$\kappa \in (1/4, 1/2), \qquad \kappa_0 < H - 1/4 - \kappa + \bar{\nu}/2, \qquad \tilde{\gamma} = \kappa_0 + \kappa.$$

PROOF. By the previous result, we have that $X^1$ has the required regularity for any $1/4 < \kappa < 1/2$, $\kappa_0 < H - 1/4 - \kappa + \bar{\nu}/2$ and $\tilde{\gamma} = \kappa_0 + \kappa < H - 1/4 + \bar{\nu}/2$. In order to check Hypothesis 1, we now need to require that $\tilde{\gamma} + \kappa > 1$. In fact, there exists $1/4 < \kappa < \kappa_0$ satisfying this inequality if and only if $\tilde{\gamma} + \kappa_0 > 1$, that is, $2H - 1/2 - \kappa + \bar{\nu} > 1$. This is equivalent to assume

$$H > 3/4 + \kappa/2 - \bar{\nu}/2 > 7/8 - \bar{\nu}/2.$$

In this latter case, it is easily seen that there exist $\tilde{\gamma}, \kappa, \kappa_0$ satisfying our requirements. $\square$

REMARK 4.9. If we are only interested in obtaining a local solution for our Young PDE, then the estimate (64) can be replaced by a bound in $\mathcal{B}_\kappa$, which will be quadratic in $y$. Hence, using the fact that

$$X^1 \in \hat{\mathcal{C}}_2^{\gamma}\mathcal{L}^{\kappa,-\kappa} \cap \hat{\mathcal{C}}_2^{\kappa_0}\mathcal{L}^{\kappa,\kappa}$$

for any $\gamma < H$ and $1/4 < \kappa < \kappa_0 < H - 1/4 - \kappa + \bar{\nu}/2$, the condition for the construction of the (local) fix-point map $\Gamma$ becomes $\gamma + \kappa > 1$. To fulfill this requirement with our fractional Brownian noise, we only have to impose $H > 3/4 - \bar{\nu}/4$. This condition is comparable (but a bit worse) with the results of [12], where the Hilbert spaces $W^{2\kappa,2}$ were considered, and where we found $H > 3/4 - \bar{\nu}/2$. One of the drawback of the approach presented in this paper is that the esimation of the random operators like $X^1$ in Banach spaces $W^{2\alpha,p}$ for $p > 2$ seems very difficult. Moreover, it seems that the estimation in the Hilbert–Schmidt norm causes another small loss of regularity, which means that even in the case of a "regular" noise $\nu = 1/2$, our bound on $H$ is $H > 5/8$ and not $H > 1/2$ as should be natural to expect and found in [12]. On the other hand, as we will see later, the operator approach seems better suited than the approach of [12] for a true rough-path expansion of SPDEs.

Now that we have checked the assumptions on $X^1$, let us turn to the hypothesis on the nonlinear coefficient $\sigma$ in equation (69). In order to deal with the Sobolev norms, it is worth mentioning that, instead of working with the spaces $\mathcal{B}_\kappa = H_\kappa$ we have used so far, characterized by their Fourier decomposition, we will consider the Sobolev spaces $W^{2\kappa,2}$, induced by the norms

$$[W_\kappa(\psi)]^2 \equiv \|\psi\|^2_{L^2([0,1])} + \int_{[0,1]^2} \frac{|\psi(\xi) - \psi(\eta)|^2}{|\xi - \eta|^{1+4\kappa}}\, d\xi\, d\eta. \tag{80}$$

These spaces are obviously more convenient than the spaces $\mathcal{B}_\kappa$ for the computations on $f$, and they are closely related to these latter spaces, since the following classical relation holds true (see [1]):

$$\mathcal{B}_{\kappa+\varepsilon} \subset W^{2\kappa,2} \subset \mathcal{B}_{\kappa-\varepsilon} \qquad \text{for any } \varepsilon > 0.$$

Using these embeddings, we can consider the operator $X^1_{ts}$ going from a space $W^{2\kappa,2}$ to a space $\mathcal{B}_\kappa$ by just loosing a little regularity in $t, s$. Then we can verify that $f$ satisfy a slight modification of Hypothesis 2.

PROPOSITION 4.10. *Let $\sigma \in C_b^2(\mathbb{R})$ be a real-valued function. Then, for any $\kappa > 1/4$, the function $f : W^{2\kappa,2} \to W^{2\kappa,2}$ defined by $[f(y)](\xi) = \sigma(y(\xi))$ is locally Lipschitz, satisfies $|f(x)|_{W^{2\kappa,2}} \le c_f(1 + |x|_{W^{2\kappa,2}})$ and is globally Lipschitz as a map $f : \mathcal{B} \to \mathcal{B}$.*

PROOF. Recall that, for our particular situation, $\mathcal{B} = L^2([0,1])$, and it is easily checked that, whenever $\sigma \in C_b^2(\mathbb{R})$, the function $f : \mathcal{B} \to \mathcal{B}$ is bounded and globally Lipschitz.

With these considerations in mind, it is readily seen that $f : W^{2\kappa,2} \to W^{2\kappa,2}$ has linear growth. In order to check that $f$ is also locally Lipschitz, note that its gradient can be computed as follows for $y, h \in W^{2\kappa,2}$:

$$\nabla f(y) : W^{2\kappa,2} \to W^{2\kappa,2}, \qquad [\nabla f(y) \cdot h](\xi) = \sigma'(y(\xi)) h(\xi).$$

Let us estimate now the norm (80) of $\nabla f(y) \cdot h$: first, if $\sigma'$ is a bounded function, then

$$\|\nabla f(y) \cdot h\|_{L^2([0,1])} \leq \|\sigma'\|_\infty \|h\|_{L^2([0,1])} \leq \|\sigma'\|_\infty \|h\|_{W^{2\kappa,2}}. \tag{81}$$

As far as the variational term of (80) is concerned, notice that we have assumed $\kappa > 1/4$, which means that $W^{2\kappa,2} \subset C([0,1])$, and for any $h \in W^{2\kappa,2}$, $\|h\|_\infty \leq c\|h\|_{W^{2\kappa,2}}$. Thus,

$$\begin{aligned} &\int_{[0,1]^2} \frac{|[\nabla f(y) \cdot h](\xi) - [\nabla f(y) \cdot h](\eta)|^2}{|\xi - \eta|^{1+4\kappa}} \, d\xi \, d\eta \\ &\quad \leq \|\sigma'\|_\infty^2 \int_{[0,1]^2} \frac{|h(\xi) - h(\eta)|^2}{|\xi - \eta|^{1+4\kappa}} \, d\xi \, d\eta \\ &\qquad + \|h\|_{W^{2\kappa,2}}^2 \int_{[0,1]^2} \frac{|\sigma'(y(\xi)) - \sigma'(y(\eta))|^2}{|\xi - \eta|^{1+4\kappa}} \, d\xi \, d\eta \\ &\quad \leq c_\sigma \|h\|_{W^{2\kappa,2}}^2 [1 + \|\sigma''\|_\infty^2 \|y\|_{W^{2\kappa,2}}^2]. \end{aligned} \tag{82}$$

Putting together (81) and (82), we have thus shown that

$$\|\nabla f(y)\|_{\mathcal{L}(W^{2\kappa,2})} \leq c_\sigma (1 + \|y\|_{W^{2\kappa,2}}),$$

which easily yields that $f : W^{2\kappa,2} \to W^{2\kappa,2}$ is locally Lipschitz. □

REMARK 4.11. Notice that, in spite of the fact that $\sigma$ is assumed to be a nicely behaved coefficient, its interpretation as an application from $W^{2\kappa,2}$ to $W^{2\kappa,2}$ does not enjoy the usual assumptions of boundedness made on coefficients in rough path theory (see, e.g., [13, 16, 17]). This is one of the major sources of problems in our computations, and in general in the extension of rough path theory to SPDEs.

Let us now summarize the considerations of the current section into the following theorem.

THEOREM 4.12. *Let $X$ be an infinite-dimensional fractional Brownian motion on $[0,T] \times [0,1]$, defined by the covariance function (43) and (44), with $H > 1/2$ and $\nu \in [0,1)$ such that $H > 7/8 - \bar{\nu}/2$ and let $\sigma \in C_b^2(\mathbb{R})$. Then, there exists $\kappa \in (1/4, 2H - 3/2 + \bar{\nu})$ such that for any initial condition $\psi \in \mathcal{B}_\kappa$, the equation*

$$Y(0,\xi) = \psi(\xi), \qquad \partial_t Y(t,\xi) = \Delta Y(t,\xi)\, dt + \sigma(Y(t,\xi)) X(dt, d\xi), \qquad t \in [0,T], \xi \in [0,1], \tag{83}$$

*with periodic boundary conditions, understood in the mild sense given by (57), has a unique global solution in $\hat{\mathcal{C}}_1^{\kappa,\kappa}$.*

PROOF. By Proposition 4.10, the map $f$ is Lipschitz and with linear growth from $\mathcal{B}_\kappa$ to $\mathcal{B}_{\kappa-\varepsilon}$ for arbitrarily small $\varepsilon$. As already noted this little mismatch of regularity can be compensated by the time-regularity of $X^1$. Then by a small modification of the arguments of Theorem 4.3 and by Proposition 4.4, we can directly solve the equation

$$\partial_t Y(t,\xi) = (\Delta - \mathrm{Id}) Y(t,\xi)\, dt + \sigma(Y(t,\xi)) X(dt, d\xi), \qquad t \in [0,T], \xi \in [0,1],$$

as a rough evolution equation in $\hat{\mathcal{C}}_1^{\kappa,\kappa}$. Now, if one wants to solve (83), it is sufficient to get an existence and uniqueness result for the equation

$$\partial_t Y(t,\xi) = (\Delta - \mathrm{Id}) Y(t,\xi)\, dt + Y(t,\xi)\, dt + \sigma(Y(t,\xi)) X(dt, d\xi), \qquad t \in [0,T], \xi \in [0,1],$$

which can be done along the same lines as for Theorem 4.3, by taking care of the additional drift term $Y dt$. This step is left to the reader. □

**5. Rough evolution equations: the linear case.** We pass now to the development of an expansion which allows to consider equation (54) in a case which goes beyond the Young theory, in terms of the Hölder regularity of the driving noise $x$. We start with a simple linear case, that is, $f \equiv \mathrm{Id}$, which will hopefully lead to a better understanding of our method.

5.1. *Strategy.* Recall that we wish to get some existence and uniqueness results for the equation

$$\hat{\delta} y_{ts} = \int_s^t S_{tu}\, dx_u\, y_u = \mathcal{J}_{ts}(\hat{d}x\, y) \quad \text{and} \quad y_0 = \psi. \tag{84}$$

Just like in the case of the Young integral, sketched at the beginning of Section 4.1, we will proceed as follows:

(1) Expand (84) as if $x$ were a regular process, until we get some terms which can be analyzed through the application of the operators $\hat{\delta}$ and $\hat{\Lambda}$.
(2) Define a natural extension of the notion of integral thanks to the first step, and show that this allows to integrate a reasonably wide class of functions.
(3) Solve the equation in the sense given by this notion of integral.

In the current section, we will mostly address the first of these three steps.

If $x$ is a regular process, equation (84) can be solved by means of the classical evolution theory. Furthermore, if $y$ designates the unique solution to (84), then according to our expansion strategy, $y$ also satisfies, for $t, s \in \mathcal{S}_2$,

$$\hat{\delta} y_{ts} = \int_s^t S_{tu}\, dx_u\, y_u = \int_s^t S_{tu}\, dx_u\, S_{us} y_s + \int_s^t S_{tu}\, dx_u\, \hat{\delta} y_{us}.$$

However, the last term of this equation cannot be defined by applying the map $\hat{\Lambda}$ when $x$ has low time regularity. In order to cope with this difficulty, let us expand again $\hat{\delta} y$ by plugging relation (84) into the previous equation. Doing this twice, we get

$$\begin{aligned}\hat{\delta} y_{ts} &= \int_s^t S_{tu}\, dx_u\, S_{us} y_s + \int_s^t S_{tu}\, dx_u \int_s^u S_{uv}\, dx_v\, y_v \\ &= \int_s^t S_{tu}\, dx_u\, S_{us} y_s + \int_s^t S_{tu}\, dx_u \int_s^u S_{uv}\, dx_v\, S_{vs} y_s \\ &\quad + \int_s^t S_{tu}\, dx_u \int_s^u S_{uv}\, dx_v \int_s^v S_{vw}\, dx_w S_{ws} y_s \\ &\quad + \int_s^t S_{tu}\, dx_u \int_s^u S_{uv}\, dx_v \int_s^v S_{vw}\, dx_w\, \hat{\delta} y_{ws}.\end{aligned}$$

Thus, going back to our notation on iterated integrals (41), we can recast (84) into

$$\hat{\delta} y = X^1 y + X^2 y + X^3 y + \mathcal{J}(\hat{d}x\, \hat{d}x\, \hat{d}x\, \hat{d}x\, y), \tag{85}$$

where, for $t, s \in \mathcal{S}_2$ and $\phi \in \mathcal{B}$, the operators $(X^i)_{i=1,2,3}$ are defined by

$$X^i_{ts}\phi := \mathcal{J}_{ts}(\hat{d}x\, X^{i-1})\phi = \int_s^t S_{tu}\, dx_u\, X^{i-1}_{us}\phi \tag{86}$$

with $X^0_{ts} = S_{ts}$. These operators are the new building block we will need in order to solve equation (84), and they play the role of the iterated integrals of rough path theory in our bilinear evolution context. Notice that the last term in equation (85) is considered as a remainder: suitable assumptions should be made to ensure that it will be *small* enough. Notice also that we stopped our expansion at the third order. We will see that this is the minimum order which allows to handle the Brownian case.

Let us say a few words now about the algebraic properties of the operators $X^i$: when $x$ is a smooth process, we have, for example,

$$\hat{\delta} X^2_{tus} = \int_u^t S_{tv}\, dx_v \int_s^v S_{vw}\, dx_w\, S_{ws} - \int_u^t S_{tv}\, dx_v \int_u^v S_{vw}\, dx_w\, S_{wu},$$

and using some elementary algebra, we end up with

$$\begin{aligned}\hat{\delta} X^2_{tus} &= \int_u^t S_{tv}\, dx_v\, S_{vu} \int_s^u S_{uw}\, dx_w\, S_{ws} + \int_u^t S_{tv}\, dx_v \int_u^v S_{vw}\, dx_w\, S_{wu}[S_{us} - \mathrm{Id}] \\ &= X^1_{tu} X^1_{us} + X^2_{tu} a_{us}.\end{aligned}$$

Thus, taking into account our algebraic convention (16) and the definition of $\tilde{\delta}$ given at (46), we have obtaind the relation $\tilde{\delta} X^2 = X^1 X^1$. In a more general way, it is not difficult to show by induction that

$$\tilde{\delta} X^n = \sum_{i=1}^{n-1} X^i X^{n-i},$$

which are exactly the Chen relations in this setting.

We can now specialize our previous program into the following:

2a. Assume that the operator-valued 1-increments $X^1, X^2, X^3$ are defined by some kind of operation which preserves the usual algebraic relations between integrals (e.g., use stochastic calculus with respect to an Hilbert space valued fractional Brownian motion or some other limiting procedure on discrete sums). They will be our (step-3) *rough path*.

2b. Using $(X^1, X^2, X^3)$ define an integration theory for a sufficiently large class of functions $\mathcal{Q}$ so that it will be possible to give a meaning to integrals of the form $z_t = \int_0^t S_{tu}\, dx_u\, y_u$ for any $y \in \mathcal{Q}$. We will call $\mathcal{Q}$ the space of paths *controlled by* $X$.

3′. Study the map $\Gamma : \mathcal{Q} \to \mathcal{Q}$ defined by $\Gamma(y)_t = \int_0^t S_{tu}\, dx_u\, y_u$, and prove that it has a fixed point $y = \Gamma(y)$ which will be then a solution of the evolution problem (84).

5.2. *Integration of weakly controlled paths.* We start by postulating some reasonable properties for $X^n$.

HYPOTHESIS 3. *We will assume that the process $x$ allows to define some operator-valued increments $X^1, X^2, X^3$, representing morally (49) and (86), respectively. This amounts for us to suppose that the $X^i$'s satisfy the algebraic relations*

$$\tilde{\delta} X^1 = 0, \qquad \tilde{\delta} X^2 = X^1 X^1, \qquad \tilde{\delta} X^3 = X^1 X^2 + X^2 X^1,$$

*and that the following Hölder-regularity properties holds true:*

$$X^i \in \hat{\mathcal{C}}_2^{\gamma+(i-1)\kappa_0} \mathcal{L}^{\eta,-\rho} \cap \hat{\mathcal{C}}_2^{i\kappa_0} \mathcal{L}^{\eta,\eta}, \qquad i = 1, 2, 3,$$

*for some $\eta, \rho \geq 0$ and $\gamma, \kappa_0$ such that $\gamma = \kappa_0 + \eta + \rho$ and $\gamma + 3\kappa_0 > 1$.*

We will define now the class $\mathcal{Q}$ of processes we wish to be able to integrate against $x$: in the current situation, it will include any process which can be decomposed into a part depending on $X^1, X^2$, plus a remainder term which is assumed to be *small* enough. For the sake of a contraction argument needed below (compare to the Young case), we fix a given time regularity $\kappa$ such that $0 < \kappa < \kappa_0$.

DEFINITION 5.1 (Weakly controlled paths). Let $\psi \in \mathcal{B}_\eta$ be a given initial condition. A path $y \in \hat{\mathcal{C}}_1^{\kappa,\eta}$ is said to be *weakly controlled* by $X^1, X^2$ if $y_0 = \psi$ and $\hat{\delta} y$ can be decomposed into

$$\hat{\delta} y = X^1 y^1 + X^2 y^2 + y^r, \qquad \hat{\delta} y^1 = X^1 y^2 + y^{1,r} \tag{87}$$

with $y^i \in \hat{\mathcal{C}}_1^{\kappa,\eta}$ $i = 1, 2$, and a regular part $y^r \in \hat{\mathcal{C}}_2^{\gamma+2\kappa,\eta}$, $y^{1,r} \in \hat{\mathcal{C}}_2^{\gamma+\kappa,\eta}$ with $\kappa < \kappa_0 \wedge \eta$. Furthermore, we asssume that the regularity of $y^1, y^2$ and $y^r, y^{1,r}$ can be related to those of $X$ by the following relation: $\gamma + 3\kappa > 1$, a condition that can be always fullfilled by a suitable choice of $\kappa$ whenever $\gamma + 3\kappa_0 > 1$. Denote this space of controlled paths by $\mathcal{Q}_{\kappa,\eta,\psi}$, or when this does not lead to an ambiguous situation, simply by $\mathcal{Q}_{\kappa,\eta}$ or $\mathcal{Q}$. Moreover, one can define a seminorm $\mathcal{N}$ on $\mathcal{Q}_{\kappa,\eta}$ in the following way:

$$\begin{aligned}\mathcal{N}[y; \mathcal{Q}_\kappa] = \mathcal{N}[y; \hat{\mathcal{C}}_1^{\kappa,\eta}] + \sum_{i=1,2} \mathcal{N}[y^i; \hat{\mathcal{C}}_1^{\infty,\eta}] + \sum_{i=1,2} \mathcal{N}[y^i; \hat{\mathcal{C}}_1^{\kappa,\eta}] \\ + \mathcal{N}[y^r; \hat{\mathcal{C}}_2^{\gamma+2\kappa,\eta}] + \mathcal{N}[y^{1,r}; \hat{\mathcal{C}}_2^{\gamma+\kappa,\eta}],\end{aligned}$$

where we recall that the notation $\mathcal{N}$ has been introduced at Section 3.3.

REMARK 5.2. Even if a weakly controlled path is, strictly speaking, given by a tuple $(y, y^1, y^2, y^r, y^{1,r})$ we will, with a slight abuse of notation, denote it with its first component, that is, simply $y$.

REMARK 5.3. The notion of weakly controlled path appeared first in [6] in the finite dimensional context as a way to linearize the space of rough paths around the driving control. Even if this linearization does not preserve the whole structure of the space of rough paths, it is enough to find solutions of rough differential equations.

With this notation at hand, we will try to implement now the strategy designed at the beginning of Section 4.1 in order to integrate a weakly controlled process $y$: let us first assume $x$ is a smooth process, and $y \in \mathcal{Q}$. Then $\mathcal{J}(\hat{d}x\, y)$ is well defined, and thanks to equations (48) and (87), we have

$$\begin{aligned}\mathcal{J}(\hat{d}x\, y) &= \mathcal{J}(\hat{d}x\, S) y + \mathcal{J}(\hat{d}x\, \hat{\delta} y) \\ &= \mathcal{J}(\hat{d}x\, S) y + \mathcal{J}(\hat{d}x\, X^1 y^1) + \mathcal{J}(\hat{d}x\, X^2 y^2) + \mathcal{J}(\hat{d}x\, y^r).\end{aligned}$$

Furthermore, for $s < t$, the term $\mathcal{J}_{ts}(\hat{d}x\, X^1 y^1)$ above only involves $y_s^1$, and hence the increment $\mathcal{J}_{ts}(\hat{d}x\, X^1 y^1)$ is equal to $\mathcal{J}_{ts}(\hat{d}x\, X^1) y_s^1$, that is, $X_{ts}^2 y_s^1$. This yields

$$\mathcal{J}(\hat{d}x\, y) = X^1 y + X^2 y^1 + X^3 y^2 + \mathcal{J}(\hat{d}x\, y^r). \tag{88}$$

Note that, in this last expression, the terms $X^1y$, $X^2y^1$ and $X^3y^2$ are well defined under Hypothesis 3. In order to have a well-defined expression for $\mathcal{J}(\hat{d}x\,y)$ in the rough case, it remains to handle the term $\mathcal{J}(\hat{d}x\,y^r)$. Then let us write

$$\mathcal{J}(\hat{d}x\,y^r) = \mathcal{J}(\hat{d}x\,y) - X^1y - X^2y^1 - X^3y^2,$$

and let us analyze this relation by applying $\hat{\delta}$ to both sides. This gives

$$\hat{\delta}[\mathcal{J}(\hat{d}x\,y^r)] = -\hat{\delta}[X^1y] - \hat{\delta}[X^2y^1] - \hat{\delta}[X^3y^2], \tag{89}$$

and notice that in the last expression, $\hat{\delta}[\mathcal{J}(\hat{d}x\,y^r)] \neq 0$, since $y^r$ belongs to $\hat{\mathcal{C}}_2$ instead of $\hat{\mathcal{C}}_1$. Moreover, a slight extension of Lemma 3.2 shows that, for $M \in \hat{\mathcal{C}}_2(\mathcal{L}(V))$ and $L \in \hat{\mathcal{C}}_1(V)$, we have

$$\hat{\delta}(ML) = \hat{\delta}M\,L - M\,\delta L = \tilde{\delta}M\,L - M\,\hat{\delta}L.$$

Applying this elementary relation to (89), we end up with

$$\begin{aligned}\hat{\delta}[\mathcal{J}(\hat{d}x\,y^r)] &= -\tilde{\delta}X^1\,y + X^1\,\hat{\delta}y - \tilde{\delta}X^2\,y^1 + X^2\,\hat{\delta}y^1 - \tilde{\delta}X^3\,y^2 + X^3\,\hat{\delta}y^2 \\ &= X^1(\hat{\delta}y - X^1y^1 - X^2y^2) + X^2(\hat{\delta}y^1 - X^1y^2) + X^3\,\hat{\delta}y^2 \\ &= X^1y^r + X^2y^{1,r} + X^3\,\hat{\delta}y^2\end{aligned} \tag{90}$$

under our hypothesis on $y$ and $X$ we have the following regularities:

$$X^1y^r \in \hat{\mathcal{C}}_3^{\gamma+3\kappa,-\rho}, \qquad X^2y^{1,r} \in \hat{\mathcal{C}}_3^{\gamma+\kappa_0+2\kappa,-\rho}, \qquad X^3\hat{\delta}y^2 \in \hat{\mathcal{C}}_3^{\gamma+2\kappa_0+\kappa,-\rho},$$

so if $\gamma + 3\kappa > 1$ we can apply the operator $\hat{\Lambda}$ and express $\mathcal{J}(\hat{d}x\,y)$ in terms of $\hat{\delta}$ and $\hat{\Lambda}$ only. Plugging (90) into (88), we get

$$\mathcal{J}(\hat{d}x\,y) = X^1y + X^2y^1 + X^3y^2 + \hat{\Lambda}(X^1y^r + X^2y^{1,r} + X^3\,\hat{\delta}y^2). \tag{91}$$

Similar to what we did in the Young case, we are now able to invert the procedure which lead to relation (91), by just invoking the assumptions made on $X^i$ and $y$:

THEOREM 5.4. *Let $x$ be a path such that $X^i, i = 1, 2, 3$, are well defined, and such that Hypothesis 3 holds true. Let also $y \in \mathcal{Q}_{\kappa,\eta,\psi}$ for $0 < \kappa < \kappa_0 < \gamma - \kappa$ and $\kappa \leq \eta$. Define $z \in \hat{\mathcal{C}}_1(\mathcal{B}_\eta)$ such that $z_0 = \psi$ and $\hat{\delta}z$ satisfies*

$$\hat{\delta}z \equiv \mathcal{J}(\hat{d}x\,y) = X^1y + X^2y^1 + X^3y^2 + \hat{\Lambda}(X^1y^r + X^2y^{1,r} + X^3\,\hat{\delta}y^2)$$

*and let $z^1 = y$, $z^2 = y^1$, $z^{1,r} = X^2y^2 + y^r$ so that $\hat{\delta}z^1 = X^1z^2 + z^{1,r}$. Then:*

(1) *$z$ is well defined as an element of $\mathcal{Q}_{\kappa,\eta}$, and coincides with the usual Riemann convolution of $y$ by $x$ in case $x$ and $y$ are smooth processes.*

(2) *The seminorm of $z$ in $\mathcal{Q}_{\kappa,\eta}$ can be estimated as*

$$\mathcal{N}[z;\mathcal{Q}_{\kappa,\eta}] \le c_X T^{\kappa_0-\kappa}(\|\psi\|_{\mathcal{B}_\kappa} + \mathcal{N}[y;\mathcal{Q}_{\kappa,\eta}]), \tag{92}$$

*for a positive constant $c_X$ depending only on $X^i$, $i=1,2,3$.*

(3) *It holds that, for any $0 \le s < t \le T$*

$$\mathcal{J}_{ts}(\hat{d}x\, y) = \lim_{|\Pi_{ts}|\to 0} \sum_{i=0}^{n} S_{tt_{i+1}}[X^1_{t_{i+1},t_i} y_{t_i} + X^2_{t_{i+1},t_i} y^1_{t_i} + X^3_{t_{i+1},t_i} y^2_{t_i}], \tag{93}$$

*where the limit is over all partitions $\Pi_{ts} = \{t_0 = t, \ldots, t_n = s\}$ of $[s,t]$ as the mesh of the partition goes to zero.*

Proof. We will divide again this proof in several steps.

*Step 1*: Let us start by evaluating the regularity of the terms in the right-hand side of (91), that is,

$$A = X^1 y, \qquad B = X^2 y^1, \qquad C = X^3 y^2,$$
$$D = \hat{\Lambda}(X^1 y^r + X^2 y^{1,r} + X^3\, \hat{\delta} y^2),$$

under our standing assumptions.

In order to bound $A$, we will first estimate $|y_s|_{\mathcal{B}}$ itself for $s \le T$: if $y \in \mathcal{Q}_\eta$, we have $y_s = S_s\psi + \hat{\delta} y_{s0}$, and hence

$$\|y\|_{0,\mathcal{B}_\eta} \le |\psi|_{\mathcal{B}_\eta} + T^\kappa \mathcal{N}[y;\hat{\mathcal{C}}_1^{\kappa,\eta}]. \tag{94}$$

In particular, $y$ is bounded in $\mathcal{B}_\eta$ on $[0,T]$. Thus, if $X^1 \in \hat{\mathcal{C}}_2^\gamma \mathcal{L}^{\eta,-\rho} \cap \hat{\mathcal{C}}_2^{\kappa_0}\mathcal{L}^{\eta,\eta}$, we have $X^1 y \in \hat{\mathcal{C}}_2^{\gamma,-\rho}$, and also $X^1 y \in \hat{\mathcal{C}}_2^{\kappa_0,\eta}$. Moreover,

$$|X^1_{ts} y_s|_{\mathcal{B}_\eta} \le \|X^1\|_{\kappa_0,\eta,\eta}(t-s)^{\kappa_0}(|\psi|_{\mathcal{B}_\eta} + T^\kappa \mathcal{N}[y;\hat{\mathcal{C}}_1^{\kappa,\eta}]),$$

and thus

$$\mathcal{N}[X^1 y;\hat{\mathcal{C}}_2^{\kappa,\eta}] \le \|X^1\|_{\kappa_0,\eta,\eta}(|\psi|_{\mathcal{B}_\eta} + T^\kappa \mathcal{N}[y;\hat{\mathcal{C}}_1^{\kappa,\eta}])T^{\kappa-\kappa_0}. \tag{95}$$

Let us estimate now the term $B$, that is $\mathcal{N}[X^2 y^1;\hat{\mathcal{C}}_2^{\gamma+\kappa,\eta}]$: since $y^1 \in \hat{\mathcal{C}}_1^{\infty,\eta}$ and $X^2 \in \hat{\mathcal{C}}_2^{2\kappa_0}\mathcal{L}^{\eta,\eta}$, we obtain again that $X^2 y^1 \in \hat{\mathcal{C}}_2^{2\kappa_0+\kappa,\eta}$, and we have

$$\mathcal{N}[X^2 y^1;\hat{\mathcal{C}}_2^{2\kappa,\eta}] \le \|X^2\|_{2\kappa_0,\eta,\eta}\mathcal{N}[y^1;\hat{\mathcal{C}}_1^{\infty,\eta}]T^{2(\kappa_0-\kappa)}. \tag{96}$$

The term $C$ can now be bounded along the same lines as for $A$ and $B$. Moreover, for the term $D$, as we already observed above, $X^1 y^r \in \hat{\mathcal{C}}_3^{\gamma+3\kappa,-\rho}$, $X^2 y^{1,r} \in \hat{\mathcal{C}}_3^{\gamma+\kappa_0+2\kappa,-\rho}$ and $X^3 \hat{\delta} y^2 \in \hat{\mathcal{C}}_3^{\gamma+2\kappa_0+\kappa,-\rho}$, and observe that we have assumed that $\gamma + 3\kappa > 1$. Thus, the operator $\hat{\Lambda}$ can be applied to $X^1 y^r +$

$X^2 y^{1,r} + X^3 \hat\delta y^2$, and invoking inequality (40), we get that

$$
\begin{aligned}
&\|\hat\Lambda(X^1 y^r + X^2 y^{1,r} + X^3\,\hat\delta y^2)\|_{\gamma+3\kappa,-\rho}\\
&\quad\le c\|X^1 y^r + X^2 y^{1,r} + X^3\,\hat\delta y^2\|_{\gamma+3\kappa,-\rho}\\
&\quad\le c(\|X^1\|_{\gamma,\eta,-\rho}\mathcal{N}[y^r;\hat{\mathcal{C}}_2^{3\kappa,\eta}] + \|X^2\|_{\gamma+\kappa_0,\eta,-\rho}\mathcal{N}[y^{1,r};\hat{\mathcal{C}}_2^{2\kappa,\eta}]\\
&\qquad\qquad + \|X^3\|_{\gamma+2\kappa_0,\eta,-\rho}\mathcal{N}[y^2;\hat{\mathcal{C}}_2^{\kappa,\eta}]).
\end{aligned}
\tag{97}
$$

Summarizing inequalities (94)–(97), we have obtained that $z$ is a well-defined element of $\hat{\mathcal{C}}_1^{\kappa,\eta}$, and that it satisfies

$$\|\hat\delta z\|_{\kappa,\eta} \le c_X T^{\kappa_0-\kappa}(|\psi|_{\mathcal{B}_\eta} + \mathcal{N}[y;\mathcal{Q}_{\kappa,\eta}]).$$

*Step 2*: Let us estimate now $z$ as an element of $\mathcal{Q}_{\kappa,\eta}$. The natural decomposition of $\hat\delta z$ is obviously $\hat\delta z = X z^1 + X^2 z^2 + z^r$, with

$$z^1 = y, \qquad z^2 = y^1 \quad\text{and}\quad z^r = X^3 y^2 + \hat\Lambda(X^1 y^r + X^2 y^{1,r} + X^3\,\hat\delta y^2).$$

It is now easily checked, along the same lines as for Step 1, that $z$ satisfies relation (92).

*Step 3*: In order to see how to get the convergence of the Riemann sums to $\mathcal{J}(\hat{d}x\,y)$ it is enough to remark that $\hat\delta z$ can be written as $\hat\delta z = (\mathrm{Id} - \hat\Lambda\hat\delta)[X^1 y + X^2 y^1 + X^3 y^2]$. Applying Corollary 3.6, we now get relation (93). $\square$

REMARK 5.5. The space of weakly controlled paths is a vector space with respect to the action of $\mathbb{R}$ but not with respect to other interesting linear endomorphisms of $\mathcal{B}$. The problem lies in the fact that for general linear $L:\mathcal{B}\to\mathcal{B}$ we can have $\hat\delta L\,y \ne L\,\hat\delta y$ since $L$ does not necessarily commute with the semigroup (which appears in the definition of $\hat\delta = \delta - a$).

5.3. *Linear evolution problem.* Let us turn now to the main aim of this section, which is to get an existence and uniqueness result for equation (84).

THEOREM 5.6. *Assume that Hypothesis 3 holds for the triple of incremental operators $X^1, X^2, X^3$ with $\gamma, \kappa_0, \kappa, \eta, \rho$ such that $\gamma = \kappa_0 + \eta + \rho$, $\gamma + 3\kappa_0 > 1$ and $\kappa < \kappa_0$. Then:*

(1) *Equation (84) admits a unique solution $y \in \mathcal{Q}_\eta$.*
(2) *The map $(\psi, X^1, X^2, X^3) \mapsto y$ is continuous.*
(3) *For $(t,s) \in \mathcal{S}_2$, the map $\Phi_{ts}:\mathcal{B}_\eta \to \mathcal{B}_\eta$, such that $\Phi_{ts}\psi = y_t$ when $y_s = \psi$ and $\hat\delta y_{ts} = \mathcal{J}_{ts}(\hat{d}x\,y)$ is a bounded linear endomorphism of $\mathcal{B}_\eta$, and it satisfies the cocycle property $\Phi_{tu}\Phi_{us} = \Phi_{ts}$.*

PROOF. Like in the Young case, the solution $y$ will be identified as the fixed point of the map $\Gamma : \mathcal{Q}_{\kappa,\eta} \to \mathcal{Q}_{\kappa,\eta}$ defined by $z = \Gamma(y)$, with $z_0 = \psi$ and $\hat{\delta} z_{ts} = \mathcal{J}_{ts}(\hat{d}x\, y)$. And here again, we will concentrate on the fact that, on a small interval $[0,T]$, the ball

$$B = \{y; y_0 = \psi, \mathcal{N}[y; \mathcal{Q}_{\kappa,\eta}] \le |\psi|_{\mathcal{B}_\eta}\}$$

is left invariant by the map $\Gamma$.

Indeed, whenever $y \in B$, then Theorem 5.4 asserts that for $z = \Gamma(y)$, the following estimate holds true:

$$\mathcal{N}[z; \mathcal{Q}_{\kappa,\eta}] \le c_X T^{\kappa_0 - \kappa}(|\psi|_{\mathcal{B}_\eta} + \mathcal{N}[y; \mathcal{Q}_{\kappa,\eta}]).$$

Hence, if one chooses a small enough $T$, so that $c_X T^{\kappa_0-\kappa} < 1/2$, it is readily checked that $\mathcal{N}[z; \mathcal{Q}_{\kappa,\eta}] \le |\psi|_{\mathcal{B}_\eta}$, which proves that $z \in B$. The contraction property is now a matter of standard arguments, and the remainder of the theorem follows easily. □

5.4. *Application: stochastic heat equation.* In the sequel of the paper, for sake of simplicity, the generic situation of a process $X$ with $\gamma$-Hölder continuity in time with $\gamma \le 1/2$ will be the case of an infinite-dimensional Brownian motion, given by the covariance function (45). For this special process, we will try to construct a pathwise solution to the linear stochastic heat equation on $[0,1]$. At the end of the section, we will give some hints about the way the fractional Brownian case should be treated.

5.4.1. *The Brownian case.* Like in the Young case, the key step in order to apply Theorem 5.6 to the Brownian setting is to define $(X^i)_{i=1,2,3}$ in a reasonable way, and then to check Hypothesis 3. We have chosen here to deal with an Itô type definition for $X^n$, and we get the following result:

PROPOSITION 5.7. *Let $X$ be an infinite-dimensional Brownian motion defined by the covariance structure (45), with $Q$ given by (44) for $\nu \in [0,1]$. For $n = 1,2,3$, let $X^n$ be the incremental operators given by (70) and (86), respectively, where the stochastic integral has to be understood in the Itô sense (see, e.g., [3, 25] for a complete definition). Then, almost surely,*

$$X^n \in \hat{\mathcal{C}}_2^{\gamma+(n-1)\kappa_0} \mathcal{L}_{\mathrm{HS}}^{\eta,-\rho} \cap \hat{\mathcal{C}}_2^{n\kappa_0} \mathcal{L}_{\mathrm{HS}}^{\eta,\eta}$$

*for any $\eta > 1/4$, $\gamma > \kappa_0 > \kappa$ satisfying*

$$\kappa_0 < 1/4 - \eta + \bar{\nu}/2 \quad \textit{and} \quad \gamma < 1/2,$$

*with $\bar{\nu} = \inf(\nu; 1/2)$.*

PROOF. We have already proved the regularity of $X^1$ in the fractional Brownian case. The proof in the current case would be similar, and we omit it. It will be enough to take $\eta = 1/4 + \varepsilon$, $\kappa_0 = \bar{\nu}/2 - 2\varepsilon$ and $\gamma = \kappa_0 + \rho + \varepsilon$, $\rho = 1/4 - \bar{\nu}/2$ for a given small $\varepsilon > 0$.

Let us concentrate then on the regularity properties of $X^2$: we will prove in fact first that $X^2 \in \hat{\mathcal{C}}_2^{\gamma+\kappa_0}\mathcal{L}_{\mathrm{HS}}^{\eta,-\rho}$, and for this, we will proceed along the same lines as for the proof of Proposition 4.4. Let us sketch the main steps which have to be followed.

*Step 1*: First of all, we have to estimate $\|A_o^{-\rho} X_{ts}^2 A^{-\eta}\|_{\mathrm{HS};\mathcal{B}\to\mathcal{B}}$, and it is readily checked that $A_o^{-\rho} X_{ts}^2 A^{-\eta}$ is represented by the kernel

$$\tilde{K}_{ts}(\xi,\eta) = \int_s^t \int_0^1 G_{t-u}^{-\rho}(\xi,\eta_1) X(du,d\eta_1) \times \int_s^u \int_0^1 G_{u-v}(\eta_1,\eta_2) X(dv,d\eta_2) G_{v-s}^{-\eta}(\eta_2,\eta).$$

Thus, when considered as an operator from $\mathcal{B}_\eta$ to $\mathcal{B}_{-\rho}$, we obtain that

$$E[\|X_{ts}^2\|_{\mathrm{HS}}^2] = \int_{[0,1]^2} E[(\tilde{K}_{ts}(\xi,\eta))^2]\, d\xi\, d\eta. \tag{98}$$

Moreover, some standard considerations about iterated integrals for Brownian noises (see, e.g., [3, 25]) yield

$$E[(\tilde{K}_{ts}(\xi,\eta))^2] = \int_s^t du \int_{[0,1]^2} G_{t-u}^{-\rho}(\xi,\eta_1) G_{t-u}^{-\rho}(\xi,\hat{\eta}_1) \times Q(\eta_1 - \hat{\eta}_1) H_{us}(\eta,\eta_1,\hat{\eta}_1)\, d\eta_1\, d\hat{\eta}_1,$$

with

$$H_{us}(\eta,\eta_1,\hat{\eta}_1) = \int_s^u dv \int_{[0,1]^2} G_{u-v}(\eta_1,\eta_2) G_{v-s}^{-\eta}(\eta_2,\eta) Q(\eta_2 - \hat{\eta}_2) \times G_{u-v}(\hat{\eta}_1,\hat{\eta}_2) G_{v-s}^{-\eta}(\hat{\eta}_2,\eta)\, d\eta_2\, d\hat{\eta}_2.$$

Plugging this equality into (98), we end up with

$$E[\|X_{ts}^2\|_{\mathrm{HS}}^2] = \int_0^\varepsilon du \int_{[0,1]^2} G_{2(\varepsilon-u)}^{-2\rho}(\eta_1,\hat{\eta}_1) \Psi_u(\eta_1,\hat{\eta}_1) \times Q(\eta_1 - \hat{\eta}_1)\, d\eta_1\, d\hat{\eta}_1, \tag{99}$$

where we have set $\varepsilon = t - s$ and

$$\Psi_u(\eta_1,\hat{\eta}_1) = \int_0^u dv \int_{[0,1]^2} G_{u-v}(\eta_1,\eta_2) G_{u-v}(\hat{\eta}_1,\hat{\eta}_2) \times G_{2v}^{-2\eta}(\eta_2,\hat{\eta}_2) Q(\eta_2 - \hat{\eta}_2)\, d\eta_2\, d\hat{\eta}_2.$$

Furthermore, using the spectral decomposition of $G_t$ and $Q$, introduced respectively, by (42) and (44), we obtain

$$\int_{[0,1]^2} G_{u-v}(\eta_1,\eta_2)G_{u-v}(\hat\eta_1,\hat\eta_2)G_{2v}^{-2\eta}(\eta_2,\hat\eta_2)Q(\eta_2-\hat\eta_2)\,d\eta_2\,d\hat\eta_2$$
$$= \sum_{i,j,k,l\in\mathbb{Z}} e_i(\eta_1)e_j(\hat\eta_1)\frac{e^{-(\lambda_i+\lambda_j)[u-v]}e^{-2\lambda_k v}}{\lambda_l^\nu\lambda_k^{2\eta}}\mathbf{1}_{\{i-k-l=0\}}\mathbf{1}_{\{j+k+l=0\}}.$$

Injecting again this value into (99) and using the fact that $\lambda_{-i}=\lambda_i$, we have that

$$E[\|X_{ts}^2\|_{\mathrm{HS}}^2] = \sum_{i,j,k,l,m,n\in E}\frac{1}{\lambda_l^\nu\lambda_m^{2\rho}\lambda_n^\nu\lambda_l^\nu\lambda_k^{2\eta}} \times\int_0^\varepsilon du\, e^{-2\lambda_m(\varepsilon-u)}\int_0^u dv\, e^{-2\lambda_j[u-v]}e^{-2\lambda_k v},$$

with

$$E=\{j,k,l,m,n\in\mathbb{Z}; m+n=j, k+l=-j\}. \tag{100}$$

Thus, we get

$$\begin{aligned}E[\|X_{ts}^2\|_{\mathrm{HS}}^2] &= \sum_{j,k,l,m,n\in E}\frac{1}{\lambda_m^{2\rho}\lambda_n^\nu\lambda_l^\nu\lambda_k^{2\eta}}\int_0^\varepsilon du\int_0^u dv\, e^{-2\lambda_m(\varepsilon-u)-2\lambda_j[u-v]-2\lambda_k v}\\ &\le c\sum_{j,k,l,m,n\in E}\frac{1}{\lambda_j^b\lambda_m^{2\rho+a}\lambda_n^\nu\lambda_k^{2\eta}\lambda_l^\nu}\int_0^\varepsilon du\int_0^u dv\,du\frac{dv}{(\varepsilon-u)^a(u-v)^b}\\ &\le c\varepsilon^{2-a-b}\sum_{j,k,l,m,n\in E}\frac{1}{\lambda_j^b\lambda_m^{2\rho+a}\lambda_n^\nu\lambda_k^{2\eta}\lambda_l^\nu}\int_0^1 du\int_0^u dv\,du\frac{dv}{(u-v)^b}.\end{aligned}$$

The double integral above is finite whenever $a,b\in(0,1)$, while the sum can be handled along the following lines: first, rewrite

$$S=\sum_{m,n,j:m+n-j=0}\frac{1}{\lambda_j^b\lambda_m^{2\rho+a}\lambda_n^\nu}\sum_{k,l:k+l=-j}\frac{1}{\lambda_k^{2\eta}\lambda_l^\nu},$$

and observe that, thanks to Lemma 4.6 and according to our hypothesis $\eta>1/4$, we have

$$\sum_{k,l:k+l=-j}\frac{1}{\lambda_k^{2\eta}\lambda_l^\nu}\le\sum_{k,l:k+l=-j}\frac{1}{\lambda_k^{2\eta}\lambda_l^{\bar\nu}}\le C\lambda_j^{-\bar\nu},$$

where $C$ stands again for a positive constant which can change from line to line. Then

$$S\le C\sum_{m,n,j:m+n-j=0}\frac{1}{\lambda_j^{b+\bar\nu}\lambda_m^{2\rho+a}\lambda_n^\nu}\le C\sum_m\frac{1}{\lambda_m^{2\rho+a}}\sum_{n,j:n-j=-m}\frac{1}{\lambda_j^{b+\bar\nu}\lambda_n^{\bar\nu}}$$

and choose $b > 1/2 - 2\bar{\nu}$, so that another application of Lemma 4.6 gives

$$S \le C \sum_m 1/\lambda_m^{2\rho+a+b+2\bar{\nu}-1/2}.$$

This latter sum is finite when $a + b > 1 - 2\rho - 2\bar{\nu}$. Then for any $\theta < \theta^*$ such that $2\theta^* < 1 + 2\rho + 2\bar{\nu}$, we have found that $E[\|X^2_{ts}\|^2_{\mathrm{HS}}] \le c\varepsilon^{2-a-b} \le c\varepsilon^{2\theta}$.

*Step 2*: One can go from $L^2$ to $L^p$ estimates for $m$ just like in Proposition 4.4 Step 4: indeed, we have

$$E[\|X^2_{ts}\|^{2p}_{\mathrm{HS}}] \le c_p \int_{[0,1]^{2p}} \prod_{i=1}^p E^{1/p}[\tilde{K}^{2p}_{ts}(\xi_i,\eta_i)]\, d\xi_1\, d\eta_1 \cdots d\xi_p\, d\eta_p.$$

Moreover, $\tilde{K}_{ts}(\xi_i,\eta_i)$ is a variable of the second chaos $\mathcal{H}_2$ with respect to the Gaussian field $X$, and invoking [19, Relation (1.61)], the $L^2$ and $L^{2p}$ norms on $\mathcal{H}_2$ are equivalent. Thus, for any integer $p \ge 1$, there exists a constant $c_p$ such that

$$E[\|X^2_{ts}\|^{2p}_{\mathrm{HS}}] \le c_p(t-s)^{2p\theta}.$$

*Step 3*: We will conclude now thanks to Lemma 3.8, which reads here as:

$$\begin{aligned}\|X^2\|_{\gamma_2,\eta,-\rho} &\le c[U_{\gamma_2+2/p,p,\eta,-\rho}(X^2) + \|\tilde{\delta}X^2\|_{\gamma_2,\eta,-\rho}] \\ &= c[U_{\gamma_2+2/p,p,\eta,-\rho}(X^2) + \|X^1X^1\|_{\gamma_2,\eta,-\rho}],\end{aligned}$$

for any integer $p \ge 1$. According to the previous step, it is then easily checked that, for any $\gamma_2 < \theta^* = 1/2 + \bar{\nu} + \rho$, and $p$ large enough, the term $U_{\gamma_2+2/p,p,\eta,-\rho}(X^2)$ can be bounded almost surely by a finite constant. Recall now that we have chosen $\gamma = \kappa_0 + \rho + \eta$, $\kappa_0 = \bar{\nu}/2 - 2\varepsilon$ and $\eta = 1/4 + \varepsilon$. Thus,

$$\gamma + \kappa_0 = 2\kappa_0 + \rho + \eta = \bar{\nu} - \eta + 1/4 + \rho < \theta^*,$$

and hence $U_{\gamma_3+2/p,p,\eta,-\rho}(X^2) < \infty$ for any $\gamma_3 \le \gamma + \kappa_0$.

Let us treat now the term $X^1X^1$. Along the same lines as in Proposition 4.4, it can be shown that $X^1 \in \hat{\mathcal{C}}_2^{\gamma}\mathcal{L}^{\eta,-\rho}$ and $X^1 \in \hat{\mathcal{C}}_2^{\kappa_0}\mathcal{L}^{\eta,\eta}$. Hence, by composition of operators, we get $X^1X^1 \in \hat{\mathcal{C}}_3^{\gamma+\kappa_0}\mathcal{L}^{\eta,-\rho}$, which means that $\|X^1X^1\|_{\gamma_2,\eta,-\rho}$ is finite for any $\gamma_3 \le \gamma + \kappa_0$. Summing up this short discussion, we have obtained that

$$\|X^2\|_{\gamma_2,\eta,-\rho} \text{ finite a.s.} \qquad \text{for any } \gamma_3 \le \gamma + \kappa_0.$$

One can proceed then to prove that $X^2 \in \hat{\mathcal{C}}_2^{2\kappa_0}\mathcal{L}^{\eta,\eta}$ by a slight elaboration of the computations above. This easy exercise is left to the reader.

The proof for the operator $X^3$ follows the same lines and will not be reported. Indeed, we prefer to concentrate on the regularity properties of higher order operators in the more complex situation of Section 6.4. $\square$

We are now able to apply our abstract results to the stochastic heat equation.

THEOREM 5.8. *Let $X$ be an infinite-dimensional Brownian motion on $[0,T]\times[0,1]$, defined by the covariance function given by (45) and (44) with $\nu>1/3$. Then there exists $\eta>1/4$, $0<\kappa<\gamma<1/2$ such that $\kappa<\kappa_0$ and $\gamma+3\kappa>1$ such that, for any $\psi\in\mathcal{B}_\eta$ the equation*

$$Y(0,\xi)=\psi(\xi), \qquad \partial_t Y(t,\xi)=\Delta Y(t,\xi)\,dt+Y(t,\xi)X(dt,d\xi), \qquad t\in[0,T],\xi\in[0,1],$$

*with periodic boundary conditions, understood as equation (84), has a unique solution in $\mathcal{Q}_{\kappa,\eta,\psi}$.*

PROOF. Like in the proof of Theorem 4.12, the claim is readily checked once we have shown that $X^n$, $n=1,2,3$ satisfy Hypothesis 3. This amount to check that there exist $\kappa_0>\kappa$ and $\gamma<1/2$ such that

$$\kappa_0<1/4-\eta+\bar{\nu}/2, \qquad \gamma+3\kappa>1.$$

However, thanks to Proposition 5.7, it is enough to take $\eta=1/4+\varepsilon$, $\kappa_0=\bar{\nu}/2-2\varepsilon$, $\kappa=\kappa_0-\varepsilon$, $\rho=1/4-\bar{\nu}/2$ and $\gamma=\kappa_0+\rho+\eta$ for some small $\varepsilon>0$. The condition $\gamma+3\kappa_0>1$ can then be read $\gamma+3\kappa_0=1/2+3\bar{\nu}/2-4\varepsilon>1$, which is possible whenever $\nu>1/3$. □

5.4.2. *The fractional Brownian case.* In order to define an integration theory for the fractional Brownian motion beyond the Young case one has to start, like in the Brownian case, by defining the operators $X^1$ and $X^2$ in a natural way. We have already seen that $X^1$ could be understood by means of Wiener integrals, and for $X^2$, two reasonable choices for the definition of (86) seem to be the use of either Skorokhod or Stratonovich integrals with respect to the fractional Brownian motion $X$. However, it turns out that these two solutions are equally unsatisfactory, for two different reasons that we proceed to detail now:

(1) When one computes moments of random variables of the second chaos defined by Stratonovich integrals, some trace terms appear, a classical phenomenon which is explained for instance in [19] in the general case, in [23] for the stochastic heat equation, or in [21] for the fBm. In the current situation, if we want these trace terms to be convergent for a fractional Brownian motion $X$ defined by (43) and (44), one has to choose $\nu>1/2$, which means in particular that $Q$ is a bounded function of $\xi\in[0,1]$. In other words, we are not allowed, even if $H>1/2$, to consider a distribution-valued noise in space, which was one of our main aim.

(2) The Skorokhod integral works better as far as convergences and regularity estimates are concerned. But one of the basic ingredients of our algebraic manipulations on integrals is the fact that one can write, under suitable hypothesis:

$$\int_s^t S_{tu}\, dx_u\, b_s = \left[\int_s^t S_{tu}\, dx_u\right] b_s,$$

an equality which is known to fail in the Skorokhod case (see [19] again for further explanations). For instance, the relation $\tilde{\delta} X^2 = X^1 X^1$, which is useful in our analysis, does not hold true when $\hat{\delta} X^2$ is understood as a Skorokhod integral.

In order to cope with these problems, one can adopt the following strategy: compute the correction term $P$, understood as a 2-increment operator-valued process, which allows to write

$$\tilde{\delta} X^2 = X^1 X^1 + P, \tag{101}$$

when $X^2$ is defined via Skorokhod integration. Notice that, since $X^2$ is an element of $\mathcal{H}_2$, the process $P$ is deterministic.

Recall now that the operator $\tilde{\delta}$ has been defined as follows: for a suitable Banach space $V$ and $M \in \hat{\mathcal{C}}_k^{\infty}(V)$, set

$$\tilde{\delta} M = \hat{\delta} M - Ma = \delta M - aM - Ma. \tag{102}$$

Then the operator $\tilde{\delta}$ enjoys the same kind of properties as $\hat{\delta}$, and in particular, $\tilde{\delta}\tilde{\delta} = 0$ and $\ker \tilde{\delta}|_{\hat{\mathcal{C}}_3} = \operatorname{Im} \tilde{\delta}|_{\hat{\mathcal{C}}_2}$. Moreover, relation (101) can be read as $\tilde{\delta} P = 0$. Thus, there exists another process $T \in \hat{\mathcal{C}}_2 \mathcal{L}$ such that $\tilde{\delta} T = P$. Consider then an explicit version of $T$ and set $\tilde{X}^2 = X^2 - T$. Then $\tilde{X}^2$ is still a Lévy area type process, such that $\tilde{\delta} \tilde{X}^2 = X^1 X^1$, which means that hopefully, $\tilde{X}^2$ will enjoy both algebraic and analytic properties allowing a nice extension of the notion of convolution integral. However, an open problem is to understand in which sense the integrals defined using this *corrected* Lévy area $\tilde{X}^2$ are useful and/or natural. We plan to report on this possibility in a further paper.

5.5. *The algebra of a rough path.* Bilinear stochastic equations, in finite or infinite dimensions, are often handled by means of chaos decomposition (see, e.g., [14, 15]). In this section, we will try to stress some relationships between our pathwise approach and this latter method.

Our basic Hypothesis 3 states that

$$\tilde{\delta} X^n = \sum_{k=1}^{n-1} X^k X^{n-k} \tag{103}$$

for $n = 1, 2, 3$, and moreover that

$$X^n \in \hat{\mathcal{C}}_2^{\gamma+(n-1)\kappa_0}\mathcal{L}^{\eta,-\rho} \quad \text{and} \quad X^n \in \hat{\mathcal{C}}_2^{n\kappa_0}\mathcal{L}^{\eta,\eta}, \tag{104}$$

with $\gamma + 3\kappa_0 > 1$, $\eta > 1/4$ and $\rho = 1/4 - \bar{\nu}/2$. Furthermore, it can be shown, along the same lines as for Theorem 3.5, that there exists an inverse $\tilde{\Lambda}$ to $\tilde{\delta}$ on $\hat{\mathcal{C}}_3^{\mu}\mathcal{L}^{0,-\rho} \cap \ker(\tilde{\delta})$ for a certain $\mu > 1$.

Let us see now how to construct an operator $X^4$ satisfying the operator Chen relation (103): by composition of operators, and invoking Hypothesis (104), it is easily checked that $X^1X^3 + X^3X^1 + X^2X^2 \in \hat{\mathcal{C}}_3^{\gamma+3\kappa_0}\mathcal{L}^{\eta,-\rho}$. Furthermore, we have assumed that $\gamma + 3\kappa_0 > 1$, and thus, by analogy with Theorem 5.4, we will set now $X^4 := \tilde{\Lambda}[X^1X^3 + X^2X^2 + X^3X^1]$, which is well defined as an element of $\mathcal{E}_2^{\gamma+3\kappa_0}\mathcal{L}^{\eta,-\rho}$ and thus that belongs to $\hat{\mathcal{C}}_2^{4\kappa_0}\mathcal{L}^{\eta,\eta}$ (since $\gamma = \kappa_0 + \eta + \rho$). It turns out that this procedure can be iterated, and we obtain the following proposition.

PROPOSITION 5.9. *Let $X$ satisfying Hypothesis 3. Then one can construct a sequence $\{X^n; n \geq 4\}$ out of $X^1, X^2, X^3$, such that, for any $\kappa < \kappa_0 = \gamma_n - \rho$, we have $X^n \in \mathcal{E}_2^{\gamma+(n-1)\kappa_0}\mathcal{L}^{\eta,-\rho}$,*

$$\|X^n\|_{\hat{\mathcal{C}}_2^{n\kappa_0}\mathcal{L}^{\eta,-\rho}} \leq C(n!)^{-\kappa_0}$$

*and such that the operator Chen relations (103) are satisfied.*

PROOF. The $X^n$ are constructed by an induction on $n$. Then it is clear that $X^n \in \mathcal{E}_2^{\gamma+(n-1)\kappa_0}\mathcal{L}^{\eta,-\rho}$. Moreover, for $n \geq 4$ we have $n\kappa_0 > 1$, so that $X^n = \tilde{\Lambda}(\sum_{k=1}^{n-1} X^{n-k}X^k)$ can be defined directly as an element of $\mathcal{L}^{\eta,\eta}$. Then the same kind of arguments as in the finite dimensional case [10] prove that we have the inductive bound

$$\|X^n\|_{\hat{\mathcal{C}}_2^{n\kappa_0}\mathcal{L}^{\eta,\eta}} \leq C_X(n!)^{-\kappa_0}. \tag{105}$$

□

In such a setting, the lifted rough path allows to express the Itô map which sends initial conditions to solution of the linear equation (84) by a convergent series of operators.

COROLLARY 5.10. *Under the conditions of Proposition 5.9, there exists an operator $T$, defined as an element of $\mathcal{E}_2^{\gamma+(n-1)\kappa_0}\mathcal{L}^{\eta,-\rho}$, given by the strongly convergent series $T := \sum_{k=1}^{\infty} X^k$, and such that the solution of the linear problem (84) satisfy the equation $\hat{\delta}y = Ty$, or written in another way*

$$y_t = S_{ts}y_s + T_{ts}y_s, \qquad (t,s) \in \mathcal{S}^2.$$

*In particular, if we define $X_{ts}^0 = S_{ts}$ and set $\hat{T} = X^0 + T$ we have that $\hat{T}$ is a cocycle of operators:*

$$\hat{T}_{ts} = \hat{T}_{tu}\hat{T}_{us}, \qquad (t,u,s) \in \mathcal{S}^3.$$

PROOF. The convergence of the series for $T$ in the operator norm follows from the bound (105) on $\|X^n\|$. The cocycle property is proven as in finite dimension. The uniqueness of the solution to the linear problem allows to identify the operator $T$ as the Itô map for the rough evolution equation. □

**6. Polynomial nonlinearities.** Going back to the general setting explained at Sections 3.2 and 3.4, we will consider now an equation of the form

$$y_t = S_t\psi + \int_0^t S_{tu}\,dx_u\,M_n(y_u^{\otimes n}), \tag{106}$$

where $y$ lives in the Hilbert space $\mathcal{B}$, where $M_n : \mathcal{B}^n \to \mathcal{B}$ is some unbounded multilinear operator from the Hilbert tensor $\mathcal{B}^n = \mathcal{B}^{\otimes n}$ to $\mathcal{B}$, and where we understand $\varphi^{\otimes n} = \varphi \otimes \cdots \otimes \varphi \in \mathcal{B}^n$ as the tensor monomial generated by $\varphi \in \mathcal{B}$. In fact, for sake of simplicity we assume that $M_n$ is symmetric and we restrict our discussion to the case $n = 2$, letting $M_2 = B$, the general situation posing no more conceptual difficulties. Then our general strategy, like in the linear case, will be first to expand equation (106) for a smooth driving process $x$ in order to guess the appropriate rough-path underlying this equation. It will be seen that those expansions involve some increments indexed by trees. Studying the algebraic and analytic properties of these increments, we shall obtain a reasonable notion of solution to our quadratic equation.

6.1. *Formal expansions and trees.* Let us first simplify a little our setting. Recall that we wish to solve an equation of the form

$$y_t = S_t\psi + \int_0^t S_{tu}\,dx_u\,B(y_u^{\otimes 2}), \tag{107}$$

where we specialize our situation in the following way: assume first that $\mathcal{B} = L^2([0,1])$, which means that we are back again to the heat equation setting of Section 3.4. Then $B : \mathcal{B}^{\otimes 2} \to \mathcal{B}$ is defined by $[B(\phi \otimes \psi)](\xi) = \phi(\xi)\psi(\xi)$ for $\xi \in [0,1]$, whenever this expression makes sense in $\mathcal{B}$. Assume for the moment that $x \in \hat{\mathcal{C}}_1^1 \mathcal{L}^{\kappa,\kappa}$ for $\kappa$ large enough. We can expand equation (107) as

$$\begin{aligned}
y_t &= S_{ts}y_s + \int_s^t S_{tu}\,dx_u B((S_{us}y_s)^{\otimes 2}) \\
&\quad + 2\int_s^t S_{tu}\,dx_u\,B\left(S_{us}y_s \otimes \int_s^u S_{uv}\,dx_v\,B((S_{vs}y_s)^{\otimes 2})\right) \\
&\quad + \int_s^t S_{tu}\,dx_u\,B\left(\int_s^u S_{uv}\,dx_v\,B((S_{vs}y_s)^{\otimes 2}) \right. \\
&\qquad\qquad \left. \otimes \int_s^u S_{uv}\,dx_v\,B((S_{vs}y_s)^{\otimes 2})\right)
\end{aligned}$$

(108)

$$
\begin{aligned}
&+ 4 \int_s^t S_{tu}\, dx_u \\
&\qquad \times B\biggl(S_{us} y_s \otimes \int_s^u S_{uv}\, dx_v \\
&\qquad\qquad \times B\biggl((S_{vs} y_s) \otimes \int_s^v S_{vw}\, dx_w\, B((S_{ws} y_s)^{\otimes 2})\biggr)\biggr) \\
&+ \text{h.o. iterated integrals.}
\end{aligned}
$$

As we see, iterated integrals appear here in combinations which are not as easy to handle as in the bilinear case of Section 5.3. A natural way to code this kind of expansion is to use planar trees, as explained below.

Without entering too much into formal definitions involving trees, let us mention that we shall consider planar binary rooted trees $\mathcal{T}$ of the form

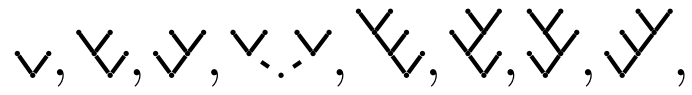 , etc...

which allow to give a compact expression of the iterated integrals appearing in the expansion (108). Observe that each tree can be constructed from the trivial tree $\tau_0 = \bullet$ by using the binary operation $V : \mathcal{T} \times \mathcal{T} \to \mathcal{T}$ consisting in gluing two trees at a newly created root, so for example

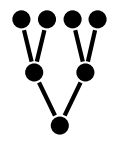 $= V(V(\tau_0, \tau_0), V(\tau_0, \tau_0)).$

For any tree $\tau \in \mathcal{T}$, we associate the function $d(\tau)$ that counts the number of leaves on the trees, so that $d(\tau_0) = 1$ and $d(V(\tau_1, \tau_2)) = d(\tau_1) + d(\tau_2)$.

Let us see now how to represent expansion (108) thanks to planar trees. Define recursively an operator-valued increment $X^\tau \in \hat{\mathcal{C}}_2 \mathcal{L}(\mathcal{B}_\eta^{d(\tau)}; \mathcal{B})$ for $\tau \in \mathcal{T}$ as

$$
X_{ts}^{\tau_0} = S_{ts} \quad \text{and} \quad X_{ts}^{V(\tau_1,\tau_2)} = \int_s^t S_{tu}\, dx_u\, B(X_{us}^{\tau_1} \otimes X_{us}^{\tau_2}). \tag{109}
$$

Notice that $X^\tau$ has always to be considered as an operator acting on $\mathcal{B}_\eta^{d(\tau)}$. For instance, we understand that, if $\tau = V(\tau_1, \tau_2)$, we have

$$
\begin{aligned}
&X_{ts}^{V(\tau_1,\tau_2)}(\varphi_1 \otimes \cdots \otimes \varphi_{d(\tau)}) \\
&\quad = \int_s^t S_{tu}\, dx_u\, B(X_{us}^{\tau_1}(\varphi_1 \otimes \cdots \otimes \varphi_{d(\tau_1)}) \\
&\qquad\qquad \otimes X_{us}^{\tau_2}(\varphi_{d(\tau_1)+1} \otimes \cdots \otimes \varphi_{d(\tau)})).
\end{aligned}
$$

This latter formula justifies also in a sense the use of planar trees, since in general $X^{V(\tau_1,\tau_2)} \neq X^{V(\tau_2,\tau_1)}$. In order to illustrate this fact, consider the simple example where $\tau_1 = \tau_0$ and $\tau_2 = V(\tau_0, \tau_0)$. Then $d(V(\tau_1, \tau_2)) = 3$ and

$$X_{ts}^{V(\tau_1,\tau_2)}(\varphi_1 \otimes \varphi_2 \otimes \varphi_3) = \int_s^t S_{tu}\, dx_u\, B(X_{us}^{\tau_0}(\varphi_1) \otimes X_{us}^{V(\tau_0,\tau_0)}(\varphi_2 \otimes \varphi_3)),$$

while

$$X_{ts}^{V(\tau_2,\tau_1)}(\varphi_1 \otimes \varphi_2 \otimes \varphi_3) = \int_s^t S_{tu}\, dx_u\, B(X_{us}^{V(\tau_0,\tau_0)}(\varphi_1 \otimes \varphi_2) \otimes X_{us}^{\tau_0}(\varphi_3)),$$

which are a priori clearly different objects.

With this notation in hand, it is now checked that our previous expansion (108) can be written in a simpler way as

$$(\hat{\delta} y)_{ts} = X_{ts}^{\vee}(y_s^{\otimes 2}) + 2X_{ts}^{\vee^{\backslash}}(y_s^{\otimes 3}) + X_{ts}^{\vee\cdot\cdot\vee}(y_s^{\otimes 4}) + 4X_{ts}^{\vee^{\backslash\backslash}}(y_s^{\otimes 4}) + r, \tag{110}$$

where $r \in \hat{\mathcal{C}}_2(\mathcal{B})$ is some remainder term, and where we took care to distinguish the various operators obtained by permuting the factors in the $\mathcal{B}$-tensors. Of course, we could have expanded the solution further, and some operators associated to larger trees would have appeared. However, in a smooth enough situation, the strategy in order to solve (107) is now clear: we can use the map $\hat{\Lambda}$ to eliminate the remainder from the equation:

$$\hat{\delta} y = (1 - \hat{\Lambda}\hat{\delta})[X^{\vee}(y^{\otimes 2}) + 2X^{\vee^{\backslash}}(y^{\otimes 3}) + X^{\vee\cdot\cdot\vee}(y^{\otimes 4}) + 4X^{\vee^{\backslash\backslash}}(y^{\otimes 4})] \tag{111}$$

and try to solve this by fixed-point method. The only condition we need to check is that

$$\hat{\delta}[X^{\vee}(y^{\otimes 2}) + 2X^{\vee^{\backslash}}(y^{\otimes 3}) + X^{\vee\cdot\cdot\vee}(y^{\otimes 4}) + 4X^{\vee^{\backslash\backslash}}(y^{\otimes 4})] \tag{112}$$

should be in the domain of $\hat{\Lambda}$, which means that its time-regularity should be greater than 1. The computation of expressions like (112) requires a little algebraic preparation.

6.2. *Algebraic computations.* To ease some computations, we introduce an "improper" increment $E_{ts} = \mathrm{Id}$ (improper because it does not vanish as $t = s$), where the Id has to be understood, according to the context, as the identity operator on the vector space under consideration. For example, we can write $\hat{\delta} h = \delta h - (S - E)h$. Moreover, we also introduce $e_{ts} = 1$ taking values in $\mathbb{R}$, so that for example, if $z \in \hat{\mathcal{C}}_2(\mathcal{B})$, then $(ze)_{tus} = z_{tu} e_{us} = z_{tu}$.

It will also be useful to extend the action of $\hat{\delta}$ to the tensors $\mathcal{B}^n$ by letting $\hat{\delta} z = \delta z - (\mathcal{S} - E)z$, where $\mathcal{S} : \mathcal{B}^n \to \mathcal{B}^n$ is defined as $\mathcal{S} = S \otimes \cdots \otimes S$ for any $n \geq 1$. If the reader is unconfortable with giving the same name at different operators, he can think that $\mathcal{S}$ is defined on the direct sum $\bigoplus_{n \geq 1} \mathcal{B}^n$; furthermore, we will write explicitly $\mathcal{S}_n$ when the context is insufficient to determine the actual space on which $\mathcal{S}$ is defined. Analogously to the case of $\hat{\delta}$, the operator $\tilde{\delta}$ defined by (102) can be allowed to act on $\hat{\mathcal{C}}_1\mathcal{L}(\mathcal{B}^n,\mathcal{B})$ as $\tilde{\delta}H = \delta H - (\mathcal{S}-E)H - H(\mathcal{S}-E)$.

We wish first to understand how the operators $\hat{\delta}$ and $\tilde{\delta}$ act on tensor products. More specifically, we shall need three relations which are summarized in the following lemmas.

Lemma 6.1. *The following relations hold true:*

(1) *Let $z,w\in\mathcal{C}_1(\mathcal{B})$. Then*

$$\hat{\delta}(z\otimes w) = Sz\otimes\hat{\delta}w + \hat{\delta}z\otimes Sw + \hat{\delta}z\otimes\hat{\delta}w. \tag{113}$$

(2) *Let $z,w\in\mathcal{C}_2(\mathcal{B})$. Then*

$$\begin{aligned}\hat{\delta}(z\otimes w) &= ze\otimes\hat{\delta}w + ze\otimes\mathcal{S}w + \hat{\delta}z\otimes we + \mathcal{S}z\otimes we\\ &\quad + \mathcal{S}z\otimes\hat{\delta}w + \hat{\delta}z\otimes\mathcal{S}w + \hat{\delta}z\otimes\hat{\delta}w.\end{aligned} \tag{114}$$

(3) *Let $Z\in\hat{\mathcal{C}}_2\mathcal{L}(\mathcal{B}^a,\mathcal{B}^b)$ and $W\in\hat{\mathcal{C}}_2\mathcal{L}(\mathcal{B}^c,\mathcal{B}^d)$. Then*

$$\begin{aligned}\tilde{\delta}(Z\otimes W) &= Z\mathcal{S}_a\otimes\tilde{\delta}W + Z\mathcal{S}_a\otimes\mathcal{S}_dW + \tilde{\delta}Z\otimes W\mathcal{S}_c\\ &\quad + \mathcal{S}_bZ\otimes W\mathcal{S}_c + \mathcal{S}_bZ\otimes\tilde{\delta}W + \tilde{\delta}Z\otimes\mathcal{S}_dW + \tilde{\delta}Z\otimes\tilde{\delta}W,\end{aligned} \tag{115}$$

*where an example of notational convention is given by*

$$(Z\mathcal{S}_a\otimes\mathcal{S}_dW)_{tus} = (Z_{tu}S_{us}^{\otimes a})\otimes(S_{tu}^{\otimes d}W_{us})\in\mathcal{L}(\mathcal{B}^{a+c},\mathcal{B}^{c+d}).$$

Proof. These relations are easily checked by elementary computations. We include the proof of the third one for sake of completeness: notice that

$$\tilde{\delta}Z = Z^{\diamond} - ZE - EZ - (\mathcal{S}_b-E)Z - Z(\mathcal{S}_a-E) = Z^{\diamond} - Z\mathcal{S}_a - \mathcal{S}_bZ,$$

where $Z^{\diamond}_{tus} = Z_{ts}$, and thus

$$\begin{aligned}\tilde{\delta}(Z\otimes W) &= Z^{\diamond}\otimes W^{\diamond} - \mathcal{S}_bZ\otimes\mathcal{S}_dW - Z\mathcal{S}_a\otimes W\mathcal{S}_c\\ &= (\mathcal{S}_bZ + Z\mathcal{S}_a + \tilde{\delta}Z)\otimes(\mathcal{S}_dW + W\mathcal{S}_c + \tilde{\delta}W)\\ &\quad - \mathcal{S}_bZ\otimes\mathcal{S}_dW - Z\mathcal{S}_a\otimes W\mathcal{S}_c,\end{aligned}$$

which yields relation (115) by a straightforward expansion. □

We also want to understand how $\hat{\delta},\tilde{\delta}$ act on the operators $X^\tau$. A first relation in this direction is to note that, according to Lemma 3.2, if $X^\tau\in\hat{\mathcal{C}}_2\mathcal{L}(\mathcal{B}^{\otimes n};\mathcal{B})$ and $h\in\hat{\mathcal{C}}_1(\mathcal{B}^{\otimes n})$,

$$\hat{\delta}[X^\tau h] = (\hat{\delta}X^\tau)h - X^\tau\delta h, \quad \text{that is,} \quad \hat{\delta}[X^\tau h] = (\tilde{\delta}X^\tau)h - X^\tau\hat{\delta}h. \tag{116}$$

It is thus useful to compute quantities of the form $\tilde{\delta} X^{\tau}$. To this purpose, consider $n \geq 1$ and define $I : \hat{\mathcal{C}}_2 \mathcal{L}(\mathcal{B}^n, \mathcal{B}^2) \to \hat{\mathcal{C}}_2 \mathcal{L}(\mathcal{B}^n, \mathcal{B})$ by

$$I(H)_{ts} = \int_s^t S_{tu}\, dx_u\, BH_{us}. \tag{117}$$

This kind of expression can be related to our tree-indexed increments by noticing, for instance, that

$$X^{\vee} = I(\mathcal{S}_2), \qquad X^{{}^{\vee}\!\vee} = I(X^{\vee} \otimes \mathcal{S}_1), \qquad X^{\vee^{\vee}} = I(\mathcal{S}_1 \otimes X^{\vee}),$$

and generally speaking, (109) can be read as

$$\begin{aligned} & X^{V(\tau_1,\tau_2)} = I(X^{\tau_1} \otimes X^{\tau_2}), \qquad X^{V(\tau_0,\tau)} = I(\mathcal{S} \otimes X^{\tau}) \quad \text{and} \\ & X^{V(\tau,\tau_0)} = I(X^{\tau} \otimes \mathcal{S}). \end{aligned} \tag{118}$$

Hence, we shall compute differentials of terms of the form $I(H)$.

Lemma 6.2. *Let $H \in \hat{\mathcal{C}}_2^1 \mathcal{L}(\mathcal{B}^n, \mathcal{B}^2)$. The following formulae hold true for the derivative of $I(H)$:*

$$(\tilde{\delta} I(H))_{tus} = I_{tu}(\mathcal{S}_2) H_{us} + I_{tu}(\tilde{\delta}(H))$$

*and*

$$(\hat{\delta} I(H))_{tus} = I_{tu}(\mathcal{S}_2) H_{us} + I_{tu}(\hat{\delta}(H)).$$

*Furthermore, if we assume that $\tilde{\delta} H$ can be decomposed as $(\tilde{\delta} H)_{tus} = \sum_{j \leq M} H_{tu}^{(1,j)} H_{us}^{(2,j)}$, for a given $M \geq 1$, $H^{(1,j)} \in \hat{\mathcal{C}}_2^1 \mathcal{L}(\mathcal{B}^2, \mathcal{B}^2)$, and $H^{(2,j)} \in \hat{\mathcal{C}}_2^1 \mathcal{L}(\mathcal{B}^n, \mathcal{B}^2)$, then we obtain*

$$\tilde{\delta} I(H) = I(\mathcal{S}_2) H + \sum_{j \leq M} I(H^{(1,j)}) H^{(2,j)}.$$

Proof. We have

$$\begin{aligned} [\tilde{\delta} I(H)]_{tus} &= I_{ts}(H) - I_{tu}(H) \mathcal{S}_{us} - \mathcal{S}_{tu} I_{us}(H) \\ &= \int_s^t S_{tw}\, dx_w\, BH_{ws} - \int_u^t S_{tw}\, dx_w\, BH_{wu} S_{us} \\ &\quad - S_{tu} \int_s^u S_{uw}\, dx_w\, BH_{ws} \\ &= \int_u^t S_{tw}\, dx_w\, BH_{ws} - \int_u^t S_{tw}\, dx_w\, BH_{wu} S_{us} \\ &= \int_u^t S_{tw}\, dx_w\, B(\mathcal{S}_{wu} H_{us}) \end{aligned}$$

$$
\begin{aligned}
&+ \int_u^t S_{tw}\, dx_w \, [BH_{ws} - B(\mathcal{S}_{wu} H_{us}) - BH_{wu} S_{us}] \\
&= I_{tu}(\mathcal{S}_2) H_{us} + I_{tu}(\tilde{\delta}(H)),
\end{aligned}
$$

which proves our first assertion. The second one is now trivially deduced. □

With these preliminaries in hand, we can now compute the action of $\tilde{\delta}$ on the tree-indexed increments we have met so far, in the following way.

Lemma 6.3. *Let $x$ be a smooth $\mathcal{L}(\mathcal{B})$-valued path. Then we have*

$$
\begin{aligned}
\tilde{\delta} X^{\vee} &= 0, \\
\tilde{\delta} X^{\vee_l} &= I(\mathcal{S}_2)(X^{\vee} \otimes S) = X^{\vee}(X^{\vee} \otimes S), \\
\tilde{\delta} X^{\vee\cdot\cdot\vee} &= X^{\vee}(X^{\vee} \otimes X^{\vee}) + X^{\vee_r}(X^{\vee} \otimes S) + X^{\vee_l}(S \otimes X^{\vee})
\end{aligned}
$$

*and*

$$
\begin{aligned}
\tilde{\delta} X^{\vee_{ll}} &= X^{\vee}(X^{\vee_l} \otimes S) + X^{\vee_l}(X^{\vee} \otimes \mathcal{S}_2), \\
\tilde{\delta} X^{\vee_{lr}} &= X^{\vee}(X^{\vee_r} \otimes S) + X^{\vee_l}(\mathcal{S}_2 \otimes X^{\vee}), \\
\tilde{\delta} X^{\vee_{rl}} &= X^{\vee}(S \otimes X^{\vee_l}) + X^{\vee_r}(X^{\vee} \otimes \mathcal{S}_2), \\
\tilde{\delta} X^{\vee_{rr}} &= X^{\vee}(S \otimes X^{\vee_r}) + X^{\vee_r}(\mathcal{S}_2 \otimes X^{\vee}).
\end{aligned}
$$

Proof. All these relations are obtained by elementary computations, and we shall only sketch the proof for some of them: first of all, invoking Lemma 6.2, we get:

$$
\tilde{\delta} X^{\vee} = \tilde{\delta} I(\mathcal{S}_2) = I(\mathcal{S}_2)\mathcal{S}_2 - I(\mathcal{S}_2)\mathcal{S}_2 = 0,
$$

where we used the fact that $\tilde{\delta}\mathcal{S} = -\mathcal{S}\mathcal{S}$. As far as $X^{\vee_l}$ is concerned, we have

$$
\tilde{\delta} X^{\vee_l} = \tilde{\delta} I(X^{\vee} \otimes S) = I(\mathcal{S}_2)(X^{\vee} \otimes S),
$$

since $\tilde{\delta}(X^{\vee} \otimes S) = 0$ by a direct computation using formula (115). Similarly, it holds that

$$
\begin{aligned}
\tilde{\delta} X^{\vee\cdot\cdot\vee} &= \tilde{\delta} I(X^{\vee} \otimes X^{\vee}) \\
&= I(\mathcal{S}_2)(X^{\vee} \otimes X^{\vee}) + I(S \otimes X^{\vee})(X^{\vee} \otimes S) + I(X^{\vee} \otimes S)(S \otimes X^{\vee}),
\end{aligned}
$$

owing to the fact that

$$
\tilde{\delta}(X^{\vee} \otimes X^{\vee}) = \mathcal{S} X^{\vee} \otimes X^{\vee} \mathcal{S} + X^{\vee} \mathcal{S} \otimes \mathcal{S} X^{\vee}.
$$

Now, invoking (118), we have $I(S \otimes X^{\vee})(X^{\vee} \otimes S) = X^{\check{\vee}}$ and $I(X^{\vee} \otimes S) = X^{\hat{\vee}}$, which yields

$$\tilde{\delta} X^{\vee\!\cdot\!\cdot\vee} = X^{\vee}(X^{\vee} \otimes X^{\vee}) + X^{\check{\vee}}(X^{\vee} \otimes S) + X^{\hat{\vee}}(S \otimes X^{\vee}),$$

as claimed in our lemma. □

It is important to note now that all the previous computations have been performed for a smooth path $x$. However, we shall ask our driving process $x$ to satisfy the following assumption:

HYPOTHESIS 4. *We assume that the path $x$ allows to define some increments $X^{\tau}$ for any $\tau \in \mathcal{T}$ such that $d(\tau) \leq 4$. We also suppose that those increments satisfy the relations of Lemma 6.3, and that the following Hölder regularities hold true: setting $|\tau| = d(\tau) - 1$, we have $X^{\tau} \in \hat{\mathcal{C}}^{|\tau|\kappa_0} \mathcal{L}(\mathcal{B}_{\eta}^{d(\tau)}, \mathcal{B}_{\eta})$ and $X^{\tau} \in \hat{\mathcal{C}}^{\gamma + (|\tau|-1)\kappa_0} \mathcal{L}(\mathcal{B}_{\eta}^{d(\tau)}, \mathcal{B}_{-\rho})$, with $\gamma + n\kappa_0 > 1$ and $\gamma = \kappa_0 + \eta + \rho$, for a given $\eta > 1/4$.*

REMARK 6.4. Here again, it is important to work in spaces of the form $\mathcal{B}_{\eta}$ with $\eta > 1/4$. Indeed, these spaces are algebras, which ensures at least that, whenever $\phi, \psi \in \mathcal{B}_{\eta}$, then $B(\phi, \psi) \in \mathcal{B}$.

REMARK 6.5. The peculiar relation between the various parameters involved in Hypothesis 4 has been suggested by the example which we treat later on and it is due to the mixing between space and time regularity due to the analytic semigroup. As operators the $X^{\tau}$'s can map to more regular spaces (with respect to the scale associated the generator $A$) at the price of loosing some time regularity. This is a phenomenon which is not peculiar of infinite-dimensional rough paths associated to random processes but it is found also in the rough-path approach to deterministic PDEs like the Korteweg–de-Vries equation or the Navier–Stokes equation [7–9].

6.3. *A space of integrable paths.* The general discussion of the bilinear equation requires a deep understanding of the algebra of $X$. It is not the aim of this paper to enter into these kind of considerations, and we prefer here to concentrate on a particular case where $\kappa$ is sufficiently large to stop the expansions at some low (but nontrivial) order. So here we assume that $\gamma + 3\kappa > 1$.

In order to solve the fixed-point problem associated to (107), we introduce a new space of weakly controlled paths, denoted by $\mathcal{Q}_{X,\kappa}$, which enjoys some nice stability properties under the map $\Gamma : y \mapsto z = S\psi + I(y \otimes y)$.

Definition 6.6. Let $\psi \in \mathcal{B}_\eta$ an initial condition, and $x$ a driving noise satisfying Hypothesis 4, with $\gamma + 3\kappa > 1$. We say that a path $y \in \hat{\mathcal{C}}_1^{*,\kappa}(\mathcal{B}_\eta)$ belongs to $\mathcal{Q}_{X,\kappa}$ if $y_0 = \psi$, and $\hat{\delta} y$ can be decomposed into

$$\hat{\delta} y = X^{\vee} y^{\vee} + X^{{}^{\vee}\!\vee} y^{{}^{\vee}\!\vee} + X^{\vee^{\vee}} y^{\vee^{\vee}} + y^{\sharp}, \tag{119}$$

where $y^{\vee}, y^{{}^{\vee}\!\vee}$ and $y^{\vee^{\vee}}$ can be written as:

$$y^{\vee} = w \otimes w, \qquad y^{{}^{\vee}\!\vee} = w^{\vee} \otimes w, \qquad y^{\vee^{\vee}} = w \otimes w^{\vee},$$

and the following regularities hold true:

$$y \in \hat{\mathcal{C}}_1^{*,\kappa}(\mathcal{B}_\eta), \qquad w, y \in \hat{\mathcal{C}}_1^{*,\kappa_1}(\mathcal{B}_\eta), \qquad w^{\vee} \in \hat{\mathcal{C}}_1^{*,\kappa_2}(\mathcal{B}_\eta^2), \qquad y^{\sharp} \in \hat{\mathcal{C}}_2^{3\kappa_2}(\mathcal{B}_\eta),$$

where $\kappa > \kappa_1 > \kappa_2$, $\kappa - \kappa_1 = \kappa_1 - \kappa_2 \equiv \mu$ and $\gamma + 3\kappa_2 > 1$. On $\mathcal{Q}_{X,\kappa}$, we define the seminorm

$$\begin{aligned} \mathcal{N}[y; \mathcal{Q}_{X,\kappa}] &= \mathcal{N}[y; \hat{\mathcal{C}}_1^{\kappa}(\mathcal{B}_\eta)] + \mathcal{N}[w; \hat{\mathcal{C}}_1^{\kappa_1}(\mathcal{B}_\eta)] \\ &\quad + \mathcal{N}[w^{\vee}; \hat{\mathcal{C}}_1^{\kappa_2}(\mathcal{B}_\eta^2)] + \mathcal{N}[y^{\sharp}; \hat{\mathcal{C}}_2^{3\kappa_2}(\mathcal{B}_\eta)]. \end{aligned}$$

Note that the constant path $y_t = S_t \psi$ is a controlled path whenever $\psi \in \mathcal{B}_\kappa$, and in this case $w, w^{\vee}, y^{\sharp}$ are all identically zero. Furthermore, the space $\mathcal{Q}$ satisfies the following useful stability property.

Theorem 6.7. *Assume that $x$ satisfies Hypothesis 4, where we recall that $\gamma + 3\kappa > 1$ and $\kappa_0 > \kappa$. For $y \in \mathcal{Q}_{\kappa,X}$, define $z \equiv \Gamma(y) \in \mathcal{Q}_{\kappa,X}$ by $z_0 = \psi$ and a decomposition of the form*

$$\hat{\delta} z = X^{\vee} z^{\vee} + X^{{}^{\vee}\!\vee} z^{{}^{\vee}\!\vee} + X^{\vee^{\vee}} z^{\vee^{\vee}} + z^{\sharp},$$

*with $z^{\vee} = w_z \otimes w_z$, $z^{{}^{\vee}\!\vee} = w_z^{\vee} \otimes w_z$, $z^{\vee^{\vee}} = w_z \otimes w_z^{\vee}$,*

*where $w_z = y$, $w_z^{\vee} = y^{\vee}$, and $z^{\sharp} \in \hat{\mathcal{C}}_2^{3\kappa_0}(\mathcal{B})$ is a remainder which can be written as*

$$\begin{aligned} z^{\sharp} &= X^{\vee\cdot\cdot\vee}(y^{\vee} \otimes y^{\vee}) + X^{{}^{\vee}\!\vee^{\vee}}(y^{{}^{\vee}\!\vee} \otimes y) + X^{{}^{\vee^{\vee}}\!\vee}(y^{\vee^{\vee}} \otimes y) \\ &\quad + X^{\vee^{{}^{\vee}\!\vee}}(y \otimes y^{{}^{\vee}\!\vee}) + X^{\vee^{\vee^{\vee}}}(y \otimes y^{\vee^{\vee}}) - \Lambda(J), \end{aligned}$$

*where $J$ is defined by relation (121). Then:*

(1) *$\Gamma : \mathcal{Q}_{\kappa,X} \to \mathcal{Q}_{\kappa_0,X}$ is well defined.*
(2) *$\hat{\delta} z$ coincides with $I(y \otimes y)$ in the smooth case.*
(3) *The following estimate holds true: for all $0 < S < T$ we have:*

$$\begin{aligned} &\mathcal{N}[z; \mathcal{Q}_{X,\kappa_0}([0,S])] \\ &\quad \le C_X \left(1 + |\psi|_\eta + |w_0|_\eta + |w_0^{\vee}|_\eta + S^{\mu} \mathcal{N}[y; \mathcal{Q}_{X,\kappa}([0,S])]\right)^4, \end{aligned} \tag{120}$$

*for a positive constant $C_X$ which only depends on the rough path $X$.*

PROOF. Start from two smooth paths $x$ and $y$. If we apply the $I$ map defined at (117), we obtain, just as in the Young case (48),

$$\begin{aligned}\hat{\delta} z &= I(y\otimes y) = I(\hat{\delta}(y\otimes y) + \mathcal{S}_2(y\otimes y)) = I(\mathcal{S}_2)(y\otimes y) + I(\hat{\delta}(y\otimes y)) \\ &= I(\mathcal{S}_2)(y\otimes y) + I(Sy\otimes \hat{\delta} y) + I(\hat{\delta} y\otimes Sy) + I(\hat{\delta} y\otimes \hat{\delta} y),\end{aligned}$$

where we have used Lemma 6.1. Expanding $\hat{\delta} y$ in this equation and invoking relation (118), we thus obtain

$$\begin{aligned}\hat{\delta} z = {}& X^{\vee}(y\otimes y) + X^{{}^{\vee}\!\vee}(y^{\vee}\otimes y) + X^{\vee^{\vee}}(y\otimes y^{\vee}) + X^{\vee\cdot\cdot\vee}(y^{\vee}\otimes y^{\vee}) \\ &+ X^{{}^{{}^{\vee}\!\vee}\!\vee}(y^{{}^{\vee}\!\vee}\otimes y) + X^{{}^{\vee^{\vee}}\!\vee}(y^{\vee^{\vee}}\otimes y) + X^{\vee^{{}^{\vee}\!\vee}}(y\otimes y^{{}^{\vee}\!\vee}) + X^{\vee^{\vee^{\vee}}}(y\otimes y^{\vee^{\vee}}) + z^{\flat},\end{aligned}$$

where $z^{\flat}$ has to be understood again as a remainder. Now, by our standard, argument we shall define $z^{\flat}$ in the nonsmooth case by $z^{\flat} = -\hat{\Lambda}(J)$, where $J$ is given by

$$\begin{aligned}J = {}& \hat{\delta}[X^{\vee}(y\otimes y) + X^{{}^{\vee}\!\vee}(y^{\vee}\otimes y) + X^{\vee^{\vee}}(y\otimes y^{\vee}) + X^{\vee\cdot\cdot\vee}(y^{\vee}\otimes y^{\vee}) \\ &+ X^{{}^{{}^{\vee}\!\vee}\!\vee}(y^{{}^{\vee}\!\vee}\otimes y) + X^{{}^{\vee^{\vee}}\!\vee}(y^{\vee^{\vee}}\otimes y) + X^{\vee^{{}^{\vee}\!\vee}}(y\otimes y^{{}^{\vee}\!\vee}) + X^{\vee^{\vee^{\vee}}}(y\otimes y^{\vee^{\vee}})].\end{aligned} \tag{121}$$

In order for this equation to be well defined, we still need to check that $J$ belongs to $\hat{\mathcal{C}}_3^{\gamma+3\kappa}(\mathcal{B}_{-\rho})$, which is in the domain of $\hat{\Lambda}$ since $\gamma + 3\kappa > 1$. Let us then compute $J$: owing to (116), we have

$$\begin{aligned}\hat{\delta}[X^{\vee}(y\otimes y)] &= [\tilde{\delta} X^{\vee}](y\otimes y) - X^{\vee}[\hat{\delta}(y\otimes y)] \\ &= -X^{\vee}(Sy\otimes \hat{\delta} y) - X^{\vee}(\hat{\delta} y\otimes Sy) - X^{\vee}(\hat{\delta} y\otimes \hat{\delta} y),\end{aligned}$$

thanks to relation (113) and Lemma 6.3. The other terms can be computed along the same lines, and here is a sample of what is obtained:

$$\begin{aligned}\hat{\delta}[X^{{}^{\vee}\!\vee}(y^{\vee}\otimes y)] &= X^{\vee}(X^{\vee} y^{\vee}\otimes Sy) - X^{{}^{\vee}\!\vee}(\hat{\delta} y^{\vee}\otimes Sy) \\ &\quad - X^{{}^{\vee}\!\vee}(Sy^{\vee}\otimes \hat{\delta} y) - X^{{}^{\vee}\!\vee}(\hat{\delta} y^{\vee}\otimes \hat{\delta} y), \\ \hat{\delta}[X^{{}^{{}^{\vee}\!\vee}\!\vee}(y^{{}^{\vee}\!\vee}\otimes y)] &= X^{\vee}(X^{{}^{\vee}\!\vee} y^{{}^{\vee}\!\vee}\otimes Sy) + X^{{}^{\vee}\!\vee}(X^{\vee}\otimes \mathcal{S}_2)(y^{{}^{\vee}\!\vee}\otimes y) \\ &\quad - X^{{}^{{}^{\vee}\!\vee}\!\vee}\hat{\delta}(y^{{}^{\vee}\!\vee}\otimes y)\end{aligned}$$

and

$$\begin{aligned}\hat{\delta}[X^{\vee\cdot\cdot\vee}(y^{\vee}\otimes y^{\vee})] &= X^{\vee}(X^{\vee} y^{\vee}\otimes X^{\vee} y^{\vee}) + X^{\vee^{\vee}}(X^{\vee} y^{\vee}\otimes \mathcal{S}_2 y^{\vee}) \\ &\quad + X^{{}^{\vee}\!\vee}(\mathcal{S}_2 y^{\vee}\otimes X^{\vee} y^{\vee}) - X^{\vee\cdot\cdot\vee}(\mathcal{S}_2 y^{\vee}\otimes \hat{\delta} y^{\vee}) \\ &\quad - X^{\vee\cdot\cdot\vee}(\hat{\delta} y^{\vee}\otimes \mathcal{S}_2 y^{\vee}) - X^{\vee\cdot\cdot\vee}(\hat{\delta} y^{\vee}\otimes \hat{\delta} y^{\vee}).\end{aligned}$$

Now, by gathering all the terms we have obtained, we obtain that $J = \sum_{k=1}^{4} J_k$, with

$$
\begin{aligned}
J_1 = X^{\vee}[&Sy \otimes (-\hat{\delta} y + X^{\vee} y^{\vee} + X^{\vee_{\ell}} y^{\vee_{\ell}} + X^{\vee_{r}} y^{\vee_{r}}) \\
&+ (-\hat{\delta} y + X^{\vee} y^{\vee} + X^{\vee_{\ell}} y^{\vee_{\ell}} + X^{\vee_{r}} y^{\vee_{r}}) \otimes Sy \\
&- \hat{\delta} y^{\vee} \otimes \hat{\delta} y^{\vee} + (X^{\vee} \otimes X^{\vee})(y^{\vee} \otimes y^{\vee})], \\
J_2 = X^{\vee_{\ell}}[&-Sy^{\vee} \otimes \hat{\delta} y - \hat{\delta} y^{\vee} \otimes Sy + (\mathcal{S}_2 \otimes X^{\vee})(y^{\vee} \otimes y^{\vee}) - \hat{\delta} y^{\vee} \otimes \hat{\delta} y \\
&+ (X^{\vee} \otimes \mathcal{S}_2)(y^{\vee_{\ell}} \otimes y) + (S \otimes X^{\vee} \otimes S)(y^{\vee_{r}} \otimes y)]
\end{aligned}
$$

and

$$
\begin{aligned}
J_3 = X^{\vee_{r}}[&-Sy \otimes \hat{\delta} y^{\vee} - \hat{\delta} y \otimes Sy^{\vee} + (X^{\vee} \otimes \mathcal{S}_2)(y^{\vee} \otimes y^{\vee}) - \hat{\delta} y \otimes \hat{\delta} y^{\vee} \\
&+ (\mathcal{S}_2 \otimes X^{\vee})(y \otimes y^{\vee_{r}}) + (S \otimes X^{\vee} \otimes S)(y \otimes y^{\vee_{\ell}})], \\
J_4 = -X^{\vee\cdot\cdot\vee}&\hat{\delta}(y^{\vee} \otimes y^{\vee}) - X^{\vee_{\ell\ell}} \hat{\delta}(y^{\vee_{\ell}} \otimes y) - X^{\vee_{\ell r}} \hat{\delta}(y^{\vee_{r}} \otimes y) \\
&- X^{\vee_{r\ell}} \hat{\delta}(y \otimes y^{\vee_{\ell}}) - X^{\vee_{rr}} \hat{\delta}(y \otimes y^{\vee_{r}}).
\end{aligned}
$$

Furthermore, notice that, using equation (119) for the increments of $y$, the quantity $J_1$ can be simplified into

$$
J_1 = X^{\vee}[-Sy \otimes y^{\sharp} - y^{\sharp} \otimes Sy - \hat{\delta} y^{\vee} \otimes \hat{\delta} y^{\vee} + (X^{\vee} \otimes X^{\vee})(y^{\vee} \otimes y^{\vee})].
$$

We are now left with the cumbersome task which consists in analyzing the regularity of all the terms we have produced so far. We shall just focus on one particular example, namely $X^{\vee_{\ell}}(\hat{\delta} y^{\vee} \otimes Sy)$, leaving the other ones to the patient reader. Invoking again Lemma 6.1, we have

$$
\begin{aligned}
X^{\vee_{\ell}}(\hat{\delta} y^{\vee} \otimes Sy) = X^{\vee_{\ell}}(\hat{\delta} w \otimes \hat{\delta} w \otimes Sy) + X^{\vee_{\ell}}(Sw \otimes \hat{\delta} w \otimes Sy) \\
+ X^{\vee_{\ell}}(\hat{\delta} w \otimes Sw \otimes Sy).
\end{aligned}
$$

Among the three terms in the right hand side of this relation, we shall analyze the first one, the other ones being similar: recall that $X^{\vee_{\ell}} \in \mathcal{C}_2^{\gamma+2\kappa_0} \mathcal{L}(\mathcal{B}_\eta^3, \mathcal{B}_{-\rho})$, $\hat{\delta} w \in \mathcal{C}_2^{\kappa}(\mathcal{B}_\eta)$ and $Sy \in \mathcal{C}_2^{0}(\mathcal{B}_\eta)$. Thus $X^{\vee_{\ell}}(\hat{\delta} w \otimes \hat{\delta} w \otimes Sy) \in \mathcal{C}_3^{\gamma+2\kappa_0+\kappa}(\mathcal{B}_{-\rho})$, which is enough regularity to apply the $\Lambda$-map. The other terms in the decomposition of $J$ can be treated similarly, which ends the proof of our first assertion.

Our second claim being immediate from the construction of our integral, let us say a few words about the last one. Here again, many terms have to be estimated, and we shall focus on a representative example, namely the term $w_z^{\vee} = y^{\vee} = w \otimes w$. In fact, the quantity $\mathcal{N}[w \otimes w; \hat{\mathcal{C}}_1^{\kappa_2}(\mathcal{B}_\eta^2)]$ has to be estimated, and recall that, according to Lemma 6.1, the following relation holds true:

$$
\hat{\delta}(w \otimes w) = Sw \otimes \hat{\delta} w + \hat{\delta} w \otimes Sw + \hat{\delta} w \otimes \hat{\delta} w.
$$

Thus, since $w \in \hat{\mathcal{C}}_1^{\kappa_1}(\mathcal{B}_\eta)$ and $S_t$ is a bounded operator on $\mathcal{B}_\eta$ for any positive $S$, we obtain

$$\begin{aligned}\mathcal{N}[w \otimes w; \hat{\mathcal{C}}_1^{\kappa_2}([0,S]; \mathcal{B}_\eta^2)] &\le c(1 + |w_0|_\eta + S^\mu \mathcal{N}[w; \hat{\mathcal{C}}_1^{\kappa_1}([0,S]; \mathcal{B}_\eta^2)])^2 \\ &\le c(1 + |w_0|_\eta + S^\mu \mathcal{N}[y; \mathcal{Q}_{X,\kappa}([0,S])])^2,\end{aligned}$$

where we used also the decomposition $w_t = S_{t0} w_0 + \hat{\delta} w_{t0}$ to bound $w_t$ in terms of $w_0$ and $\hat{\delta} w$:

$$|w_t|_\eta \le |w_0|_\eta + S_1^\kappa \mathcal{N}[w; \hat{\mathcal{C}}_1^{\kappa_1}([0,S]; \mathcal{B}_\eta^2)].$$

The other terms defining $\hat{\delta} z$ can be treated along the same lines, which proves relation (120). □

We can turn now to the main goal of this section, which is to get an existence and uniqueness result for equation (107).

THEOREM 6.8. *Assume that $x$ allows to define some incremental operators $X^\tau$ for any $\tau \in \mathcal{T}$ such that $d(\tau) \le 4$, and that these increments satisfy Hypothesis 4, for $\gamma, \kappa_0, \kappa, \eta$ such that $\kappa < \kappa_0 < \gamma$, $\gamma + 3\kappa > 1$ and $\eta > 1/4$. Then there exists a strictly positive $T_0 = T_0(X^\tau; d(\tau) \le 4)$ such that equation (107) admits a unique solution $y \in \mathcal{Q}_{\kappa,X}([0,T_0])$.*

PROOF. The proof of this result is very similar to those of Theorems 4.3 and 5.6, and we shall omit the details here. Just notice that inequality (120) allows to construct an invariant ball for the map $\Gamma$ in $\mathcal{Q}_{\kappa,X}([0,T_0])$, whenever $T_0$ is small enough. The contraction argument can then be written in a standard way. □

6.4. *The Brownian case.* In this section, we investigate the behavior of the operators $X^\tau$ defined above, when $x = X$ is an infinite-dimensional Brownian motion, defined at Section 3.4. Our aim is of course to show that, under certain conditions, $X$ satisfies Hypothesis 4. To this purpose, for the remainder of the section, we will mainly consider some applications on the space $\mathcal{B}_\eta$ for $\eta = 1/4 + \varepsilon$ and a small $\varepsilon > 0$. Let us also introduce an additional notation: for the remainder of the article, we will write $A \lesssim B$ for two real quantities $A$ and $B$ when $A \le cB$ for a universal constant $c$.

PROPOSITION 6.9. *Let $X$ be an infinite-dimensional Brownian motion defined by the covariance structure (45), with $Q$ given by (44) for $\nu \in (1/3, 1/2]$. Recall also that, for a planar binary tree, we have set $|\tau| = d(\tau) - 1$. Let $X^\tau$ be the incremental operator given by (109) where the stochastic integrals have to be understood in the Itô sense. Then, almost surely, $X^\tau \in$*

$\hat{\mathcal{C}}_2^{\kappa_0|\tau|}\mathcal{L}_{\mathrm{HS}}(\mathcal{B}_\eta^{d(\tau)};\mathcal{B}_\eta)$ *for* $\tau = \vee, \vee\!\!\vee, \vee\!\!\vee, \vee\cdot\vee, \vee\!\!\vee\!\!\vee, \vee\!\!\vee\!\!\vee, \vee\!\!\vee\!\!\vee, \vee\!\!\vee\!\!\vee$ *and for any* $\kappa_0$ *satisfying* $0 < \kappa_0 < 1/4 - \eta + \bar{\nu}/2$. *Moreover,* $X^\tau \in \hat{\mathcal{C}}_2^{\gamma+\kappa_0(|\tau|-1)-1/4}\mathcal{L}_{\mathrm{HS}}(\mathcal{B}_\eta^{d(\tau)};\mathcal{B}_{-\rho})$ *for* $\gamma = \kappa_0 + \eta + \rho < 1/2$. *Theorem 6.8 can then be applied in this situation.*

PROOF. In all the cases, the line of the proof is the same. We obtain an $L^2$ estimate on the Hilbert–Schmidt norm of $X^\tau$, which by Gaussian tools can be boosted to an $L^p$ bound for any $p$. Applying Lemma 3.8, the result is then easily deduced.

Admitting for the moment the results of Lemma 6.11 below, let us give some details about our method. Since our incremental processes always belong to a finite chaos of the infinite-dimensional Brownian motion $X$, it is easily deduced from Lemma 6.11 that

$$E[\|X_{ts}^\tau\|^p_{\mathrm{HS},\mathcal{L}(\mathcal{B}_\eta^{d(\tau)};\mathcal{B}_\eta)}] \lesssim (t-s)^{p|\tau|(\kappa_0+\varepsilon+(1-1/|\tau|)/2)} \lesssim (t-s)^{p|\tau|(\kappa_0+\varepsilon)}$$

for any $0 < \kappa_0 < 1/4 - \eta + \bar{\nu}/2 < \eta/2$. Moreover, it also holds that:

$$\|\tilde{\delta} X_{tus}^\tau\| \leq \sum_{\tau^1,\tau^2} \|X_{tu}^{\tau^1}\|\|X_{us}^{\tau^2}\| \lesssim (t-u)^{|\tau^1|\kappa_0}(u-s)^{|\tau^2|\kappa_0} \lesssim (t-s)^{|\tau|\kappa_0},$$

where $\tau^1, \tau^2$ denote the trees appearing in the expansion for $\tilde{\delta} X_{tus}^\tau$ given at Lemma 6.3, for which we have always $|\tau^1|, |\tau^2| \geq 1$ and $|\tau^1| + |\tau^2| = |\tau| + 1$. So it is clear that using the extended G–R–R Lemma 3.8, we obtain

$$\|X_{tus}^\tau\| \lesssim (t-s)^{|\tau|\kappa_0},$$

for any $\tau$ such that $d(\tau) \leq 4$. Finally, observe that the conditions $\gamma < 1/2$, $\kappa_0 = \bar{\nu}/2 - \varepsilon$ and $3\kappa_0 = \gamma > 1$ force us to choose $\nu > 1/3$. $\square$

An easy consequence of the last estimations is an existence result for a Brownian SPDE in the rough-path sense:

THEOREM 6.10. *Let* $X$ *be an infinite-dimensional Brownian motion on* $[0,T] \times [0,1]$, *defined by the covariance function given by (45) and (44) with* $\nu > 1/3$. *Then there exists* $\eta > 1/4$, $0 < \kappa < \gamma < 1/2$ *satisfying* $\kappa < \kappa_0$ *and* $\gamma + 3\kappa > 1$ *such that, for any* $\psi \in \mathcal{B}_\eta$ *the equation*

$$Y(0,\xi) = \psi(\xi), \qquad \partial_t Y(t,\xi) = \Delta Y(t,\xi)\,dt + Y(t,\xi)^2 X(dt,d\xi), \qquad t \in [0,T], \xi \in [0,1],$$

*with periodic boundary conditions, understood as equation (84), has a unique local solution in* $\mathcal{Q}_{\kappa,\eta,\psi}$ *up to a time* $T_*$ *which depends on the initial condition and on the operators* $X^\tau$, $|\tau| \leq 3$.

PROOF. Like in the proof of Theorem 5.8, the proof amounts to check the validity of Hypothesis 4 in the light of Proposition 6.9. □

The rest of the paper is dedicated to the $L^2$ estimations for the operators $X^\tau$. In fact, we will obtain a slightly stronger result than the one we claimed at Proposition 6.9.

LEMMA 6.11. *For the trees considered at Proposition 6.9, we have the following $L^2$ bounds:*

$$E[\|X_{ts}^\tau\|^2_{\mathrm{HS},\mathcal{L}(\mathcal{B}_\eta^{d(\tau)};\mathcal{B}_\eta)}] \lesssim (t-s)^{|\tau|\Delta-1/2}, \tag{122}$$

*where $\Delta = 1 - 2\eta + \nu - \varepsilon$ for some arbitrary small $\varepsilon > 0$.*

PROOF. It is conceptually easy to generalize the arguments of Proposition 5.7 to reduce the problem to an estimation of a mixed sum (over eigenvalues of $A_o$) and integral over time variables (after contraction of the stochastic integrals). This long and tedious task is left to the reader. We prefer to give a diagrammatical algorithm which allows to go from the kernel on $L^2$ (associated to each operator) to a simple sum estimation. This will be detailed in the next two subsections. □

REMARK 6.12. We can interpret this result by the following heuristic considerations. The situation more similar to the finite dimensional theory is when $\eta$ is slightly larger than $1/4$ and $\nu$ slightly larger than $1/2$. In this case, $\Delta$ is arbitrarily near to 1 which would give the classical scaling of Brownian increments if we could ignore the factor $-1/2$ apperaring at the exponent in the r.h.s. of equation (122). This further loss of regularity is due to the need of estimating the Hilbert–Shmidt norm. We conjecture that some technique which would allow to estimate directly the operator norms of the increments would give a better time regularity which would improve the overall theory. Apart from this technical difficulty there is an intrinsic departure from the Brownian regularity due to a loss in space regularity which must be compensated via a transfer from time to space regularity (allowed by the convolutions). This loss in space regularity has two sources: one is the nonlinear operation which start to be badly behaved when $\eta$ is smaller than $1/4$, the other is the presence of the noise which degrades the spatial regularity of the result when $\nu$ is smaller than $1/2$. Both contributions are clearly accounted in the formula for the effective time regularity $\Delta$.

REMARK 6.13. As we noted also elsewhere in the context of rough-paths associated to deterministic PDEs [7–9], in the infinite-dimensional setting objects like the semigroup $S$ (and in general the unbounded linear operators appearing in the equations) must be considered at the same level of the driving stochastic processes in the sense that the ensemble of these objects form a rough-path. In this perspective, the fact that the rough path $X^\tau$ which we construct depends both on the Gaussian noise and on the specific convolution semigroup $S$ should not be considered more unusual that the fact that in the finite-dimensional theory the higher order iterated integrals depends on the vector of all irregular components driving the equation.

6.5. *Diagrammatica.* We will first show, at a heuristic level and on a simple example, how to pass from an incremental operator to a graph for the computation of Hilbert–Schmidt norms.

(1) *Case of the operator* $X^{\vee}$. Consider first the operator $X^{\vee}:\mathcal{B}_\eta^2 \to \mathcal{B}_\eta$. Recalling the notation of Section 3.4, an orthonormal basis for $\mathcal{B}_\eta$ is given by $\{\tilde{e}_j; j \geq 1\}$, where $\tilde{e}_j = \lambda_j^{-\eta} e_j$. Furthermore, the particular form (44) we have assumed on the covariance function $Q$ implies that $x$ can be decomposed as

$$x_u = \sum_{p\in\mathbb{Z}} \lambda_p^{-\nu/2} e_p \beta_u^p, \qquad u \geq 0,$$

where $\{\beta^p; p \in \mathbb{Z}\}$ is a sequence of independent Brownian motions. Hence, setting $\langle\cdot,\cdot\rangle$ for the inner product in $L^2(S)$, the matrix elements of $X^{\vee}$ are given by

$$\begin{aligned}
[X_{ts}^{\vee}]_{i,jk} &= \langle e_i, X_{ts}^{\vee}(\tilde{e}_j \otimes \tilde{e}_k)\rangle = \left\langle e_i; \int_s^t S_{tu}\, dx_u (S_{us}\tilde{e}_j)(S_{us}\tilde{e}_k)\right\rangle \\
&= \int_s^t \langle S_{tu} e_i; dx_u (S_{us}\tilde{e}_j)(S_{us}\tilde{e}_k)\rangle \\
&= \sum_{p\in\mathbb{Z}} \lambda_j^{-\eta}\lambda_k^{-\eta}\lambda_p^{-\nu/2} \int_s^t d\beta_u^p \langle S_{tu} e_i, e_p (S_{us} e_j)(S_{us} e_k)\rangle \\
&= \sum_{p; i=p+j+k} \lambda_j^{-\eta}\lambda_k^{-\eta}\lambda_p^{-\nu/2} \int_s^t e^{-\lambda_i(t-u)-\lambda_j(u-s)-\lambda_k(u-s)}\, d\beta_u^p.
\end{aligned}$$

Thus, the Hilbert–Schmidt norm of $X_{ts}^{\vee}$ in $\mathcal{L}(\mathcal{B}^2;\mathcal{B}_\eta)$ can be written as

$$\|X_{ts}^{\vee}\|^2_{\mathrm{HS},\mathcal{L}(\mathcal{B}_\eta^2;\mathcal{B}_\eta)} = \sum_{i,j,k\in\mathbb{Z}} \lambda_i^{2\eta} |[X_{ts}^{\vee}]_{i,jk}|^2.$$

From this simple computation, the following rules appear:

- Some multiple sums (involving terms of the form $\lambda_i^\alpha$) with constraints on the indices appear, due to the fact that $\{e_i; i \in \mathbb{Z}\}$ is the trigonometric basis of $L^2(S)$.

- Some contractions in the sums take place, because of the Brownian stochastic integrals.

With the above considerations in mind, we can associate to $X^{\vee}_{ts}$ the following graphical representation:

$$X^{\vee}_{ts} = \int_s^t d\beta^p_u$$

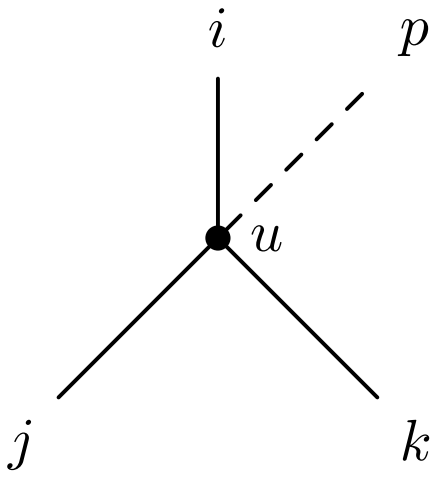

where the solid lines represent factors of $S$. This is a bookkeeping device for the relation between the various indices and time parameters. The computation of $E[\|X^{\vee}_{ts}\|^2_{\mathrm{HS},\mathcal{L}(\mathcal{B}^2_\eta;\mathcal{B}_\eta)}]$ corresponds to putting side by side two specular copies of this graph and connecting the corresponding top and bottom lines (to compute the HS norm), while contracting in all the allowed ways the dashed lines (to compute the contractions of the stochastic integrals). Doing so, we obtain the graph

$$E\|X^{\vee}_{ts}\|^2_{\mathrm{HS},\mathcal{L}(\mathcal{B}^2_\eta;\mathcal{B}_\eta)} = \int_s^t du$$

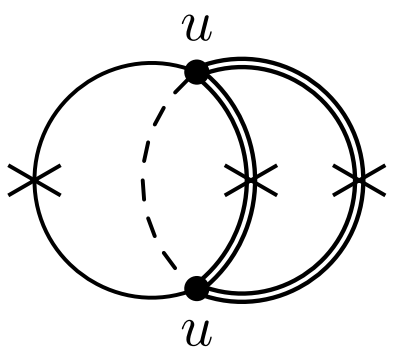

where we use the following convention: solid lines correspond to factors of $S$, time parameters are attached to vertices, crossed solid lines correspond to factors of $A^{2\eta}S$ (coming form contraction of output lines), crossed double lines correspond to factors of $A^{-2\eta}S$ (coming form contraction of input lines), dashed lines are associated to factors of $Q$ (coming from the Itô contraction of the noise). We fix an orientation for each edge of the graph and associate an index to each oriented edge. To each vertex corresponds a constraint that the sum of indexes of incoming edges minus indexes of outgoing edges should be zero. According to these rules, the formula of the mean squared norm is then

$$E[\|X^{\vee}_{ts}\|^2_{\mathrm{HS},\mathcal{L}(\mathcal{B}^2_\eta;\mathcal{B}_\eta)}] = \int_s^t du \sum_{i+j+k+l=0} \lambda_i^{2\eta}\lambda_j^{-\nu}\lambda_k^{-2\eta}\lambda_l^{-2\eta} e^{-2\lambda_i(t-u)-2\lambda_k(u-s)-2\lambda_l(u-s)},$$

and the reader can check that this is indeed the expression for the mean value of the HS norm of $X^{\vee}_{ts}$.

Consider now the expression

$$A = \sum_{i+j+k+l=0} \lambda_i^{2\eta}\lambda_j^{-\nu}\lambda_k^{-2\eta}\lambda_l^{-2\eta} e^{-2\lambda_i(t-u)-2\lambda_k(u-s)-2\lambda_l(u-s)}.$$

We trivially have

$$A \le \sum_{i+j+k+l=0} \lambda_i^{2\eta}\lambda_j^{-\nu}\lambda_k^{-2\eta}\lambda_l^{-2\eta} e^{-2\lambda_i(t-u)}.$$

Furthermore, setting $q = k + l$, applying Lemma 4.6, and recalling that we have chosen $\eta > 1/4$, we can bound $A$ as

$$A \lesssim \sum_{i+j+q=0} \lambda_i^{2\eta}\lambda_j^{-\nu}\lambda_q^{-2\eta} e^{-2\lambda_i(t-u)},$$

where we used the relation $\sum_{k+l=q} \lambda_k^{-2\eta}\lambda_l^{-2\eta} \lesssim \lambda_q^{-2\eta}$. Moreover, assuming that $\nu \le 2\eta$ we can use again Lemma 4.6 to get

$$A \lesssim \sum_i \lambda_i^{2\eta-\nu} e^{-2\lambda_i(t-u)} \lesssim (t-u)^{-a} \sum_i \lambda_i^{2\eta-\nu-a},$$

and choosing $a = 1/2 + 2\eta - \nu + \varepsilon$ so that the sum is convergent, we obtain $\lesssim (t-u)^{\Delta-3/2}$, with $\Delta = 1 - 2\eta + \nu - \varepsilon$. So we proved the graphical equation

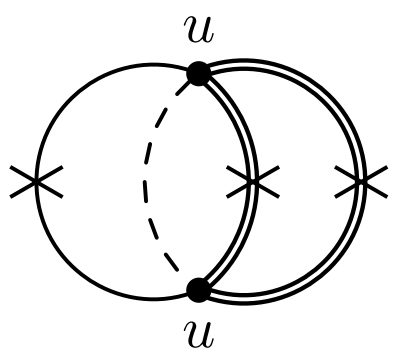

(123) $\lesssim (t-u)^{\Delta-3/2}$

and hence, if we suppose that $\Delta - 3/2 > -1$, that is, $\Delta > 1/2$, we have obtained that

$$E[\|X^{\vee}_{ts}\|^2_{\mathrm{HS},\mathcal{L}(\mathcal{B}^2_\eta;\mathcal{B}_\eta)}] \lesssim \int_s^t du\,(t-u)^{\Delta-3/2} \lesssim (t-s)^{\Delta-1/2},$$

which is the desired bound for Lemma 6.11. Let us say a few words about the condition $\Delta > 1/2$: if $\eta = 1/4 + \hat{\varepsilon}$, then one has $\Delta = 1/2 + \nu - \varepsilon + \hat{\varepsilon}$. This means that the condition $\Delta > 1/2$ can be met as soon as $\nu > 0$, which simply rules out the possibility of considering a space–time white noise at this stage.

(2) *Case of the operator* $X^{\curlyvee}$. With the same kind of considerations as for $X^{\vee}$, it can be shown that the the matrix elements of the operator $X^{\curlyvee}$ are

given by

$$[X^{\curlyvee}_{ts}]_{i,jkl} = \sum_{p+j+n=i} \sum_{q+k+l=n} \lambda_j^{-2\eta}\lambda_k^{-2\eta}\lambda_l^{-2\eta}\lambda_p^{-\nu/2}\lambda_q^{-\nu/2} \times \int_s^t e^{-\lambda_i(t-u)}\, d\beta_u^p\, e^{-\lambda_j(u-s)} \times \int_s^u e^{-\lambda_n(u-v)}\, d\beta_v^q\, e^{-\lambda_k(v-s)} e^{-\lambda_l(v-s)}.$$

Thus, its Hilbert–Schmidt norm in $\mathcal{L}(\mathcal{B}_\eta^3;\mathcal{B}_\eta)$ can be written as

$$\|X^{\curlyvee}_{ts}\|^2_{\mathrm{HS},\mathcal{L}(\mathcal{B}_\eta^3;\mathcal{B}_\eta)} = \sum_{i,j,k,l} \lambda_i^{2\eta} |[X^{\curlyvee}_{ts}]_{i,jkl}|^2,$$

and the following graphical representation can be associated to this last expression:

$$X^{\curlyvee}_{ts} = \int_s^t d\beta_u^p \int_s^u d\beta_v^q$$

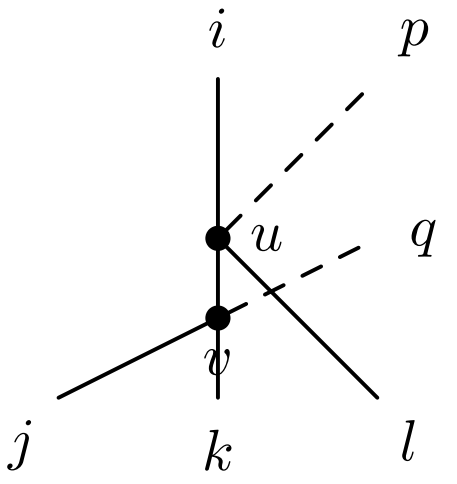

Thus, for the Hilbert–Schmidt norm of $X^{\curlyvee}_{ts}$, we obtain the graph

$$E[\|X^{\curlyvee}_{ts}\|^2_{\mathrm{HS},\mathcal{L}(\mathcal{B}_\eta^3;\mathcal{B}_\eta)}] = \int_s^t du \int_s^u dv$$

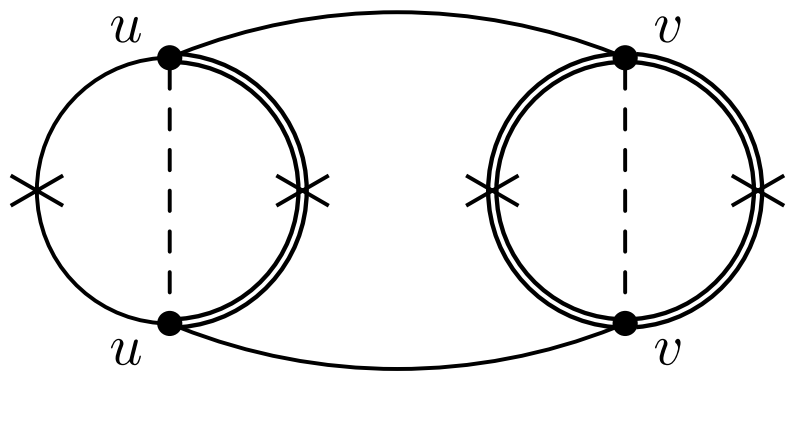

and the corresponding formula

$$E[\|X^{\curlyvee}_{ts}\|^2_{\mathrm{HS},\mathcal{L}(\mathcal{B}_\eta^3;\mathcal{B}_\eta)}] = \int_s^t du \int_s^u dv \sum_{i+j+k+l=0} \lambda_i^{2\eta}\lambda_j^{-\nu}\lambda_k^{-2\eta} \times \sum_{n+m+o=l} \lambda_m^{-\nu}\lambda_n^{-2\eta}\lambda_o^{-2\eta} e^{-2\lambda_i(t-u)-2\lambda_l(u-v)-2\lambda_n(v-s)-2\lambda_o(v-s)-2\lambda_k(u-s)}.$$

The strategy to control this expression is now straightforward: bounding the exponential and performing the time integrations gives, for two positive constants $a, b$ such that $a+b<2$,

$$E[\|X^{\curlyvee}_{ts}\|^2_{\mathrm{HS},\mathcal{L}(\mathcal{B}^3_\eta;\mathcal{B}_\eta)}] \lesssim (t-s)^{2-(a+b)} \sum_{i+j+k+l=0} \lambda_i^{2\eta-a}\lambda_j^{-\nu}\lambda_k^{-2\eta}\lambda_l^{-b} \sum_{n+m+o=l} \lambda_m^{-\nu}\lambda_n^{-2\eta}\lambda_o^{-2\eta}.$$

Using the fact that $2\eta > 1/2$, we can reduce this to the bound

$$E[\|X^{\curlyvee}_{ts}\|^2_{\mathrm{HS},\mathcal{L}(\mathcal{B}^3_\eta;\mathcal{B}_\eta)}] \lesssim (t-s)^{2-(a+b)} \sum_{i+j+l+k=0} \lambda_i^{2\eta-a}\lambda_j^{-\nu}\lambda_k^{-2\eta}\lambda_l^{-b} \sum_{n+m=l} \lambda_m^{-\nu}\lambda_n^{-2\eta}.$$

Assuming moreover that $\nu \le 2\eta$, we get

$$E[\|X^{\curlyvee}_{ts}\|^2_{\mathrm{HS},\mathcal{L}(\mathcal{B}^3_\eta;\mathcal{B}_\eta)}] \lesssim (t-s)^{2-(a+b)} \sum_{i+j+l+k=0} \lambda_i^{2\eta-a}\lambda_j^{-\nu}\lambda_k^{-2\eta}\lambda_l^{-\nu-b}.$$

At this point, choose $b = 2\eta - \nu$ so that $\nu + b = 2\eta$ and

$$E[\|X^{\curlyvee}_{ts}\|^2_{\mathrm{HS},\mathcal{L}(\mathcal{B}^3_\eta;\mathcal{B}_\eta)}] \lesssim (t-s)^{2-(a+b)} \sum_{i+k+l=0} \lambda_i^{2\eta-a}\lambda_l^{2\eta}\lambda_k^{-2\eta}.$$

Hence, this sum is finite if we choose $a = 1/2 + 2\eta - \nu + \varepsilon$, and we get

$$E[\|X^{\curlyvee}_{ts}\|^2_{\mathrm{HS},\mathcal{L}(\mathcal{B}^3_\eta;\mathcal{B}_\eta)}] \lesssim (t-s)^{2\Delta-1/2}.$$

Before proceeding to the estimation of the other more complex operators, let us make a useful observation. Consider the following subgraph on the left of the previous graph:

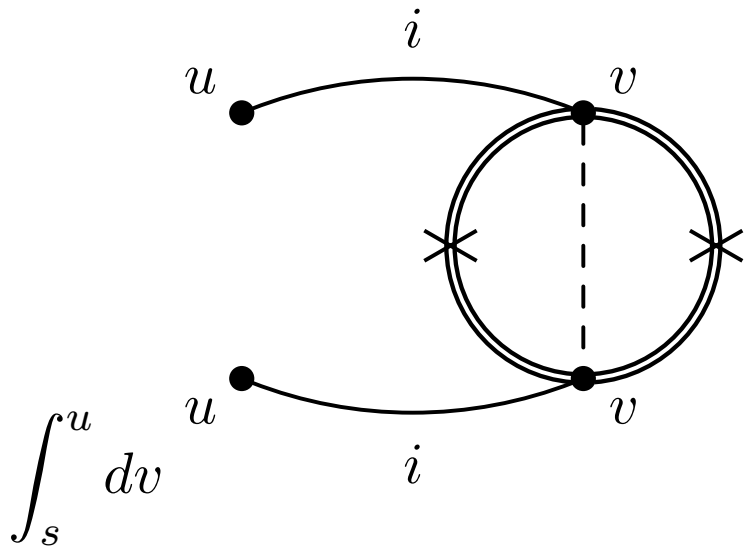

After reduction of the crossed double lines (carrying the factors due to $A^{-\eta}$) by an iterated application of Lemma 4.6, we obtain the following expression which corresponds to a bound for this graph:

$$\lesssim \int_s^u dv\, \lambda_i^{-\nu} e^{-2\lambda_i(u-v)} \lesssim \int_s^u \frac{dv}{(u-v)^b} \lambda_i^{-\nu-b},$$

so that choosing $b = 2\eta - \nu$, we get an estimate of the form $\lesssim (v-s)^{\Delta} \lambda_i^{-2\eta}$. Summarizing these considerations, we have obtained the graphical equation

(124) $$\int_s^u dv$$ 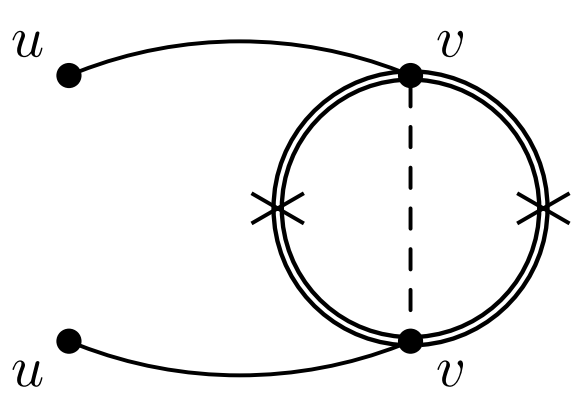

$$\lesssim (u-s)^{\Delta}$$ 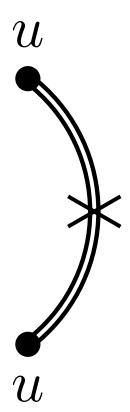

which we will use multiple times below.

6.6. *More complex graphs.* The tools we have introduced so far will allow us to treat the two remaining cases we are left with, namely $X^{\vee\vee\vee}$ and $X^{\vee\cdot\cdot\vee}$.

(1) *Case of the operator* $X^{\vee\vee\vee}$. By using the same kind of arguments as in the previous subsection, we obtain a representation of the form

$$X_{ts}^{\vee\vee\vee} = \int_s^t d\beta_u^p \int_s^u d\beta_v^q \int_s^v d\beta_w^r$$ 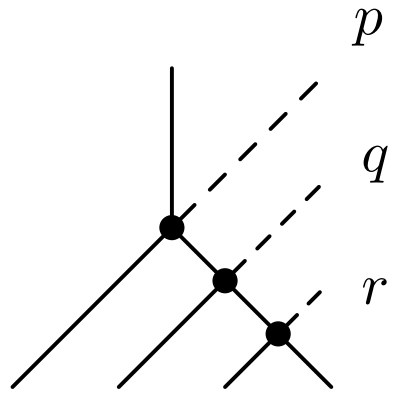

Thus, for the computation of $E[\|X_{ts}^{\vee\vee\vee}\|^2_{\mathrm{HS},\mathcal{L}(\mathcal{B}_\eta^4;\mathcal{B}_\eta)}]$, we obtain the graph

$$E\|X_{ts}^{\vee\vee\vee}\|^2_{\mathrm{HS},\mathcal{L}(\mathcal{B}_\eta^3;\mathcal{B}_\eta)}$$

$$= \int_s^t du \int_s^u dv \int_s^u dw$$ 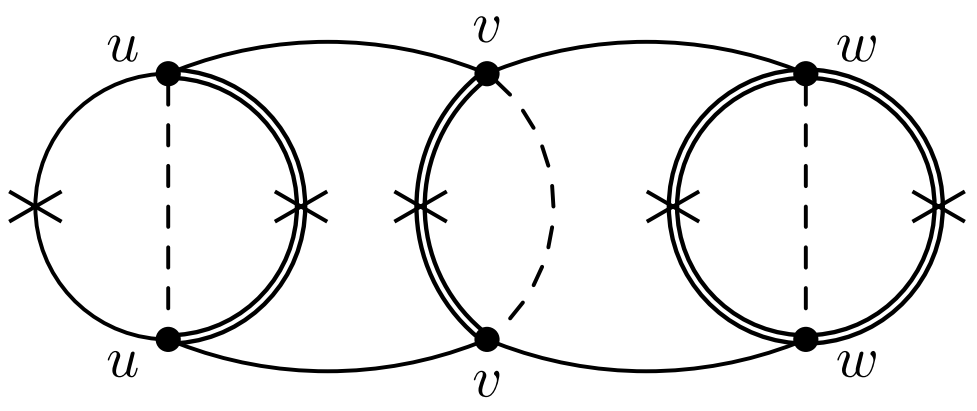

Now, invoking repeatedly relation (124), we can iteratively reduce the above graph to obtain a bound for $E[\|X_{ts}^{\vee\vee\vee}\|^2_{\mathrm{HS},\mathcal{L}(\mathcal{B}_\eta^4;\mathcal{B}_\eta)}]$ of the form

$$\lesssim \int_s^t du (u-s)^{2\Delta}$$ 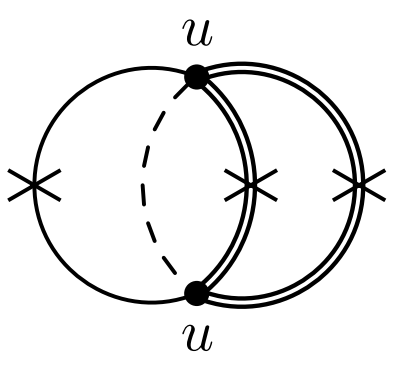

$$\lesssim \int_s^t du\,(t-u)^{\Delta-3/2}(u-s)^{2\Delta} \lesssim (t-s)^{3\Delta-1/2},$$

which is again what is needed for our Lemma 6.11.

(2) *Case of the operator* $X^{\vee\cdot\cdot\vee}$. Using the same conventions as before, the operator $X^{\vee\cdot\cdot\vee}$ can be represented as:

$$X_{ts}^{\vee\cdot\cdot\vee} = \int_s^t d\beta_u^p \int_s^u d\beta_v^q \int_s^u d\beta_w^r$$

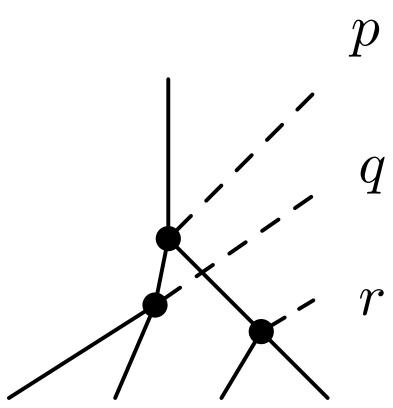

Now our current situation is slightly different from the previous ones, since in the triple Brownian integral above, the last two are not iterated. This means that we have to handle some sums of the form $E[(\sum_\alpha Y_\alpha Z_\alpha)^2]$ for some centered Gaussian random variables $(Y_\alpha, Z_\alpha)_\alpha$ forming a Gaussian vector. The standard way to compute such sums is to write

$$\begin{aligned} E\Big[\Big(\sum_\alpha Y_\alpha Z_\alpha\Big)^2\Big] &= \sum_{\alpha,\beta} E[Y_\alpha Z_\alpha Y_\beta Z_\beta] \\ &= \sum_{\alpha,\beta} E[Y_\alpha Z_\alpha]E[Y_\beta Z_\beta] + E[Y_\alpha Y_\beta]E[Z_\alpha Z_\beta] \\ &\quad + E[Y_\alpha Z_\beta]E[Y_\beta Z_\alpha]. \end{aligned}$$

By extrapolating this elementary consideration to our situation, this implies that, in the computation of $E\|X_{ts}^{\vee\cdot\cdot\vee}\|^2_{\mathrm{HS},\mathcal{L}(\mathcal{B}_\eta^4;\mathcal{B}_\eta)}$, three different kind of contractions are involved. Hence, we also obtain three different graphs:

$$E\|X_{ts}^{\vee\cdot\cdot\vee}\|^2_{\mathrm{HS},\mathcal{L}(\mathcal{B}_\eta^3;\mathcal{B}_\eta)}$$

$$= \int_s^t du \int_s^u dv \int_s^u dw \qquad (125)$$

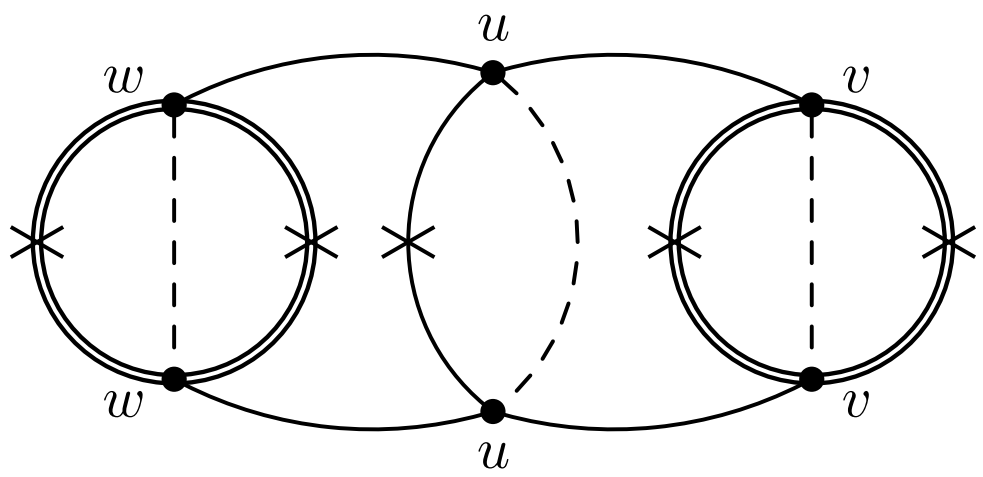

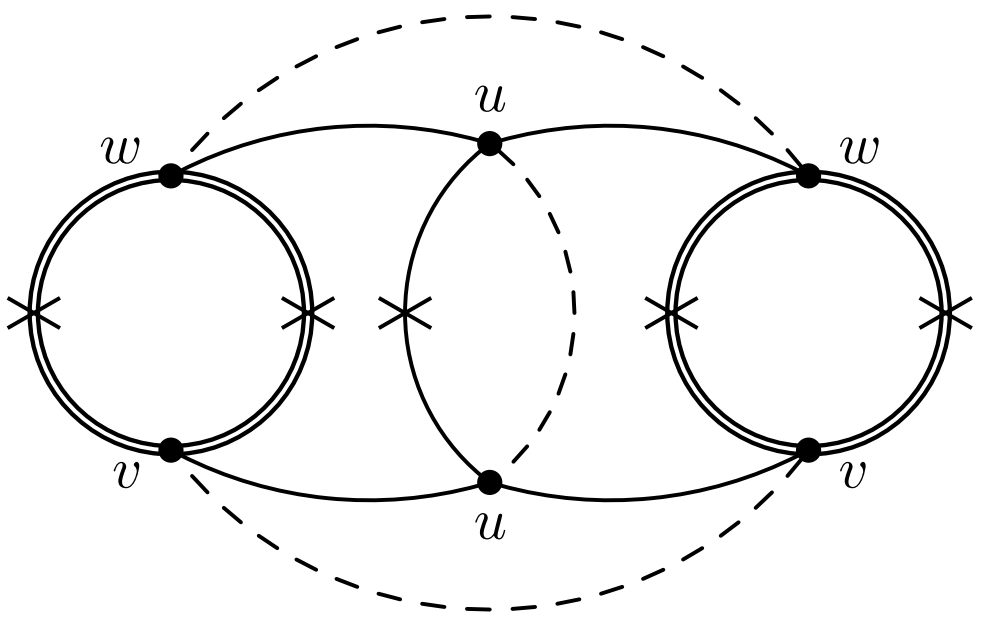

$$+\int_s^t du \int_s^u dv \int_s^u dw$$

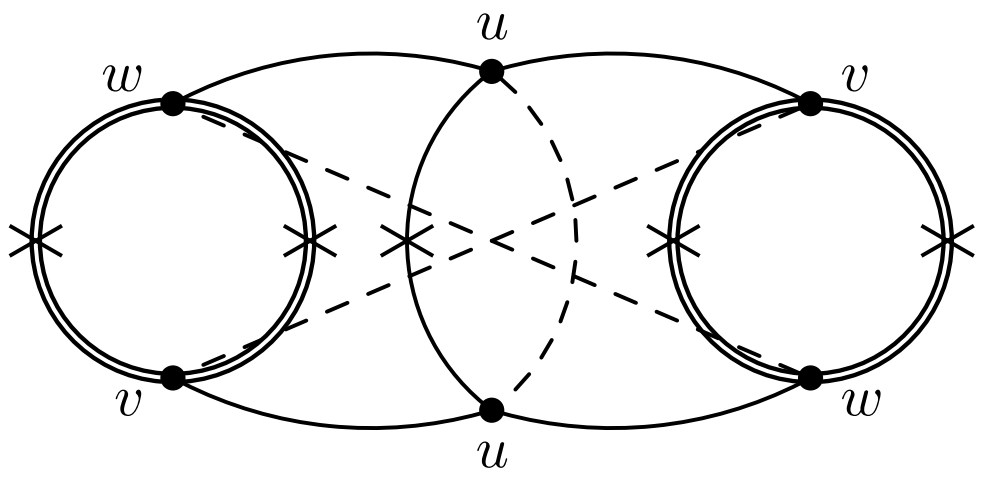

$$+\int_s^t du \int_s^u dv \int_s^u dw$$

Observe that the first of those graphs already appeared in the study of $X^{\vee}$. We will then focus on the two other ones.

The analysis of the second graph above can be started by reducing again the crossed double lines, which gives a new graph of the form

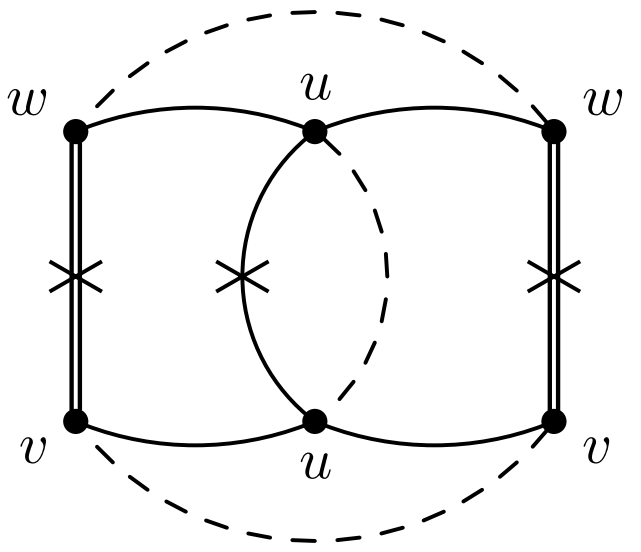

$$\int_s^t du \int_s^u dv \int_s^u dw \tag{126}$$

However, at this point, we cannot proceed as in the previous cases, with a sequence of reduction of subgraphs, in order to prove the convergence. Indeed, in the current situation, some irreducible triangular structures like the following appear, at the top and bottom of the graph (126):

(127)

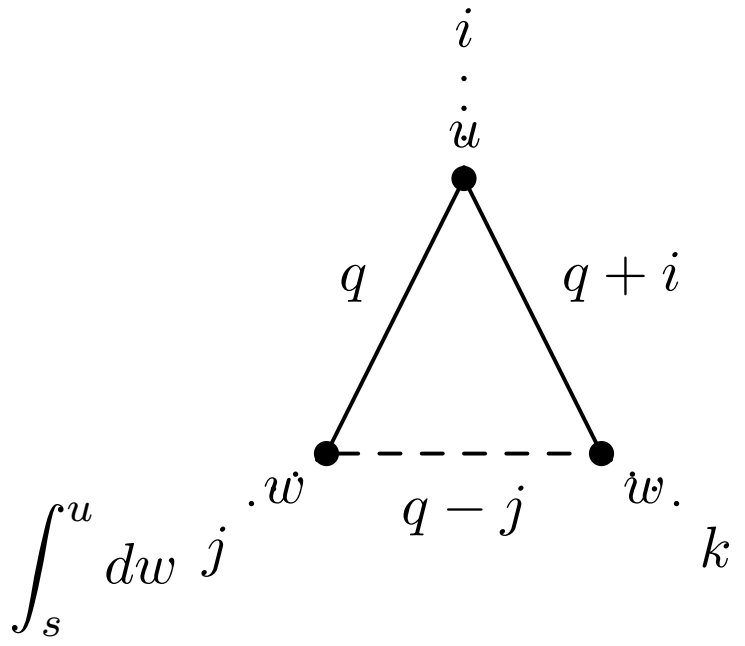

where the dotted lines stand for the remaining part of the graph, and where we put explicit indexes on outgoing lines and on the edges of the subgraph. Note that by the rules we have imposed on the graph, the constraint $i + j + k = 0$ (we consider the dotted lines directed inwards) holds true. The contribution of the triangular structure is thus given by

$$\delta_{i+j+k=0} \int_s^u dw \sum_q \lambda_{q+k}^{-\nu} e^{-\lambda_q (u-w) - \lambda_{q-i}(u-w)},$$

and we can bound this last expression by

$$\lesssim \delta_{i+j+k=0} \int_s^u \frac{dw}{(u-w)^b} \sum_q \lambda_{q+k}^{-\nu} \lambda_q^{-b/2} \lambda_{q-i}^{-b/2},$$

for some $b \in (0,1)$. The latter sum is finite when $\nu + b > 1/2$, which means that, by choosing $b = 1/2 - \nu + \varepsilon$, the triangular structure yields a bound of order

$$\lesssim \delta_{i+j+k=0} (u-s)^{1-b}.$$

Summarizing the previous discussion, it is now easily seen that the structure (127) behaves like a simple vertex, up to an appropriate factor of $(u - s)$:

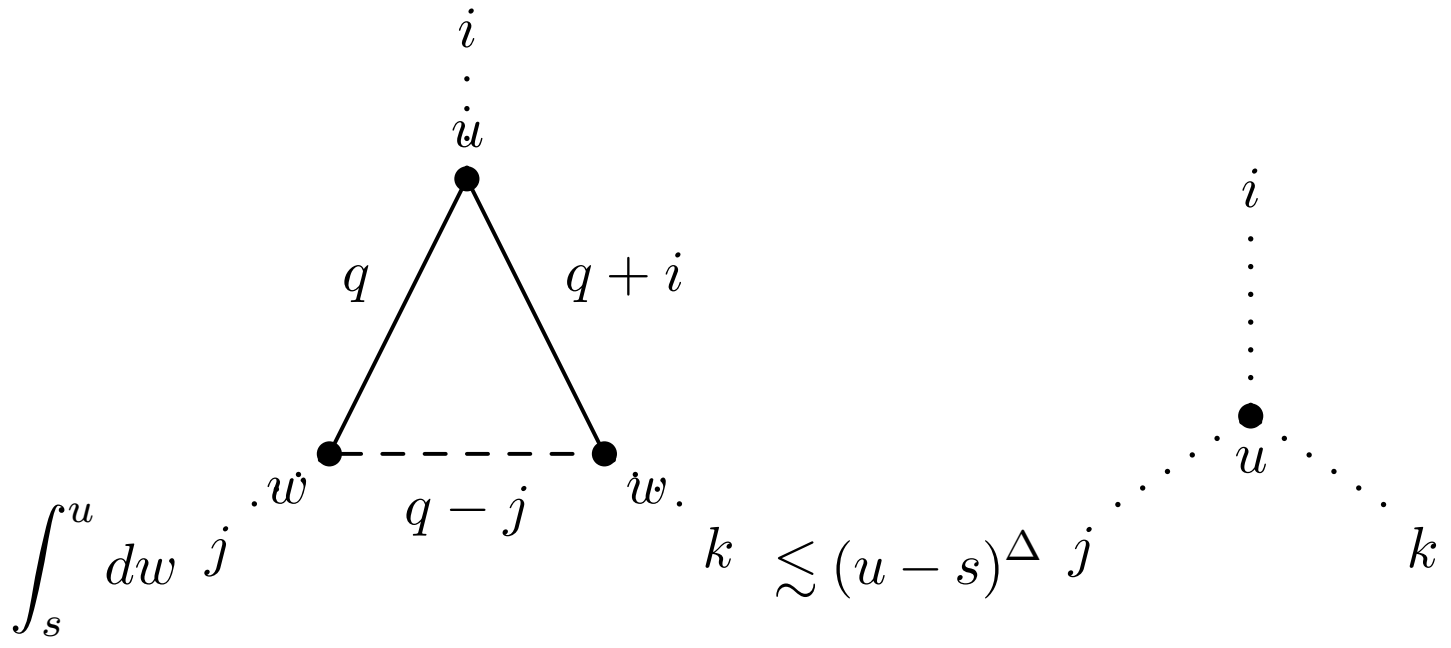

Using this fact, we can reduce our graph (126) to the following simpler structure:

$$\lesssim \int_s^t du (u-s)^{2\Delta}$$

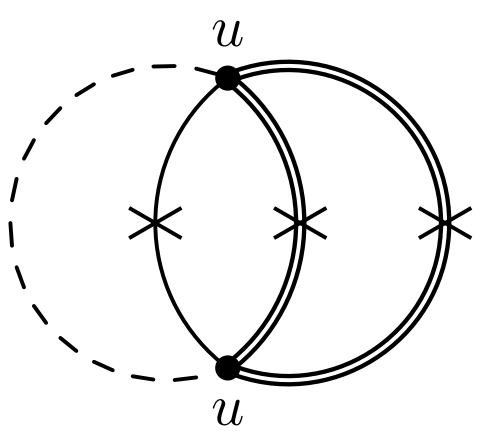

$$\lesssim (t-s)^{3\Delta-1/2},$$

where the last bound has been obtained similarly to (123).

Finally, let us associate a bound to the third graph in (125). First, notice that this third graph can be reduced to

(128) $$\int_s^t du \int_s^u dv \int_s^u dw$$

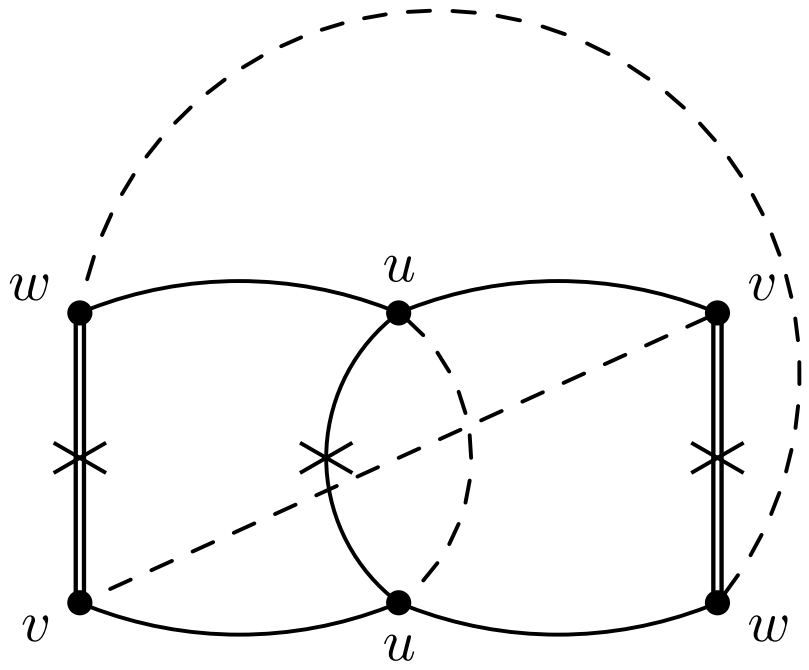

Furthermore, this graph contains the irreducible subgraph

$$\int_s^u dv \int_s^u dw$$

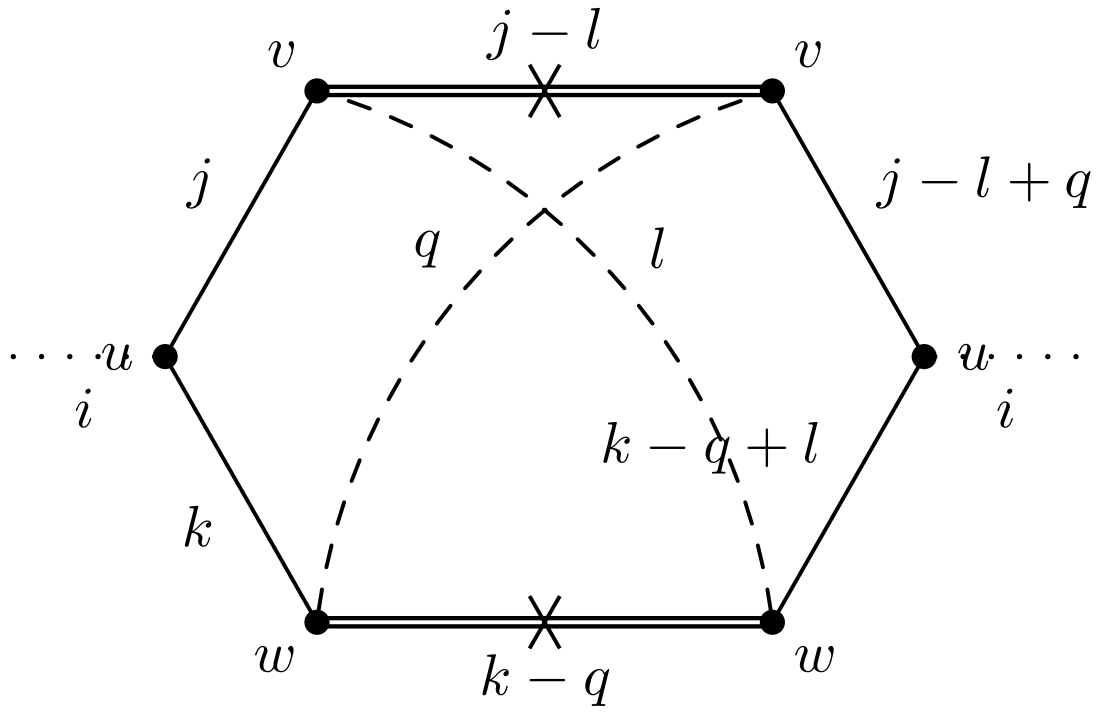

which corresponds to the expression

$$D \equiv \int_s^u dv \int_s^u dw$$
$$\times \sum_{j+k=i} \sum_{q,l} \exp\{-\lambda_j(u-v) - \lambda_k(u-w)$$

$$
\begin{aligned}
&\qquad - \lambda_{j-l+q}(u-v) - \lambda_{k-q+l}(u-w) \\
&\qquad - \lambda_{j-l}(w+v-2s) - \lambda_{k-q}(w+v-2s)\} \\
&/(\lambda_l^\nu \lambda_q^\nu \lambda_{j-l}^{2\eta} \lambda_{k-q}^{2\eta}).
\end{aligned}
$$

Introducing an additional parameter $b$ and bounding the exponential terms as usual gives

$$
D \lesssim \int_s^u \frac{dv}{(u-v)^b} \int_s^u \frac{dw}{(u-w)^b} \sum_{j+k=i} \sum_{q,l} \frac{1}{\lambda_j^b \lambda_k^b \lambda_l^\nu \lambda_q^\nu \lambda_{j-l}^{2\eta} \lambda_{k-q}^{2\eta}}.
$$

Now, using Lemma 4.6, we can bound the sums over $q$ and $l$ in order to obtain

$$
D \lesssim \int_s^u \frac{dv}{(u-v)^b} \int_s^u \frac{dw}{(u-w)^b} \sum_{j+k=i} \frac{1}{\lambda_j^{b+\nu} \lambda_k^{b+\nu}}.
$$

Thus, choosing $b = 2\eta - \nu$, we end up with

$$
D \lesssim (u-s)^{2\Delta} \lambda_i^{-2\eta},
$$

which means that we have obtained the graphical inequality

$$\int_s^u dv \int_s^u dw$$ 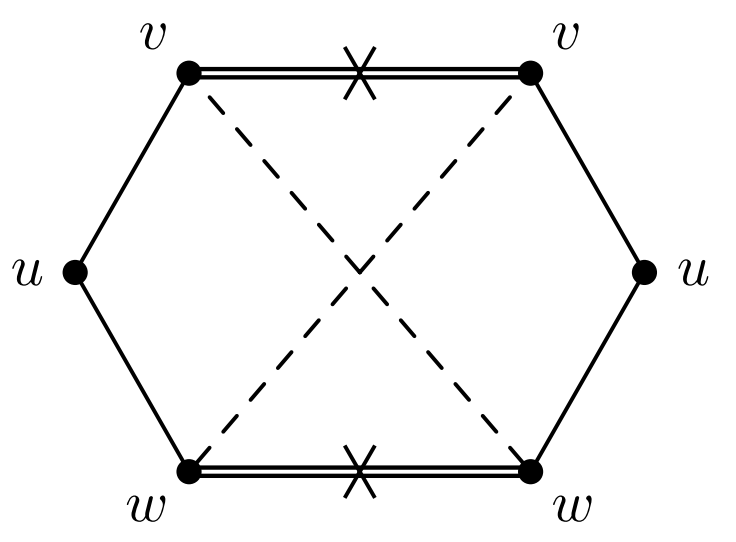 $$\lesssim (u-s)^{2\Delta}$$ 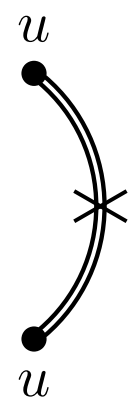

Plugging this representation in the complete graph (128), we obtain

$$\lesssim \int_s^t du (u-s)^{2\Delta}$$ 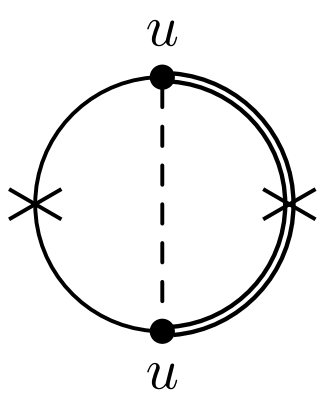 $$\lesssim (t-s)^{3\Delta - 1/2}.$$

Going back to Lemma 6.11, we should still treat the case of $X^\tau$ for $\tau =$ 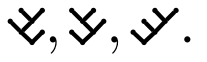. But these estimates are now mere variations of the previous ones, and are left to the reader for sake of conciseness.

CEREMADE
Université de Paris-Dauphine
Place du Maréchal De Lattre De Tassigny
75775 Paris cedex 16
France
E-mail: massimiliano.gubinelli@ceremade.dauphine.fr

Institut Élie Cartan Nancy
Université de Nancy
B.P. 239
54506 Vandœuvre-lès-Nancy Cedex
France
E-mail: tindel@iecn.u-nancy.fr